\documentclass[11pt]{amsart}
\usepackage{graphicx}
\usepackage{latexsym,amsmath,amsfonts,amscd, amsthm}
\usepackage{color}
\usepackage[colorlinks=true]{hyperref}
\hypersetup{urlcolor=blue, citecolor=red}
\usepackage{tikz}
\usetikzlibrary{arrows,backgrounds,snakes}
\usepackage{epsfig}
\usepackage{dsfont}
\usepackage{amssymb}
\usepackage{xcolor}
\usepackage{comment}
\usepackage{subfigure}

\newtheoremstyle{plainNoItalics}{}{}{\normalfont}{}{\bfseries}{.}{ }{}
\theoremstyle{plain}
\newtheorem{thm}{Theorem}[section]

\newtheorem{defn}[thm]{Definition}

\newtheorem{prop}[thm]{Proposition}

\renewcommand{\theequation}{\thesection.\arabic{equation}}

\newcommand{\beq}{\begin{equation}}
\newcommand{\eeq}{\end{equation}}
\newcommand{\beqa}{\begin{eqnarray}}
\newcommand{\eeqa}{\end{eqnarray}}
\newcommand{\bit}{\begin{itemize}}
\newcommand{\eit}{\end{itemize}}
\newcommand{\bedef}{\begin{defn}}
\newcommand{\edefn}{\end{defn}}
\newcommand{\bpro}{\begin{prop}}
\newcommand{\epro}{\end{prop}}

\newcommand{\mZ}{\mathcal{Z}}

\newcommand{\eps}{\varepsilon}

\newcommand{\mQ}{{\mathcal Q}}

\newcommand{\xL}{{x_{i-\frac{1}{2}}}}
\newcommand{\xR}{{x_{i+\frac{1}{2}}}}
\newcommand{\yL}{{y_{j-\frac{1}{2}}}}
\newcommand{\yR}{{y_{j+\frac{1}{2}}}}
\newcommand{\iL}{{i-\frac{1}{2}}}
\newcommand{\iR}{{i+\frac{1}{2}}}
\newcommand{\jL}{{j-\frac{1}{2}}}
\newcommand{\jR}{{j+\frac{1}{2}}}
\newcommand{\pt}{\partial_t}
\newcommand{\px}{\partial_x}
\newcommand{\py}{\partial_y}
\newcommand{\nablax}{\nabla_{\bf x}}

\newcommand{\bfx}{{\bf x}}
\newcommand{\bfv}{{\bf v}}

\newcommand{\bfg}{{\bf g}}
\newcommand{\bfh}{{\bf h}}
\newcommand{\bfu}{{\bf u}}
\newcommand{\bfn}{{\bf n}}

\newcommand{\bfA}{{\bf A}}
\newcommand{\bfB}{{\bf B}}
\newcommand{\bfD}{{\bf D}}
\newcommand{\bfF}{{\bf F}}

\newcommand{\bfM}{\mathcal{\bf M}}
\newcommand{\bfK}{{\bf K}}
\newcommand{\bfR}{{\bf R}}
\newcommand{\bfS}{{\bf S}}
\newcommand{\bfU}{{\bf U}}
\newcommand{\bfV}{{\bf V}}
\newcommand{\bfI}{{\bf I}}
\newcommand{\bfZ}{{\bf Z}}

\newcommand\ds{ \displaystyle }

\email{francis.filbet@math.univ-toulouse.fr}
\email{txiong@xmu.edu.cn}

\title[Hybrid discontinuous Galerkin for multiscale kinetic equations]{A Hybrid Discontinuous Galerkin Scheme for Multi-scale Kinetic Equations} 

\author[Francis Filbet and Tao Xiong]{}



\begin{document}
\maketitle

\centerline{\scshape Francis Filbet}
\medskip
{\footnotesize
    % please put the address of the author
    \centerline{Institut de Math\'ematiques de Toulouse, UMR5219
      Universit\'e de Toulouse; CNRS \& IUF}
    \centerline{UPS, F-31400, Toulouse, France}
}

\medskip
\centerline{\scshape Tao Xiong}
\medskip
{\footnotesize
   % please put the address of the author
\centerline{School of Mathematical Sciences, Fujian Provincial Key Laboratory of Mathematical Modeling} 
\centerline{ and High-Performance Scientific Computing, Xiamen University}
\centerline{Xiamen, Fujian, 361005, P.R. China}
}

\bigskip

\begin{abstract}
We develop a multi-dimensional hybrid discontinuous Galerkin method
for multi-scale kinetic equations. This method is  based on moment
realizability matrices, a concept  introduced by D. Levermore,
        W. Morokoff and B. Nadiga in \cite{levermore1998moment} and used in
\cite{filbet2015hierarchy} for one dimensional problem. The main issue
addressed in this paper is to  provide a simple indicator to select the
most appropriate model and to apply a compact numerical scheme to reduce
the interface region between different models. We also construct a
numerical flux for the fluid  model obtained as the asymptotic limit of the
flux of the kinetic equation. Finally we perform several numerical simulations for time evolution and stationary problems.   
\end{abstract}

\vspace{0.1cm}

\tableofcontents

\section{Introduction}
\label{sec1}
\setcounter{equation}{0}
\setcounter{figure}{0}
\setcounter{table}{0}

Many physical  problems, such as micro-electro-mechanical systems,  involve boundary layers or transitional
          regimes which cannot be solved by using standard fluid
          models. Hence a kinetic model is needed to accurately
          describe complex phenomena occurring around the boundary. However, a kinetic description is computationally very expensive to simulate and we desire to use it only locally in space. An interesting approach is to design hybrid kinetic/fluid schemes, with domain decomposition which can automatically and accurately identify fluid and kinetic zones. In this case we only solve the kinetic model in a minimized area around the kinetic boundary layers or shocks, while take advantage of the low computational cost of numerical methods for the fluid system elsewhere.

Taking into account of collisions, the multi-scale kinetic equation for the description of particle dynamics in a dilute gas is given by
\beq
\label{eq:bgk}
\left\{
	\begin{array}{l}
\ds	\frac{\partial f}{\partial t} + \bfv \cdot \nablax f  \,= \, \frac{1}{\eps}\,\mathcal{Q}(f) , 
	\\ \, \\
	f(0,\bfx,\bfv)  \,= \, f_0(\bfx,\bfv),
	\end{array}
\right.
\eeq
with $\bfx \in \Omega \subset \mathbb{R}^{d_x}, \bfv \in
\mathbb{R}^3$. The open set $\Omega$ is a bounded Lipschitz-continuous
domain in $\mathbb{R}^{d_x}$ and supplied with some boundary
conditions, whereas the particle distribution function $f:=f(t,\bfx,
\bfv)$ depends on time and phase space variables $(\bfx,\bfv) \in
\Omega \times\mathbb{R}^3$ and the initial datum $f_0$ is a non-negative function. The parameter $\eps > 0$ is the dimensionless Knudsen number, which is defined as the ratio of the mean free path of particles over a typical length scale such as the size of the spatial domain. It measures the rarefaction of the gas, namely, the gas is in a rarefied or kinetic regime if $\eps \sim 1$ and in a dense or fluid regime if $\eps \ll 1$. 
The collision operator $\mQ(f)$ may be given by the full Boltzmann operator \cite{cercignani1969mathematical,cercignani1988boltzmann}
\beq
\mQ(f)(\bfv) \,=\,\int_{\mathbb{R}^3}\int_{\mathbb{S}^2}B(|\bfv-\bfv_\star|,\cos\theta)(f'_\star\,f'-f_\star\,f)\,d{\bf \sigma}d\bfv_\star,
\label{boltzoper}
\eeq
where we used the shorthand $f = f(v)$, $f_* = f(v_*)$, $f ^{'} = f(v')$, $f_* ^{'} = f(v_* ^{'})$.   
			The velocities of the colliding pairs $(v,v_*)$ and $(v',v'_*)$ are related by
			\begin{equation*}
			  \bfv' \,=\, \frac{\bfv+\bfv_*}{2} \,+\, \frac{|\bfv-\bfv_*|}{2} \,\sigma, \qquad \bfv'_* \,=\, \frac{\bfv+\bfv^*}{2} \,-\, \frac{|\bfv-\bfv_*|}{2} \,\sigma.
			\end{equation*}
			The \emph{collision kernel} $B$ is a
                        non-negative function which by physical
                        arguments of invariance only depends on
                        $|\bfv-\bfv_*|$ and	$\cos \theta = {\widehat \bfu}
                        \cdot \sigma$, where ${\widehat \bfu}
                        =(\bfv-\bfv_*)/|\bfv-\bfv_*|$. A simpler choice is the BGK
                        operator \cite{bhatnagar1954model} given by
\beq
\mQ(f)(\bfv) \,=\,\nu(\rho,T)\big(\mathcal{M}[f]-f),
\label{bgkoper}
\eeq
or a modified ES-BGK operator \cite{andries2000gaussian}
\beq
\mQ(f)(\bfv) \,=\,\nu(\rho,T)\big(\mathcal{G}[f]-f),
\label{esbgkoper}
\eeq
where $\nu$ is a collision frequency. In the ES-BGK operator, $\mathcal{G}(f)$ is a Gaussian defined as
\[
\mathcal{G}[f]=\frac{\rho}{\sqrt{{\rm det}(2\pi \mathcal{T})}}\exp\left(-\frac{(\bfv-\bfu)\mathcal{T}^{-1}(\bfv-\bfu)}{2}\right),
\]
with a corrected stress tensor
\[
\mathcal{T}\,:=\, (1-\beta) T \, {\bf I} + \beta \, \Theta, \quad \rho\, \Theta = \int_{\mathbb{R}^3}(\bfv-\bfu)\otimes(\bfv-\bfu)\, f\, d\bfv,
\]
where $\bfI$ is the identity matrix. The parameter $ \beta\in(-\infty, 1)$ is used to modify the value of the Prandtl number through the formula
\[
0< {\rm Pr} \,=\,\frac{1}{1-\beta} \le +\infty.
\]
The correct Prandtl number for a monoatomic gas of hard spheres is equal to $2/3$ for $\beta=-1/2$, whereas
the classical BGK operator has a Prandtl number equal to $1$ for $\beta =0$. We refer to \cite{filbet2015hierarchy} for more discussions about these collision operators. These Boltzmann-like collision operators share the fundamental properties of conserving mass, momentum and energy, that is
\beq
\int_{\mathbb{R}^3}\mQ(f)(\bfv) \, m(\bfv) \,d\bfv \,=\, {\bf 0}^T_{\mathbb{R}^5}, 
\label{eq:cons}
\eeq
where $m(\bfv) = (1,\bfv,\frac12|\bfv|^2)^T$ and super index $T$ denotes the transpose of the corresponding vector. The equilibrium of the collision operator, when $\mQ(f)=0$, is given by the local Maxwellian distribution function $\mathcal{M}[f]$, which is defined as :
$$
\mathcal{M}[f] \, := \, \frac{\rho}{(2\pi T)^{3/2}} \exp\Big(-\frac{|\bfv-\bfu|^2}{2T}\Big).
$$
The density $\rho$, mean velocity $\bfu$ and temperature $T$ are macroscopic moments of the distribution function $f$, which can be computed as
$$
\rho = \int_{\mathbb{R}^3} \, f(\bfv) \, d \bfv, \quad \bfu = \frac{1}{\rho} \int_{\mathbb{R}^3} \bfv f(\bfv) \, d \bfv, \quad T =  \frac{1}{3\rho} \int_{\mathbb{R}^3} |\bfv-\bfu|^2 f(\bfv) \, d \bfv.
$$

There are already several works about hybrid methods for multi-scale kinetic equations in the literature, which mostly rely on the same domain decomposition technique.  I. Boyd, G. Chen and G. Candler \cite{boyd1995predicting} first proposed a macroscopic criterion based on the local Knudsen number, they pass from a hydrodynamic description to a kinetic one when the quantity is below a (problem-dependent) threshold. It was practically used by V. Kolobov et. al. \cite{kolobov2007unified} with a discrete
velocity model of direct numerical solution (DNS) for the Boltzmann
equation and a gas-kinetic scheme for the hydrodynamic part, and
recently by P. Degond and G. Dimarco \cite{degond2012fluid} 
with a Monte Carlo solver for the kinetic equation and a finite volume method for the macroscopic ones. Another hydrodynamic breakdown indicator based on the viscous and
heat flux of the Navier-Stokes equations is introduced by S. Tiwari \cite{tiwari1998coupling}
and has been used with various deterministic kinetic and hydrodynamic solvers \cite{degond2010multiscale,tiwari2009particle,tiwari2014simulations,alaia2012hybrid}.
Recently, F. Filbet and T. Rey \cite{filbet2015hierarchy} proposed a domain
decomposition indicator  based on moment realizability matrices, a
concept first introduced by D. Levermore, W. Morokoff and B. Nadiga
\cite{levermore1998moment}. They have shown that criteria from/to hydrodynamic
to/from kinetic are both needed, and finite volume schemes for
both the Boltzmann equation and the hydrodynamic equations are
explored. This criterion has lately been used by T. Xiong and J.-M. Qiu
\cite{xiong2017hierarchical} to form a hierarchy high order
discontinuous Galerkin schemes for the BGK equation  under a
micro-macro decomposition framework. 
%%%%% FJF :  Ajouter une explication sur ce critere

However, this criterion requires some computational
effort since it is based on the comparison of local eigenvalues on
each cell, hence the cost is not negligible in high
dimension. Furthermore, high order finite volume schemes often need  
wide stencils. This may easily be done in 1D in space
\cite{filbet2015hierarchy}, however, it is not very convenient or
robust to be generalized to multi-dimensional problems, as it would require each region as wide as the stencil along each space dimension. It will cause some trouble when switching from
one region to the other dynamically with time evolution.

The aim of this paper is to develop a hybrid discontinuous Galerkin scheme, which
extends the work in \cite{filbet2015hierarchy} in high dimension. The
well recognized discontinuous Galerkin method
\cite{karniadakis2000discontinuous}, due to its advantages of
compactness, $h$-$p$ adaptivity, high efficiency on parallelization,
and flexibility on complicated geometries, has been widely applied to
physical and engineering problems. For compressible Euler or
Navier-Stokes systems, F. Bassi and S. Rebay first introduced a discontinuous
Galerkin method with macroscopic numerical flux such as Godunov flux
for the advection and a primal formulation with centered numerical
flux for the diffusion
\cite{bassi1997high1,bassi1997high,bassi2002numerical}. It has
provided a good framework for solving the compressible Navier-Stokes
system and many works are followed, {\it e.g.} a positivity preserving
high order discontinuous Galerkin scheme recently proposed by X. Zhang
\cite{zhang2017positivity}. Some other type of discontinuous Galerkin
schemes can be referred to
\cite{lomtev1999discontinuous,luo2010reconstructed,cheng2016hybrid}
and many references therein. For the kinetic Boltzmann, ellipsoidal
statistical BGK  equations, A. Alexeenko {\it et al.}
\cite{alexeenko2008high} proposed a high order explicit Runge-Kutta
discontinuous Galerkin method, with Newton's iteration solving the
nonlinear collisional source term. It has recently been generalized to
2D in space by W. Su {\it et al.} in \cite{su2017stable}. On the other
hand, to take care of the stiff collisional term and avoid a stringent
time step size, T. Xiong {\it et al.} \cite{xiong2015high} proposed a
high order asymptotic preserving nodal discontinuous Galerkin method
coupled with an implicit-explicit scheme for the BGK equation. There
are also some other works which would like to build a connection from
the kinetic equation to the hydrodynamic system, and discontinuous
Galerkin methods for compressible Navier-Stokes equations with kinetic flux-vector splitting are proposed, such as the kinetic flux-vector splitting flux following Chapman-Enskog velocity distribution function \cite{chou1997kinetic,chandrashekar2013discontinuous}, or gas-kinetic schemes
mimicking a Hilbert expansion \cite{liu2007runge,ren2015multi}, etc.

In this paper, we propose a simplified indicator, based on
\cite{filbet2015hierarchy}, to determine the kinetic and fluid regions and develop a discontinuous Galerkin method,
which is well suited to couple the kinetic and fluid model at the interface. Due to the compactness of the discontinuous Galerkin method, the interface coupling
condition for two different regions in the hybrid discontinuous Galerkin scheme only requires consistent numerical
fluxes defined at the cell interfaces. For the domain decomposition,
it allows isolated hydrodynamic or kinetic cells, which will be very
convenient especially when extended to high-dimensional space problems
and/or on unstructured meshes. We will solve the kinetic equation in
the kinetic region with a discontinuous Galerkin scheme under the
asymptotic preserving
framework \cite{filbet2010class}, while solving the compressible Navier-Stokes equations,
which can be (formally) obtained by Chapman-Enskog expansion
\cite{cercignani1969mathematical,cercignani2000rarefied,chapman1970mathematical,bardos1991fluid},
by following the discontinuous Galerkin scheme of F. Bassi and
S. Rebay \cite{bassi1997high}. We apply an upwind numerical flux for
the kinetic solver, while to be consistent, especially in order to
match the numerical fluxes at the cell interface between two different
regions, we define a flux for the hydrodynamic solver obtained as the asymptotic limit of  the kinetic
flux. Let us emphasize that this particular hydrodynamic flux has already been introduced
in \cite{chou1997kinetic} as a kinetic flux splitting scheme. Furthermore,  we will simplify the computation of the moment realizability criterion for domain decomposition, by extracting some main derivatives in the moment realization matrix. For illustration in our numerical example, we take the BGK equation and reduce the 3D in velocity by the technique of Chu reduction \cite{chu1965kinetic}. Our new proposed hybrid discontinuous Galerkin scheme will be more robust on $h$-$p$ adaptivity, parallelization, and on high dimensional problems with unstructured meshes. We will perform some numerical tests on some physical relevant problems with shocks or boundary layers, such as flow caused by evaporation and condensation, 2D Riemann problem as well as 2D ghost effect. The results will demonstrate the efficiency and effectiveness of our proposed approach.

The rest of the paper is organized as follows. In Section \ref{sec2},
we recall the hydrodynamic limit of the multi-scale kinetic equation
based on Chapman-Enskog expansion, and describe the domain
decomposition from moment realizability criterion. In Section
\ref{sec4}, a discontinuous Galerkin scheme for the kinetic equation
and a discontinuous Galerkin scheme for the compressible Navier-Stokes equations will be formed, and consistent numerical fluxes for the interface coupling condition between two different regions are stated. Numerical tests are followed in Section \ref{sec5}. Conclusion and future work are made in Section \ref{sec6}.

%%%%%%%%%%%%%%%%%%%%%%%%%%%%%%%%%%%%%%%%%%
%
%%%%%%%%%%%%%%%%%%%%%%%%%%%%%%%%%%%%%%%%%%
\section{Hydrodynamic limit and domain decomposition}
\label{sec2}
\setcounter{equation}{0}
\setcounter{figure}{0}
\setcounter{table}{0}

In this section, we will recall the hydrodynamic limit of the kinetic equation \eqref{eq:bgk}. By formally doing the Chapman-Enskog expansion, the equation \eqref{eq:bgk}, to the zeroth order limit as $\eps\rightarrow 0$ will converge to the compressible Euler equations, while to the first order limit when $0<\eps\ll 1$, it is going to be the compressible Navier-Stokes equations.  In the following, we will briefly review the derivation of compressible Euler and Navier-Stokes limit \cite{filbet2015hierarchy,bardos1991fluid}, and describe the domain decomposition criteria to divide the domain into kinetic and hydrodynamic regions, which will be used to define our hybrid scheme.

\subsection{The hydrodynamic limit}
\label{sec2.1}
According to the conservation property \eqref{eq:cons}, by integrating in the velocity space, without any closure, we have
\beq
\label{eq:Ueps}
\left\{
	\begin{array}{l}
		\ds\pt \rho + \nablax \cdot(\rho \bfu)  \,= \, 0, 
		\\
		\,
		\\
		\ds\pt (\rho \bfu)+ \nablax \cdot \left(  \int_{\mathbb{R}^3}\bfv\otimes\bfv \,f(\bfv) \, d\bfv \right) \, = \, {\bf 0}_{\mathbb{R}^3}, 
		\\
		\,
		\\
		\ds\pt E + \nablax \cdot \left( \int_{\mathbb{R}^3}\frac12 |\bfv|^2 \, \bfv \, f(\bfv) \, d\bfv \right) \,= \, 0.
	\end{array}
\right.
\eeq
Let $\bfU$ denote the conservative macroscopic components of density, momentum and energy, which correspond to the first three moments of the distribution function $f$:
$$
\bfU \, := \, (\rho, \rho \bfu, E)^\text{T} = \int_{\mathbb{R}^{3}} m(\bfv) \, f(\bfv) \, d\bfv,
$$
where the energy is $E=\frac12\rho|\bfu|^2+\frac{3}{2}\rho\,T$. $f$ can be approximated by the Chapman-Enskog expansion, which is
\beq
\label{eq:feps}
f^\eps(\bfv) \, := \, \mathcal{M} \big [1 + \eps \, g^{(1)}  + \eps^2 \, g^{(2)} + \cdots \big],
\eeq
where $\mathcal{M}$ is in short of $\mathcal{M}[f]$ and the fluctuations satisfy
\[
\int_{\mathbb{R}^3} g^{(i)}(\bfv) \,m(\bfv) \, d\bfv = {\bf 0}^T_{\mathbb{R}^5}, \quad i = 1, 2, \cdots.
\]
Substituting $f^\eps(\bfv)$ into \eqref{eq:Ueps}, we get
\beq
\label{eq:U1}
\left\{
	\begin{array}{l}
		\ds\pt \rho + \nablax \cdot(\rho \bfu) \, = \, 0, 
		\\
		\,
		\\
		\ds\pt (\rho \bfu)+ \nablax \cdot \left( \rho \, \bfu
          \otimes \bfu + \rho\,T(\bfI + \bar{\bfA}^\eps) \right)  \,= \, {\bf 0}_{\mathbb{R}^3}, 
		\\
		\,
		\\
		\ds\pt E + \nablax \cdot \Bigg( \frac12\rho\,|\bfu|^2\,\bfu + \rho\,T\,\Big(\frac{3+2}{2}\bfI + \bar\bfA^\eps\Big)\bfu +\rho\,T^{3/2}\bar\bfB^\eps \Bigg) \,= \, 0,
	\end{array}
\right.
\eeq
where the traceless matrix $\bar\bfA^\eps\in\mathbb{M}^3$ and vector $\bar\bfB^\eps\in\mathbb{R}^3$ are 
\beq
\label{eq:AB}
\left\{
\begin{array}{l}
\ds\bar\bfA^\eps\,:=\,\frac{1}{\rho}\int_{\mathbb{R}^3}\bfA(\bfV)f^\eps(\bfv)d\bfv,\\ \,\\ 
\ds\bar\bfB^\eps\,:=\, \frac{1}{\rho}\int_{\mathbb{R}^3}\bfB(\bfV)f^\eps(\bfv)d\bfv,
\end{array}
\right.
\eeq
with
$$
\bfA(\bfV) = \bfV\otimes\bfV-\frac{|\bfV|^2}{3}\bfI, \quad  \bfB(\bfV)=\frac12\big[|\bfV|^2-(3+2)\big]\bfV
$$
and
\[
\bfV\,:=\,\frac{\bfv-\bfu}{\sqrt{T}}.
\]

In the zeroth order limit when $f^\eps(\bfv) = \mathcal{M}$, since the matrix $\bar\bfA^\eps$ is traceless, and due to $\bar\bfB^\eps$ only involving odd, centered moments of $f$, we obtain
\beq
\label{eq:AB:Euler}
\left\{
\begin{array}{l}
\ds\bar\bfA^\eps_{Euler} \,:=\, \frac{1}{\rho}\int_{\mathbb{R}^3}\bfA(\bfV)\mathcal{M}(\bfv)d\bfv \,= \, {\bf 0}^T_{\mathbb{M}^3}, 
\\ \, \\
\ds\bar\bfB^\eps_{Euler} \,:=\, \frac{1}{\rho}\int_{\mathbb{R}^3}\bfB(\bfV)\mathcal{M}(\bfv)d\bfv  \,= \, {\bf 0}^T_{\mathbb{R}^3},
\end{array}
\right.
\eeq
hence \eqref{eq:U1} gives us the compressible Euler system
\beq
\label{eq:Euler}
\left\{
	\begin{array}{l}
		\pt \rho + \nablax \cdot(\rho \bfu) \, = \, 0, 
		\\
		\,
		\\
		\pt (\rho \bfu)+ \nablax \cdot \left( \rho \, \bfu \otimes \bfu + \rho\,T\,\bfI \right) \, = \, {\bf 0}_{\mathbb{R}^3}, 
		\\
		\,
		\\
		\pt E + \nablax \cdot \Big( \bfu \big(E + \rho\,T\big) \Big) \,= \, 0.
	\end{array}
\right.
\eeq

In the first order limit, by plugging \eqref{eq:feps} into \eqref{eq:U1}, since the Maxwellian distribution $\mathcal{M}$ is an equilibrium of the collision operator $\mQ(\mathcal{M})=0$, the fluctuation $g^{(1)}$ is given by
\beq
\pt \mathcal{M} + \bfv\cdot\nablax \mathcal{M} \,=\, \mathcal{L}_\mathcal{M}g^{(1)} + \mathcal{O}(\eps),
\label{eq:Lg}
\eeq
where $\mathcal{L}_\mathcal{M}$ is the linearized collision operator around the Maxwellian distribution $\mathcal{M}$. Besides, we also have
$$
\pt \mathcal{M} + \bfv\cdot\nablax \mathcal{M}\,=\,\mathcal{M}\Big[\pt \rho + \bfv\cdot\nablax \rho + \frac{1}{\sqrt{T}}\left(\bfV\cdot\pt \bfu + \bfV\otimes\bfv:\nablax \bfu\right) + \frac{1}{2T}(|\bfV|^2-3)(\pt T + \bfv\cdot\nablax T)\Big].
$$
Replacing the time derivatives by spatial ones from \eqref{eq:U1}, and dropping all terms of order $\eps$ in \eqref{eq:Lg}, after some computations, we get
\beq
\label{eq:Lg1}
\mathcal{L}_\mathcal{M}g^{(1)}=\mathcal{M}\,\left[\bfA(\bfV):\bfD(\bfu)+2\bfB(\bfV)\cdot \frac{\nablax T}{\sqrt{T}}\right],
\eeq
where $\bfA$ and $\bfB$ are defined in \eqref{eq:AB} and the deformation tensor $\bfD(\bfu)$ is given by
\beq
\label{eq:Du}
\bfD(\bfu)\,:=\, \nablax \bfu + (\nablax \bfu)^T - \frac23(\nablax\cdot\bfu)\bfI.
\eeq
The linear operator $\mathcal{L}_\mathcal{M}$ is invertible on the orthogonal of its kernel \cite{filbet2015hierarchy},
where its kernel is 
\[
{\rm ker} \,\mathcal{L}_\mathcal{M}  \,=\, {\rm Span}\left\{\frac{1}{\rho},\frac{\bfV}{\rho},\frac{1}{2\rho}\big(|\bfV|^2-3\big)\right\},
\]
and $\bfA(\bfV), \,\bfB(\bfV) \,\in\, {\rm
  ker}^\perp\left(\mathcal{L}_\mathcal{M}\right)$. From \eqref{eq:Lg1}, it yields
\beq
\label{eq:g1}
g^{(1)}=\mathcal{L}^{-1}_\mathcal{M}\big(\mathcal{M}\bfA(\bfV)\big):\bfD(\bfu)+2\mathcal{L}^{-1}_\mathcal{M}\big(\mathcal{M}\bfB(\bfV)\big)\cdot \frac{\nablax T}{\sqrt{T}}.
\eeq
Plugging this expression into \eqref{eq:AB}, it gives
\beq
\label{eq:AB:NS}
\left\{
\begin{array}{l}
\ds\bar \bfA^\eps_{NS} \,:=\, \frac{1}{\rho}\int_{\mathbb{R}^3}\bfA(\bfV)\mathcal{M}(\bfv)\big[1+\eps g^{(1)}(\bfv)\big]d\bfv \,=\, -\eps\frac{\mu}{\rho\,T}\bfD(\bfu),
\\ \,\\
\ds\bar\bfB^\eps_{NS} \,:=\, \frac{1}{\rho}\int_{\mathbb{R}^3}\bfB(\bfV)\mathcal{M}(\bfv)\big[1+\eps g^{(1)}(\bfv)\big]d\bfv\,=\, -\eps\frac{\kappa}{\rho\,T^{3/2}}\nablax \,T.
\end{array}
\right.
\eeq
The viscosity and thermal conductivity coefficients are given by
\beq
\label{eq:coef}
\left\{
\begin{array}{l}
\ds\mu \,:=\,-T\int_{\mathbb{R}^3}\mathcal{M}(\bfv) \bfA(\bfV):\mathcal{L}^{-1}_\mathcal{M}\big(\mathcal{M}\bfA)(\bfv)d\bfv, \\ \, \\
\ds\kappa \,:=\,-T\int_{\mathbb{R}^3}\mathcal{M}(\bfv)\bfB(\bfv)\cdot\mathcal{L}^{-1}_\mathcal{M}\big(\mathcal{M}\bfB)(\bfv)d\bfv.
\end{array}
\right.
\eeq
Substituting \eqref{eq:AB:NS} into \eqref{eq:U1}, it gives us
the compressible Navier-Stokes system
\beq
\label{eq:cns}
\left\{
	\begin{array}{l}
		\ds\pt \rho + \nablax \cdot \, (\rho \, \bfu) \, = \, 0, \\
		\,
		\\
		\,
		\ds\pt (\rho \, \bfu) + \nablax \cdot \, (\rho \, \bfu\otimes\bfu + \rho \, T \, \bfI) \,=\, \eps \, \nablax \cdot({\bf \sigma}(\bfu)), \\
		\,
		\\
		\,
		\ds\pt E + \nablax \cdot \, (\bfu \, (E+\rho\, T)) \, = \, \eps \, \nablax \cdot({\bf \sigma}(\bfu) \cdot \bfu + {\bf q}),
	\end{array}
\right.
\eeq
where ${\bf \sigma} := -\mu \,\bfD(\bfu)$ is the viscosity tensor and ${\bf q} := -\kappa \nablax T$ is the heat flux. 

For the Boltzmann equation in the hard sphere case, $\mu$ and $\kappa$ can be expressed as \cite{golse2005boltzmann}
\beq
\label{eq:vis_heat:boltz}
\mu = \mu_0 \sqrt{T}, \qquad \kappa = \kappa_0 \sqrt{T},
\eeq
for some positive constants $\mu_0$ and $\kappa_0$. In the BGK case, $\mu$ and $\kappa$ are related to the macroscopic quantities of $\rho$ and $T$ by \cite{struchtrup2005macroscopic}
\beq
\label{eq:vis_heat:bgk}
\mu = \frac{1}{1-\beta}\frac{\rho T}{\nu}, \qquad \kappa = \frac52 \frac{\rho T}{\nu},
\eeq
and the collision frequency can be taken as $\nu = \frac{2}{\sqrt{\pi}}\rho$ \cite{aoki1994gas}.

The compressible Navier-Stokes equations \eqref{eq:cns} can be written in a vector form as
\beq
\label{eq:cnsvform}
\pt \bfU + \nablax \cdot\bfF^a(\bfU) \,=\, \eps \, \nablax \cdot \bfF^d(\bfU,\nablax \bfU),
\eeq
with the advection flux $\bfF^a(\bfU)$ and the viscous diffusion flux $\bfF^d(\bfU,\nablax \bfU)$
to be
$$
\bfF^a(\bfU) \,= \, \begin{pmatrix}
\rho \bfu \\

\rho \bfu \otimes \bfu + \rho\, T \, \bfI \\

\bfu \, (E + \rho \, T)
\end{pmatrix},
\qquad
\bfF^d(\bfU,\nablax \bfU) \, = \, \begin{pmatrix}
0 \\

{\bf \sigma}(\bfu) \\

{\bf \sigma}(\bfu)\cdot \bfu + {\bf q}
\end{pmatrix}.
$$
Moreover, if we define
\beq
\bfF(\bfU,\nablax \bfU) \, = \, \bfF^a(\bfU) - \eps \, \bfF^d(\bfU, \nablax \bfU),
\eeq
it becomes
\beq
\label{eq:cns1}
\pt \bfU + \nablax \, \bfF(\bfU, \nablax \bfU) \,=\, {\bf 0}.
\eeq
The flux function $\bfF(\bfU,\nablax \bfU)$ corresponds to the kinetic description \eqref{eq:Ueps} by taking the first order distribution function $f^{\eps}(\bfv)$ \eqref{eq:feps}, that is
\beq
\label{eq:kflux}
\bfF(\bfU,\nablax \bfU) \,=\, \int_{\mathbb{R}^{3}} \bfv \, m(\bfv)  \, \big[\mathcal{M}(1+\eps g^{(1)})\big](\bfv) \, d\bfv.
\eeq
Formally by letting $\eps\rightarrow 0$,  \eqref{eq:cns1} corresponds to the compressible
Euler equations
\beq
\label{eq:eulervform}
\pt \bfU + \nablax \cdot\bfF^a(\bfU) \,=\, \bf0.
\eeq

\subsection{Domain decomposition}
\label{sec2.2}
We now follow the domain decomposition indicator as defined in \cite{filbet2015hierarchy},
to divide the computation domain into kinetic and hydrodynamic regions, where in the kinetic region (mainly around shocks or boundary layers) the kinetic equation \eqref{eq:bgk} is used,
otherwise in the hydrodynamic region, the low computational cost
compressible Navier-Stokes equations \eqref{eq:cns} or \eqref{eq:cns1} as the hydrodynamic equation will be solved. We review the main principle on how to define the criteria, however we would like to simplify the tedious computation of the eigenvalues for the moment realizability matrix.  

From fluid to kinetic, the hydrodynamic breakdown criterion is defined by the moment realizability matrix. The deviation of the eigenvalue for the moment realizability matrix measures the appropriateness of the corresponding hydrodynamic limit. The compressible Euler system has all eigenvalues to be $1$ with an identity moment realizability matrix.  Here since we use
the compressible Navier-Stokes equations as the fluid model, we only
concern about the deviation of the eigenvalues between the first order
compressible Navier-Stokes equations and the second order Burnett equations. We recall the moment realization matrix from \cite{filbet2015hierarchy}, which is given by
\beq
\label{eq:eigen}
\mathcal{V} \,:=\, \bfI \, + \, \bar{\bfA}^\eps \,-\, \frac{2}{3\bar{C}^\eps} \bar{\bfB}^\eps \otimes \bar{\bfB}^\eps.
\eeq
For compressible Navier-Stokes equations, $\bar{\bfA}^\eps_{NS}$ and $\bar{\bfB}^\eps_{NS}$ are defined in \eqref{eq:AB:NS},
while for the second order Burnett equations,
$\bar{\bfA}^\eps_{Burnett}$ and $\bar{\bfB}^\eps_{Burnett}$ are \cite{filbet2015hierarchy} 
\beq
	\label{eq:AB:burnett}
	\left\{
	\begin{array}{ll}
	\ds\bar{\bfA}^\eps_{Burnett} \,:= \,& \ds\, -\eps \frac{\mu}{\rho
          \,T} \bfD(\bfu) - 2\eps^2\frac{\mu^2}{\rho^2\,T^2}\Bigg\{
          -\frac{T}{\rho}\textrm{Hess}_\bfx(\rho) + \frac{T}{\rho^2}
          \nablax \rho \otimes \nablax \rho-\frac{1}{\rho}\nablax T
          \otimes \nablax \rho \\ \, \\ \, &\ds \,+\, \nablax \bfu \,\nablax
         \ds \bfu^T-\frac13\bfD(\bfu)\,
          \textrm{div}_\bfx(\bfu)+\frac{1}{T}\nablax T \otimes \nablax
          T \Bigg\},  \\ \, \\
	\ds\bar{\bfB}^\eps_{Burnett} \,:= &\ds\, -\eps \frac{\kappa}{\rho
          \, T^{3/2}}\nablax T - \eps^2 \frac{\mu^2}{\rho^2 \,
          T^{5/2}} \Bigg\{ \frac{25}{6}\, \textrm{div}_\bfx \bfu\,
          \nablax T - \frac{5}{3}\Big[ T\textrm{div}_\bfx(\nablax
          \bfu) + \textrm{div}_\bfx \bfu \,\nablax T  \\ \, \\
	\,&\ds\,+\, 6(\nablax \bfu) \, \nablax T\Big]+\frac{2}{\rho}\bfD(\bfu)\nablax(\rho\,T)+2T\,\textrm{div}_\bfx(\bfD(\bfu))+16\bfD(\bfu)\,\nablax T \Bigg\}. 
	\end{array}
\right.
\eeq
Noticing that for the main components $\bar{\bfA}^\eps$ and $\bar{\bfB}^\eps$ in \eqref{eq:eigen}, the major differences between $\bar{\bfA}^\eps_{NS}$ and $\bar{\bfA}^\eps_{Burnett}$ as well as  $\bar{\bfB}^\eps_{NS}$ and $\bar{\bfB}^\eps_{Burnett}$ are the $\mathcal{O}(\eps^2)$ terms appeared in \eqref{eq:AB:burnett}. As an indicator, we may only need to extract the crucial spatial derivatives appeared in those terms and roughly measure their magnitudes, to simplify the tedious computations of all terms, thus avoiding too much computational cost on the domain decomposition indicator along this side. From \eqref{eq:vis_heat:boltz} and \eqref{eq:vis_heat:bgk}, either $\mu/\sqrt{T}$ and $\kappa/\sqrt{T}$ for the Boltzmann operator or $\nu \mu/(\rho\,T)$ and $\nu \kappa/(\rho\,T)$ for the BGK operator are $\mathcal{O}(1)$ if we assume macroscopic quantities $\rho$ and $T$ are $\mathcal{O}(1)$, a new indicator from fluid to kinetic we are now using is defined as
$$
\lambda_{\eps^2}(t,\bfx) \, := \,  \eps^2 \Bigg( \frac{|\nablax T|^2}{T} + |\nablax \bfu|^2 + \sqrt{\Big(|\Delta_\bfx \bfu|^2 + \Big|\Delta_\bfx \rho/\rho\Big|^2\Big)(1+T^2)} \Bigg).
$$
We take the fluid model of compressible Navier-Stokes equations to be not appropriate at some point $(t,\bfx)$ if 
\beq
\label{eq:crit1}
|\lambda_{\eps^2}(t,\bfx)| \, > \, \eta_0.
\eeq
$\lambda_{\eps^2}(t,\bfx)$ is measured at the cell center in the numerical section and we take $\eta_0 = 10^{-3}$. 

From kinetic to fluid, the criterion is that a kinetic description at some point $(t,\bfx)$ corresponds to a hydrodynamic equations if
\beq
\label{eq:crit2}
\| f(t,\bfx,\cdot) - f_2(t,\bfx,\cdot) \|_{L^2} \, \le \, \delta_0,
\eeq
where $\delta_0$ is a small parameter and we take $\delta_0=10^{-3}$
and the function $f_2=\mathcal{M}[1+\eps g^{(1)}]$ corresponds to the Chapman-Enskog expansion \eqref{eq:feps} up to first order.

%%%%%%%%%%%%%%%%%%%%%%%%%%%%%%%%%%%%%%%%%%
%
%%%%%%%%%%%%%%%%%%%%%%%%%%%%%%%%%%%%%%%%%%
\section{Hybrid discontinuous Galerkin scheme}
\label{sec4}
\setcounter{equation}{0}
\setcounter{figure}{0}
\setcounter{table}{0}

In this section, we will propose our hybrid discontinuous Galerkin
scheme for the multi-scale kinetic equations, namely we define either an approximation of the kinetic equation
\eqref{eq:bgk} according to the value of the indicator
\eqref{eq:crit1}  or another one of
the fluid equations \eqref{eq:cns1} in the  region where
\eqref{eq:crit2} is satisfied. Due to the compactness of discontinuous Galerkin schemes, only consistent numerical fluxes need to be defined at the cell interface between two regions, which is also important to ensure mass conservation. 

We adopt the discontinuous Galerkin scheme for spatial discretization, as it is more
convenient than a finite volume scheme used in
\cite{filbet2015hierarchy}, especially when extending to high
dimensions in space and in the domain decomposition. Indeed, the finite volume
method requires either the kinetic or fluid region to be at least as
wide as the stencil of the scheme. Besides, the discontinuous Galerkin
scheme is $h$-$p$ adaptive and very flexible to nonuniform meshes,
making it more suitable to physical problems with boundary layers (see
Section \ref{sec5}).

Our hybrid discontinuous Galerkin scheme can be defined for the
multi-scale kinetic equation with either the full Boltzmann collision operator or the (ES-) BGK operator. The Boltzmann operator \eqref{boltzoper} in a bounded domain can be discretized by
the method introduced in \cite{filbet2012deterministic} and references therein, and in an asymptotic preserving framework can be penalized by the BGK operator \cite{filbet2010class}. We will only take the BGK operator \eqref{bgkoper} in the kinetic equation and describe the discontinuous Galerkin scheme for the BGK equation \eqref{eq:bgk}. The ES-BGK operator \eqref{esbgkoper} can be done similarly and we can use the technique in \cite{filbet2010class} for the Boltzmann operator. The discontinuous Galerkin scheme for the hydrodynamic equations \eqref{eq:cns1} is followed. Only the discontinuous Galerkin scheme in 2D in space is presented, the 1D case can be easily deduced from the 2D ones. We take the kinetic flux-vector splitting \cite{chou1997kinetic} to define the hydrodynamic numerical flux at the cell interface between two hydrodynamic cells or two cells between two regions. A nodal discontinuous Galerkin scheme \cite{hesthaven2008nodal}  is implemented to produce the numerical results.

For the velocity space, we will take a large enough cut-off domain $\Omega_v=[-V_c, V_c]$ and discretize it uniformly with $N_v$ points,  $\{v^j\}_{j=1}^{N_v}$, along each $v$ direction. We integrate the velocity space by mid-point rule, which is spectrally accurate for smooth functions with periodic
boundary conditions or compact supports \cite{boyd2001chebyshev}, as there is no differential operator acting on $\bfv$. However, we note that the domain cut-off and the discretization in $\bfv$ will lead to conservation of mass, momentum and energy not exactly but approximately up to the integral error along $\bfv$ direction. In the following, for easy presentation, we still keep $\bfv$ to be continuous and only discuss the discretizations in space and in time.

\subsection{Preliminary}
\label{sec4.1}
We consider the space domain in 2D to be a rectangular $\Omega_\bfx=[a,\, b] \times [c,\, d]$, and divide it by $a=x_{1/2}<x_{{3}/{2}}<\cdots <x_{N_x+1/2}=b$ and $c=y_{1/2}<y_{{3}/{2}}<\cdots <y_{N_y+1/2}=d$. Let $I_{i,j}=I_i \times I_j=[\xL, \xR]\times[\yL,\yR]$ denote an element with its length $\Delta x_i = \xR - \xL$ and $\Delta y_j = \yR - \yL$. Let $h_x=\max_{i=1}^{N_x} \Delta x_i$, $h_y=\max_{j=1}^{N_y} \Delta y_j$ and $h=\max(h_x,\, h_y)$. Given any non-negative integer vector $\bfK=(K_1,K_2)$, we define a finite dimensional discrete space,
\begin{equation}
\mZ_h^\bfK \,=\, \left\{w \in L^2(\Omega_\bfx)\,:\,\, w|_{I_{i,j}}\in
  Q^\bfK(I_{i,j}), \, 1\leq i\leq N_x,\, 1\leq j\leq N_y\right\}.
\label{eq:DiscreteSpace2D}
\end{equation}
Its vector version is denoted as
\begin{equation}
\bfZ_h^\bfK \,=\, \left\{{\bf w}=(w_1,w_2,w_3,w_4)^T \,:\,\, w_l \in \mZ_h^\bfK, \, 1 \leq l \leq 4\right\}.
\end{equation}
For simplicity, we take the 2D local space $Q^\bfK(I_{i,j})$ as a tensor product space of $P^{K_1}(I_i)\otimes P^{K_2}(I_j)$, where $P^K(I)$ consists of polynomials of degree at most $K$ on $I$. Functions in $\mZ_h^\bfK$ are piecewise defined and may be discontinuous across cell interfaces. The left and right limits of a function $u \in \mZ_h^\bfK$ at the interface $(x_{i+\frac{1}{2}},y)$ along $y$ direction are denoted as $u(x_\iR^\pm,y) =\lim_{\varepsilon\rightarrow \pm 0}u(x_{i+\frac{1}{2}}+\varepsilon, y)$, similarly for $u(x, y^\pm_\jR)$.

We will use a nodal basis to represent functions in the discrete space $P^K(I)$, and approximate integrals by numerical quadratures. Specifically, we choose the local nodal basis (also called Lagrangian basis) $\{\phi_i^k(x)\}_{k=0}^K$ associated with the $K+1$ Gaussian quadrature points $\{x^k_i\}^K_{k=0}$ on $I_i$, defined as below
\begin{equation}
\phi_i^k(x) \in P^K(I_i), \quad\textrm{and} \quad \phi_i^k (x_i^{k'})=\delta_{k k'},\quad k, k'\,=\,0,\,1,\cdots,\, K,
\label{bases}
\end{equation}
where $\delta_{k k'}$ is the Kronecker delta function. We denote by
$\{\omega_k\}^K_{k=0}$ the corresponding quadrature weights on the
reference element $(-{1}/{2}, {1}/{2})$. For the two dimensional local
nodal basis, we choose a tensor product of the one dimensional local
nodal basis \eqref{bases} but associated with $K_1+1$ and $K_2+1$
Gaussian quadrature points on $I_{i}$ and $I_{j}$ along $x$ and $y$
directions respectively. 

\subsection{Discontinuous Galerkin for the kinetic equation}
\label{sec4.2}
For the kinetic equation \eqref{eq:bgk} with the BGK operator, we take
an implicit-explicit (IMEX) time discretization as in
\cite{filbet2010class,filbet2015hierarchy}. Considering the continuous
problem in phase space $(\bfx,\bfv)$, a first order implicit/explicit
scheme is defined as follows to pass from time level $t^n$ to $t^{n+1}$
\beq
\label{eq:imex12d}
\begin{cases}
	\begin{array}{ll}
		f^{n+1}(\bfv) & = \,\ds \frac{\eps}{\eps+\nu^{n+1}\Delta t} \Big( f^n(\bfv) - \Delta t \, \bfv \cdot \nablax \,f^n(\bfv) \Big) + \frac{\nu^{n+1}\Delta t}{\eps+\nu^{n+1}\Delta t} \, \mathcal{M}(\bfv,\bfU^{n+1}), \\
		\,
		\\
		\,
		f^0(\bfv) & = \,\ds f(0,\bfx,\bfv).
	\end{array}
\end{cases}
\eeq
The implicit Maxwellian $\mathcal{M}(\bfv,\bfU^{n+1})$ is first computed
from the macroscopic quantity $\bfU^{n+1}$, where
\beq
\label{eq:imex1U2d}
\bfU^{n+1} \,:=\, (\rho^{n+1}, (\rho \, \bfu)^{n+1}, E^{n+1})^\text{T} \, = \, \int_{\mathbb{R}^3} m(\bfv) \, ( f^n - \Delta t \, \bfv \cdot  \nablax \, f^n) \,d\bfv.
\eeq
For each $\bfv$, if we define
\beq
\label{eq:R2d}
\bfR(\bfv) \, := \, f(\bfv) - \Delta t \, \bfv \cdot \nablax \, f(\bfv),
\eeq
in vector form we have 
\beq
\label{eq:fv}
\begin{cases}
\begin{array}{ll}
\bfR^{n+1}(\bfv) & = \,\ds f^n(\bfv) - \Delta t \, \bfv \cdot \nablax \, f^n(\bfv), \\ \,\\
\bfU^{n+1}   & = \,\ds \int_{\mathbb{R}^3} m(\bfv) \, \bfR^{n+1}(\bfv) \, d\bfv, \\ \, \\
f^{n+1}(\bfv) & = \,\ds \frac{\eps}{\eps+\nu^{n+1}\Delta t} \, \bfR^{n+1}(\bfv) + \frac{\nu^{n+1}\Delta t}{\eps+\nu^{n+1}\Delta t} \, \mathcal{M}(\bfv, \bfU^{n+1}).
\end{array}
\end{cases}
\eeq
Following this strategy and for simplicity keeping $\bfv
\in\mathbb{R}^3$, a  two dimensional in space discontinuous Galerkin
scheme for \eqref{eq:fv} is defined as follows: for given
$f^n_h(\bfv)$, we find $f^{n+1}_h(\bfv)$, by computing a discrete
approximation $\bfR_h^{n+1}$, which solves the
following problem :   for any $\zeta\in \mZ_h^\bfK $ and for $1\leq  i\leq N_x$,
$1\leq j\leq N_y$,
\begin{eqnarray}
\label{eq2d:DGf:a}
\int_{I_{i,j}} \bfR^{n+1}_h(\bfv) \zeta(\bfx) \, d\bfx &= &
                                                            \int_{I_{i,j}} f^n_h(\bfv) \zeta(\bfx)  d\bfx  \,+\, \Delta t \, \int_{I_{i,j}} \bfv \cdot \nablax \zeta(\bfx) \, f^n_h(\bfv) \, d\bfx \\   
&-&  \Delta t \, \int_{I_i} \widetilde{(v_1 f)}(x_{\iR},y)\zeta(x^-_\iR, y) -
   \widetilde{(v_1 f)}(x_\iL,y) \zeta(x^+_\iL, y) \, dy    \notag
\nonumber
\\   
&-& \Delta t \, \int_{I_j} \widetilde{(v_2 f)}(x,y_{\jR})\zeta(x, y^-_\jR) -
   \widetilde{(v_2 f)}(x,y_\jL) \zeta(x, y^+_\jL) \, dx,
   \nonumber
\end{eqnarray}
then we compute  $\bfU^{n+1}_h$ given for any $\beta \in \mZ_h^\bfK $ and for $1\leq  i\leq N_x$,
$1\leq j\leq N_y$, by
\beq
	\int_{I_{i,j}}\bfU^{n+1}_h \beta(\bfx) \, d\bfx \,  =  \,
        \int_{I_{i,j}} \int_{\mathbb{R}^3} m(\bfv)
          \bfR^{n+1}_h(\bfv) d\bfv \, \beta(\bfx) \, d\bfx.
\label{eq2d:DGf:b} 
\eeq
Finally the discrete distribution function $f_h^{n+1}$ is defined such that  for any $\alpha \in \mZ_h^\bfK $ and for $1\leq  i\leq N_x$,
$1\leq j\leq N_y$, 
\begin{eqnarray}
\label{eq2d:DGf:c} 
\int_{I_{i,j}} f^{n+1}_h(\bfv) \alpha(\bfx) \, d\bfx &=&
  \int_{I_{i,j}}\frac{\eps}{\eps+\nu^{n+1}\Delta t}\bfR^{n+1}_h(\bfv)\alpha(\bfx)
  d\bfx \\
 &+&\int_{I_{i,j}} \frac{\nu^{n+1}\Delta t}{\eps+\nu^{n+1}\Delta t} \mathcal{M}(\bfv, \bfU^{n+1}_h) \, \alpha(\bfx) \, d\bfx,  
\nonumber
\end{eqnarray}
where $\widetilde{v f}$ is an upwind numerical flux along its direction, 
\beq
\label{eq2d:upwind}
\widetilde{v f}:=
\left\{
\begin{array}{ll}
	v \,  f^-, &\mbox{if}\;  v \ge 0,\\
	v \,  f^+, &\mbox{if}\;  v <0,
\end{array}
\right.
\eeq
whereas $f^\pm$ are the left and right limits of $f^n_h$ at the cell interface of two adjacent cells respectively. 

For a nodal discontinuous Galerkin scheme with Lagrangian bases, the
scheme \eqref{eq2d:DGf:a} becomes
\begin{eqnarray*}
%\label{eq2d:DGf:nodal:a}
\bfR^{n+1}_{i,k_1,j,k_2}(\bfv) &= &  f^n_{i,k_1,j,k_2}(\bfv)  + \frac{1}{\omega_{k_1}\omega_{k_2}}\frac{ \Delta t }{\Delta x_i \Delta y_j} \Bigg(\sum_{l_1=0}^{K_1}\sum_{l_2=0}^{K_2} \bfv \cdot \nablax (\phi_i^{k_1}(x_i^{l_1}) \phi_j^{k_2}(y_j^{l_2}) ) \, \bfg^n_{i,l_1,j,l_2}(\bfv) \\ 
	& -& \Delta y_j \sum_{l_2=0}^{K_2} \omega_{l_2}
           \Big(\widetilde{(v_1 f)}(x_{\iR}, y_j^{l_2}) \,
           \zeta(x^-_\iR, y_j^{l_2})- \widetilde{(v_1
           f)}(x_\iL,y_j^{l_2}) \,\zeta(x^+_\iL,y_j^{l_2}) \Big) 
%\nonumber
	\\ & -& \Delta x_i \sum_{l_1=0}^{K_1} \omega_{l_1}
              \Big(\widetilde{(v_2 f)}(x_i^{l_1}, y_\jR) \,
              \zeta(x_i^{l_1}, y^-_\jR)- \widetilde{(v_2
              f)}(x_i^{l_1},y_\jL) \, \zeta(x_i^{l_1},y^+_\jL) \Big)
              \Bigg),  
%\nonumber 
\end{eqnarray*}
whereas $\bfU$ is given by
%\beq
$$
\bfU^{n+1}_{i,k_1,j,k_2} \, =  \,  \int_{\mathbb{R}^3} m(\bfv) \bfR^{n+1}_{i,k_1,j,k_2}(\bfv) \,d\bfv,
%	\label{eq2d:DGf:nodal:b} 
$$
%\eeq 
and the discrete approximation $f_h^{n+1}$ is defined by
%\beq
$$
f^{n+1}_{i,k_1,j,k_2}(\bfv) \, =  \,
\frac{\eps}{\eps+\nu^{n+1}_{i,k_1,j,k_2}\Delta t}
\bfR^{n+1}_{i,k_1,j,k_2}(\bfv) + \frac{\nu^{n+1}_{i,k_1,j,k_2}\Delta
  t}{\eps+\nu^{n+1}_{i,k_1,j,k_2}\Delta t}
\mathcal{M}(\bfv,\bfU^{n+1}_{i,k_1,j,k_2}), 
%\label{eq2d:DGf:nodal:c} 
%\eeq
$$
where $f^n_{i,k_1,j,k_2}(\bfv), \bfR^{n+1}_{i,k_1,j,k_2}(\bfv),
\bfU^{n+1}_{i,k_1,j,k_2}$ with subindex $(i,k_1,j,k_2)$ are the
corresponding numerical values at the $(k_1,k_2)$-th Gaussian
quadrature point in cell $I_{i,j}$, for $0\le k_1 \le K_1$ and $0\le
k_2 \le K_2$. 

\subsection{Discontinuous Galerkin for the compressible Navier-Stokes equations}
\label{sec4.3}
We will follow F. Bassi and S. Rebay \cite{bassi1997high} to define a
discontinuous Galerkin scheme for the compressible Navier-Stokes equations \eqref{eq:cns1}. 
Letting $\nablax \bfU:=(\bfS_1, \bfS_2)$, with an Euler forward time
discretization, the discontinuous Galerkin scheme for \eqref{eq:cns1}
is defined as follows:  we seek $\bfU^{n}_h\in {\bfZ}_h^\bfK$, such that for any $\eta(\bfx) \in \mZ_h^{\bfK}$ and $1\leq i\leq N_x$, $1\leq j\leq N_y$, we have
\begin{eqnarray}
\label{eq2d:SDG:a}
\int_{I_{i,j}}  \bfU^{n+1}_h \, \eta(\bfx) \, d\bfx&=& \int_{I_{i,j}}
                                                        \bfU^n_h \,
                                                        \eta(\bfx)
                                                        \, d\bfx \\
&+&
                                                        \Delta t
                                                        \left(\int_{I_{i,j}}
                                                        \bfF(\bfU^n_h,
                                                        \bfS^n_{h})
                                                        \cdot \nablax
                                                        \eta(\bfx)
                                                        \, d\bfx -
                                                        \oint_{\partial
                                                        I_{i,j}}
                                                        \eta(\bfx)
                                                        \, \hat{\bfF}
                                                        \cdot \bfn \,
                                                        d\sigma
                                                        \right). 
\nonumber
\end{eqnarray}
whereas $\nablax \bfU^n_h=(\bfS^n_{1,h},\bfS^n_{2,h})$ is such that
$\bfS^n_{1,h}, \bfS^n_{2,h}\in {\bfZ}_h^\bfK$  for any $\eta(\bfx) \in \mZ_h^{\bfK}$ and $1\leq i\leq N_x$, $1\leq j\leq N_y$, we have
\begin{eqnarray}
	\label{eq2d:SDG:b}
	\int_{I_{i,j}} \bfS^n_{1, h} \, \eta(\bfx) \, d\bfx&=& -
        \int_{I_{i,j}} \bfU^n_h \, \px \eta(\bfx) \,d\bfx 
\\
&+&
        \int_{I_{j}} \eta(x^-_\iR, y) \, \hat {\bfU}(x_\iR,y) -
        \eta(x^+_\iL,y) \, \hat {\bfU}(x_\iL, y)  \, dy,  
\nonumber
\end{eqnarray}
and
\begin{eqnarray}
\label{eq2d:SDG:c}
\int_{I_{i,j}} \bfS^n_{2, h} \, \eta(\bfx) \, d\bfx &=& -
                                                        \int_{I_{i,j}}
                                                        \bfU^n_h \,
                                                        \py \eta(\bfx)
                                                        \,d\bfx 
\\
\nonumber
&+& \int_{I_{i}} \eta(x, y^-_\jR) \, \hat {\bfU}(x,y_\jR) -
    \eta(x,y^+_\jL) \, \hat {\bfU}(x, y_\jL)  \, dx,  
\end{eqnarray}
where $\hat{\bf U}$ is taken to be a central flux 
\[
\hat{\bfU} \,:= \, \frac12(\bfU^++\bfU^-),
\]
along $x$ or $y$ direction. $\bfU^\pm$ are the left and right limits
of $\bfU^n_h$ at the cell interface.  Finally to construct  a nodal discontinuous Galerkin scheme, we apply
a quadrature formula as it has been done in the previous section. 

The flux $\hat{\bf F}$ is defined
in the next subsection. An important issue is to construct a numerical
flux which is consistent with the kinetic flux in the asymptotic limit
$\varepsilon \rightarrow 0$ in order to avoid spurious oscillations of
order $\varepsilon$ or  $\varepsilon^2$.

\subsection{Interface coupling condition}
\label{sec4.4}
For a discontinuous Galerkin scheme, due to its compactness, we only need to pass an interface coupling condition between two cells of different regions, that is, to define consistent numerical fluxes at such a cell interface for both the kinetic equation and the hydrodynamic system. Let us describe how to define the numerical fluxes 
$\widetilde{v_1 f}$ and $\hat{\bfF}_1$ in the $x$ direction, while $\widetilde{v_2 f}$ and $\hat{\bfF}_2$ in the $y$ direction can be defined similarly.

We assume $I_{i,j}$ is a kinetic cell and $I_{i+1,j}$ is a hydrodynamic cell, and the cell interface
is at $\{x_\iR\}\times [\yL, \yR]$:
\begin{itemize}
	\item For the discontinuous Galerkin scheme of the kinetic
          equation \eqref{eq2d:DGf:a}-\eqref{eq2d:DGf:c}, an upwind
          flux \eqref{eq2d:upwind} in \eqref{eq2d:DGf:a} is needed,
	\beq
	\label{eq2d:upwind2}
	\widetilde{v_1 f}:=
	\left\{
	\begin{array}{ll}
		v_1 \,  f^-,  &\mbox{if}\; v_1 \ge 0,\\
		v_1 \,  f^+, &\mbox{else},
	\end{array}
	\right.
	\eeq 
at each Gaussian quadrature point within $[\yL, \yR]$. For $v_1\ge 0$, $f^-$ in cell $I_{i,j}$ can be defined from $f^n_h$, while $f^+$ for $v_1 < 0$ in cell $I_{i+1,j}$ cannot.  In this case, we take $f^n_h$ in cell $I_{i+1,j}$ to be the first order truncated distribution functions as defined in \eqref{eq:feps} with $g^{(1)}$ defined in \eqref{eq:g1}, so that $f^+$ can be recovered from $v_1$, $\bfU(x^+_\iR,\cdot)$ and $\bfS(x^+_{\iR},\cdot)$, where $\bfS = \nablax \bfU$. 
	
	\item For the discontinuous Galerkin scheme of the fluid
          equations \eqref{eq2d:SDG:a}-\eqref{eq2d:SDG:c}, for
          conservation, the flux function $\hat{\bfF}_1$ in \eqref{eq2d:SDG:a} is
          defined as the integral of the kinetic upwind flux
          \eqref{eq2d:upwind2} on $\bfv$, 
	\beq
	\label{eq2d:Hflx}
	\hat \bfF_1(x_\iR, \cdot) = \int_{v_1 \ge 0} \, v_1 \, m(\bfv) \, f(x^-_{\iR},\cdot,\bfv) \, d\bfv + \int_{v_1 < 0} \, v_1 \, m(\bfv) \, \, f(x^+_{\iR},\cdot,\bfv) \,d\bfv.
	\eeq
	Based on $f^\pm$ in \eqref{eq2d:upwind2}, it can be obtained. 
\end{itemize}
The other way when $I_{i,j}$ is a hydrodynamic cell and $I_{i+1,j}$ is a kinetic cell can be defined symmetrically.
The hydrodynamic flux $\hat \bfF_1(x_\iR, \cdot)$ at the interface of
two cells inside the hydrodynamic region is defined in the same form
as \eqref{eq2d:Hflx}, with $f^\pm$ both taken as the first order
truncated distribution function \eqref{eq:feps} in cell $I_{i,j}$ and
$I_{i+1,j}$ respectively, which is known as the kinetic flux-vector
splitting flux \cite{chou1997kinetic}. 

This choice is particularly important to avoid that some spurious
oscillations appear when the coupling between kinetic and fluid
regions occur. Furthermore, this numerical  flux  is consistent with the continuous flux
given in \eqref{eq:kflux}. Indeed, it is given by,
\begin{eqnarray*}
\hat \bfF_1(x_\iR, \cdot) &=& \int_{v_1 \ge 0} \, v_1 \, m(\bfv) \,
                              \big[\mathcal{M}\,(1\,+\,\eps
                              g^{(1)})\big](x^-_{\iR},\cdot,\bfv) \, d\bfv 
\\
&+& \int_{v_1 < 0} \, v_1 \, m(\bfv) \, \, \big[\mathcal{M}(1+\eps g^{(1)})\big](x^+_{\iR},\cdot,\bfv) \,d\bfv.
\end{eqnarray*}
When we consider the BGK equation, we have 
\begin{eqnarray*}
\hat \bfF_1(x_\iR, \cdot) &=& \int_{v_1 \ge 0} \, v_1 \, m(\bfv) \,
                              \left[\mathcal{M}\,\left(1+\frac{\eps}{\nu}\left(\bfA(\bfV):\bfD(\bfu)+2\bfB(\bfV)\cdot \frac{\nablax T}{\sqrt{T}}
                              \right)\right)\right](x^-_{\iR},\cdot,\bfv) \, d\bfv 
\\
&+& \int_{v_1 < 0} \, v_1 \, m(\bfv) \, \left[\mathcal{M}\left(1+\frac{\eps}{\nu}\left(\bfA(\bfV):\bfD(\bfu)+2\bfB(\bfV)\cdot \frac{\nablax T}{\sqrt{T}}
                              \right)\right)\right](x^+_{\iR},\cdot,\bfv) \,d\bfv,
\end{eqnarray*}
where $\bfA$ and $\bfB$ are given in (\ref{eq:AB}). The resulting integrals on $v_1$, either on $[0,\infty]$ or on $[-\infty,0]$, can be expressed in terms of the
Gauss error function 
\[
\textrm{erfc}(x) \, = \, \frac{2}{\sqrt{\pi}} \int_x^\infty e^{-t^2} dt
\]
and some of its related functions. We can get explicit expressions in terms of $\textrm{erfc}(x)$ to avoid integration on $v_1$.
We mainly need to do some integrals in the following form
$$
%\beq
\frac{\rho}{\sqrt{2\pi T}} \int_{v\ge 0} v^n \, \left(\frac{v-u}{\sqrt{T}}\right)^m e^{-\frac{(v-u)^2}{2T} }\, dv
%\label{eq:A1}
%\eeq
$$
for $0 \le n \le 3$ and $0 \le m \le 3$.  First by letting
\[
v \,=\, z \, \sqrt{2T} + u,
\]
It can be transformed to
$$
%\beq
%\label{eq:A2}
\frac{\rho}{2}\,\frac{2}{\sqrt{\pi}}\int_{-\frac{u}{\sqrt{2T}}}^\infty (z \sqrt{2T}+u)^n \, (\sqrt{2} z)^m \, e^{-z^2}\,dz.
%\eeq
$$
Except the coefficient $\frac{\rho}{2}$, it is in the form of
\beq
\label{eq:A3}
\frac{2}{\sqrt{\pi}} \int_s^\infty (a\,z+b)^n \, (c\,z)^m\,e^{-z^2}\,dz,
\eeq
which is the Gauss error function $\textrm{erfc}(s)$ in case of $n=0$ and $m=0$.
Expanding on $(a\,z+b)^n$ for $n=0, 1,2,3$, what we need are
$$
%\beq
%\label{eq:A4}
\int_s^\infty z^n \,e^{-z^2}\,dz
%\eeq
$$
for $1\le n \le 6$. From integration by parts,  they are
$$
\left\{
%\begin{subequations*}
%	\label{eq:A5}
	\begin{array}{ll}
	\ds\int_s^\infty z   \, e^{-z^2} \, = \, & \ds\frac12 e^{-s^2}, \\
          \, \\
	\ds\int_s^\infty z^2 \, e^{-z^2} \, = \, & \ds\frac12 (s \,e^{-s^2}
                                                + \frac{\sqrt{\pi}}{2}
                                                \, \textrm{erfc}(s)),
          \\ \, \\
	\ds\int_s^\infty z^3 \, e^{-z^2} \, = \, & \ds\frac12 (1+s^2)
                                                e^{-s^2}, \\ \, \\
	\ds\int_s^\infty z^4 \, e^{-z^2} \, = \, & \ds\frac12
                                                \left(\left(s^2+\frac32\right)s
                                                \,e^{-s^2} +
                                                \frac{\sqrt{\pi}}{2}
                                                \,
                                                \textrm{erfc}(s)\right),
          \\ \, \\
	\ds\int_s^\infty z^5 \, e^{-z^2} \, = \, & \ds\left(1+s^2+\frac12
                                                s^4\right) e^{-s^2},\\
          \, \\
	\ds\int_s^\infty z^6 \, e^{-z^2} \, = \, & \ds\frac12 \left(\left(s^4+\frac52 s^2+\frac{15}{4}\right)s\,e^{-s^2} + \frac{15}{4}\,\frac{\sqrt{\pi}}{2} \, \textrm{erfc}(s)\right).
	\end{array} \right.
%\end{subequations*}
$$
For the integral 
$$
%\beq
%\label{eq:A6}
\frac{\rho}{\sqrt{2\pi T}} \int_{v < 0} v^n \, \left(\frac{v-u}{\sqrt{T}}\right)^m e^{-\frac{(v-u)^2}{2T} }\, dv
%\eeq
$$
on the lower half plane, by letting $v =z\, \sqrt{2T} + u$, it can be transformed to
$$
%\beq
%\label{eq:A7}
\frac{\rho}{2}\,\frac{2}{\sqrt{\pi}}\int^{-\frac{u}{\sqrt{2T}}}_{\infty} (z \sqrt{2T}+u)^n \, (\sqrt{2} z)^m \, e^{-z^2}\,dz,
%\eeq
$$
replacing $z$ by $-z$, it becomes to
$$
%\beq
%\label{eq:A8}
\frac{\rho}{2}\,\frac{2}{\sqrt{\pi}}\int_{\frac{u}{\sqrt{2T}}}^\infty (- z \sqrt{2T}+u)^n \, (-\sqrt{2} z)^m \, e^{-z^2}\,dz,
%\eeq
$$
except the coefficient ${\rho}/{2}$, it is in the form of \eqref{eq:A3}.

%%% FJF : a modifier

%%%%%%%%%%%%%%%%%%%%%%%%%%%%%%%%%%
%
%%%%%%%%%%%%%%%%%%%%%%%%%%%%%%%%%%
\section{Numerical examples}
\label{sec5}
\setcounter{equation}{0}
\setcounter{figure}{0}
\setcounter{table}{0}

In this section, we will test the hybrid discontinuous Galerkin scheme
with first order time discretizations for some physical relevant
problems. For simplicity and reduce the high dimension in velocity, we
take the reduced BGK system described in Appendix \ref{sec3} as our kinetic equation. No limiters on the discontinuous Galerkin method are applied. We take a second order discontinuous Galerkin method for 1D problem, while for 2D problems specify it in the examples. We take a cut-off domain with $V_c = 8$ and discretize it with $N_v=32$ uniform points along each direction if not specified.

\subsection{Flow caused by evaporation and condensation.}
\label{ex2}

In this example, we consider a rarefied gas flow caused by evaporation and condensation between two parallel planes with condensed phases \cite{aoki1994gas}. We consider two cases: one is a weak evaporation and condensation with phase conditions to be
\beq
	T_{wl} = 1, \quad p_{wl} = 1, \quad T_{wr} = 1.002, \quad p_{wr} = 1.02,
	\label{eq:evapw}
\eeq 
the other is a strong evaporation and condensation 
\beq
	T_{wl} = 0.5, \quad p_{wl} = 0.01, \quad T_{wr} = 1, \quad p_{wr} = 1.
	\label{eq:evaps}
\eeq 
The boundary conditions are fixed wall temperature and pressure, and we take the initial conditions to be a linear function connecting the wall temperature or pressure. We take a nonuniform mesh with $N_x/4$ cells on the width of $0.05$ at the boundary, while with $N_x/2$ cells in the center on a width of $0.9$. For the hybrid scheme, we force the cells within a width of $0.1$ at the boundary always to be in the kinetic region.
	
\begin{figure}[ht]
\centering
\includegraphics[width=2in]{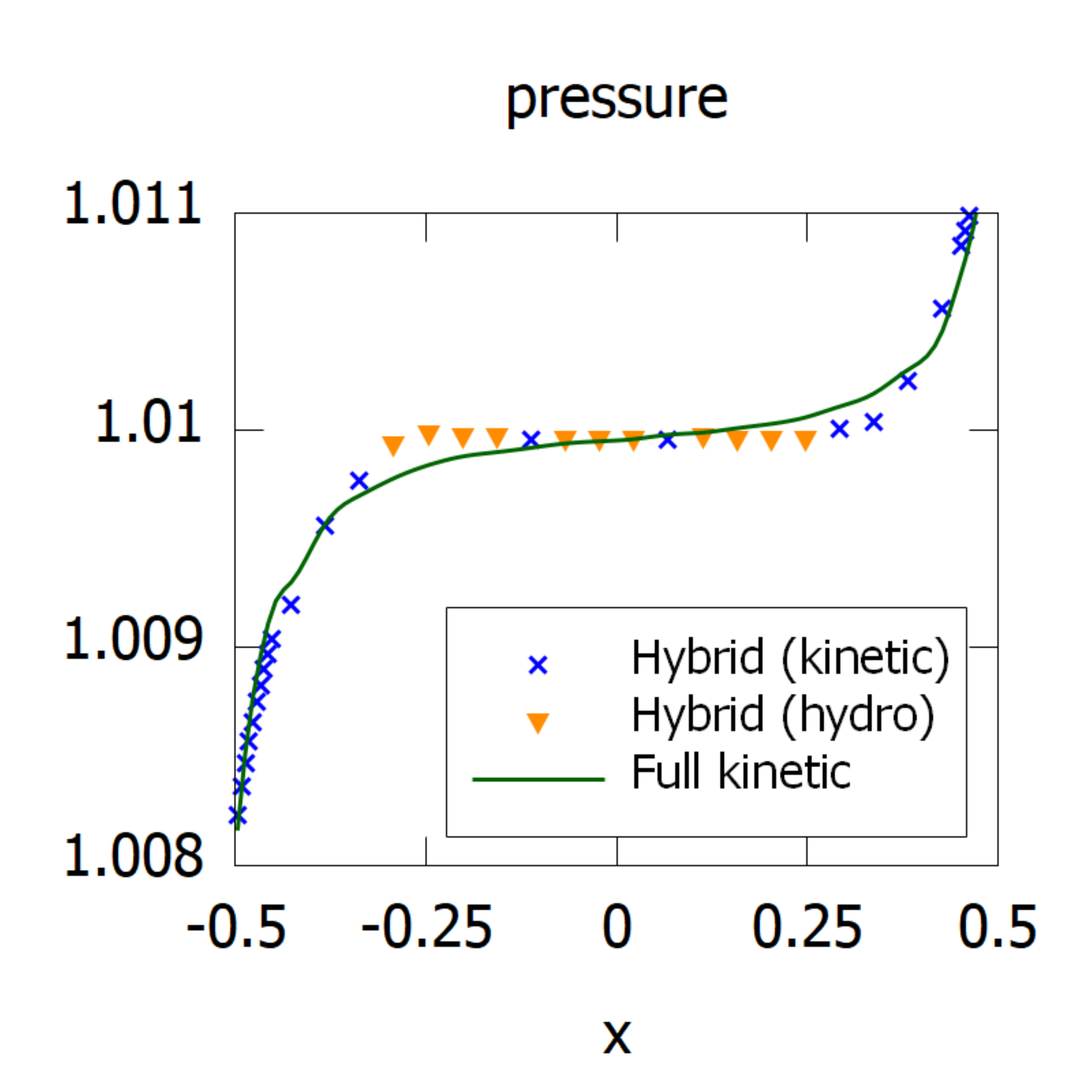}
\includegraphics[width=2in]{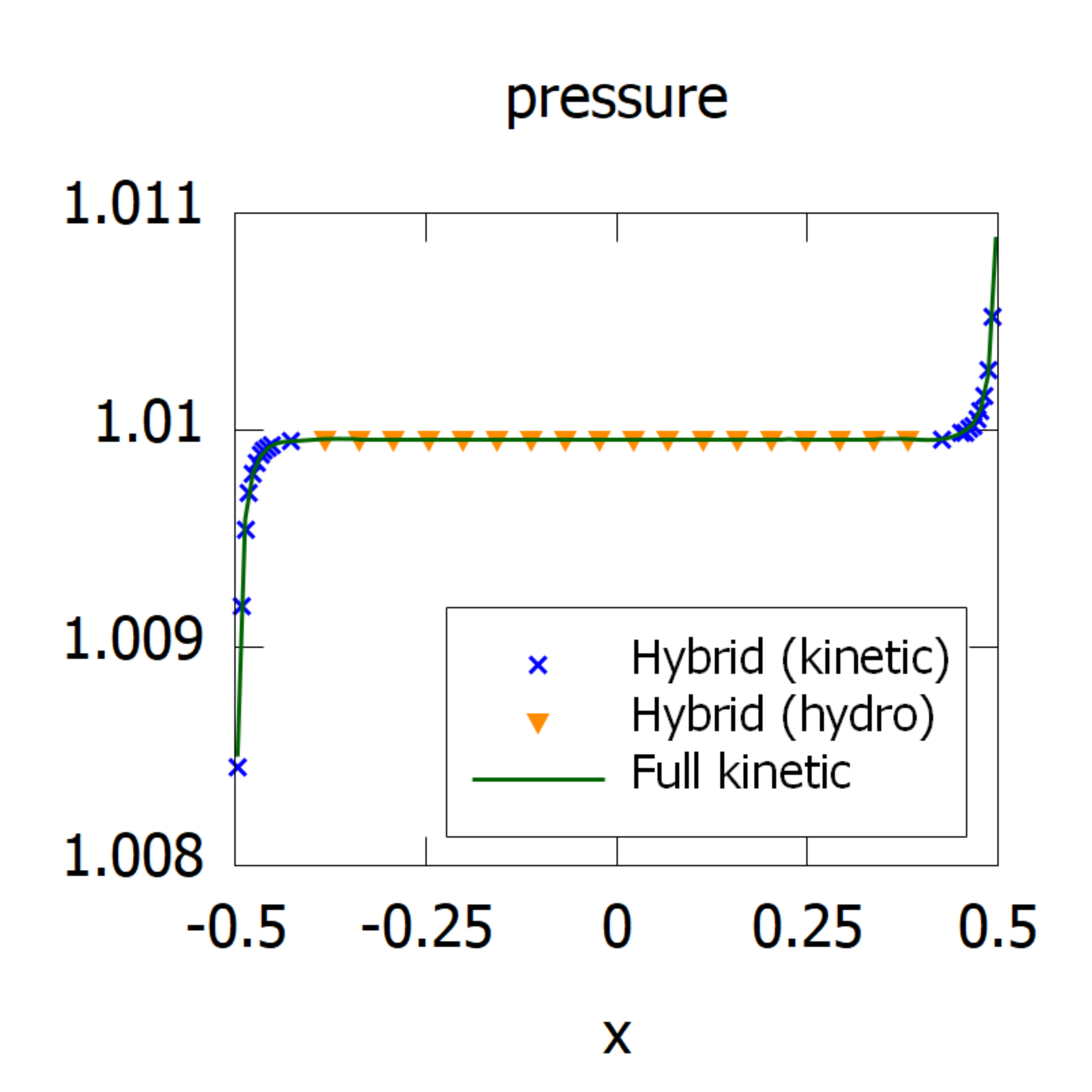}
\includegraphics[width=2in]{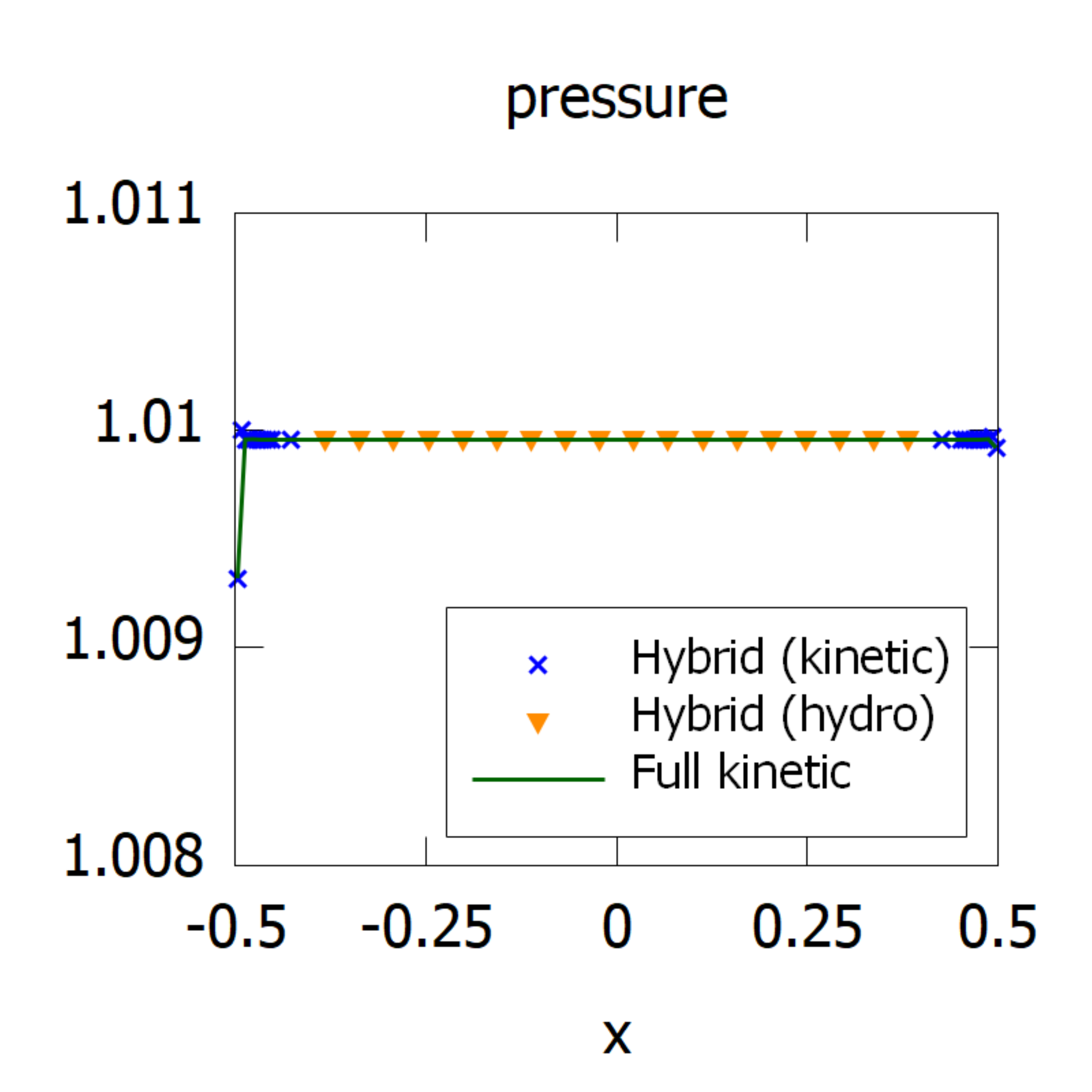}\\
\includegraphics[width=2in]{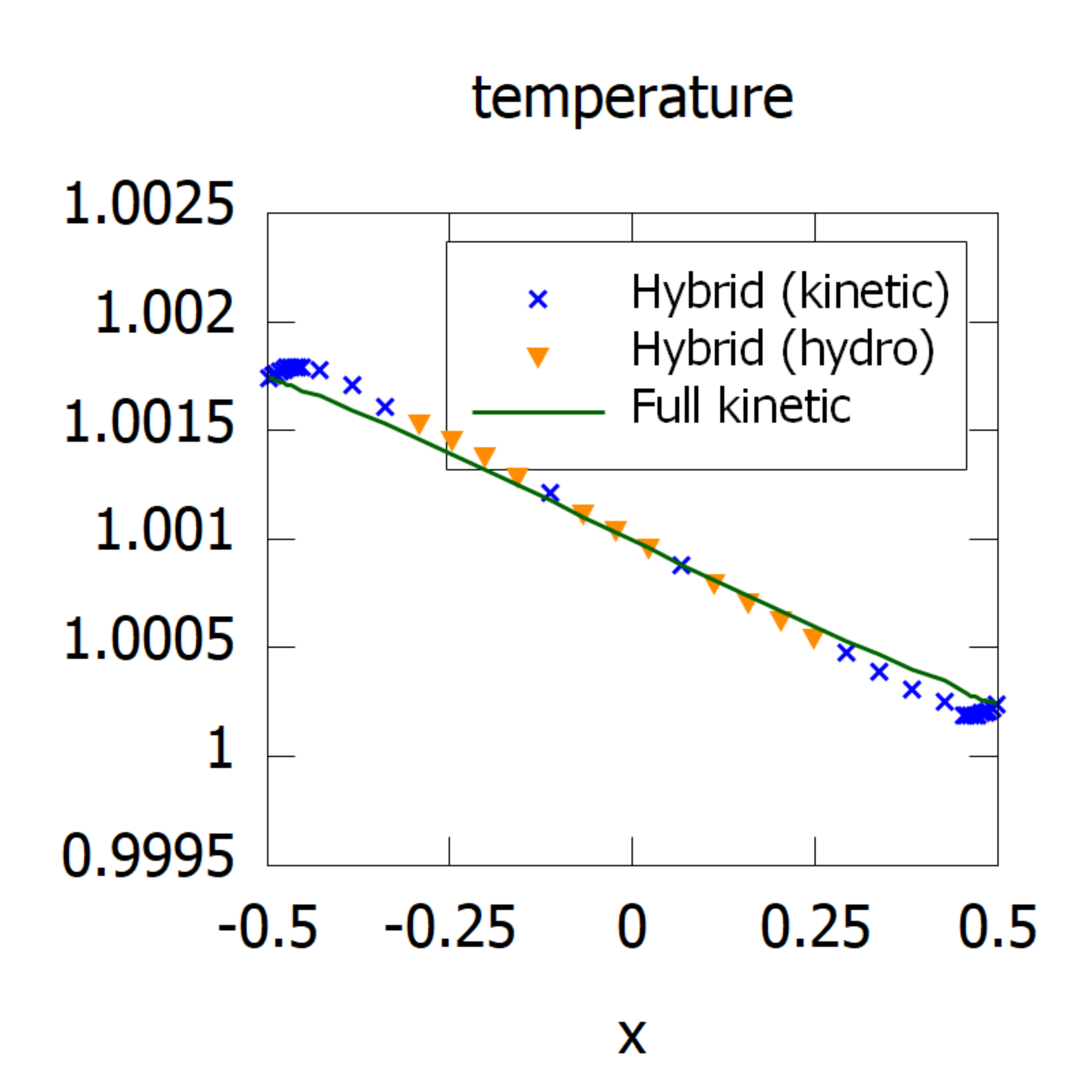}
\includegraphics[width=2in]{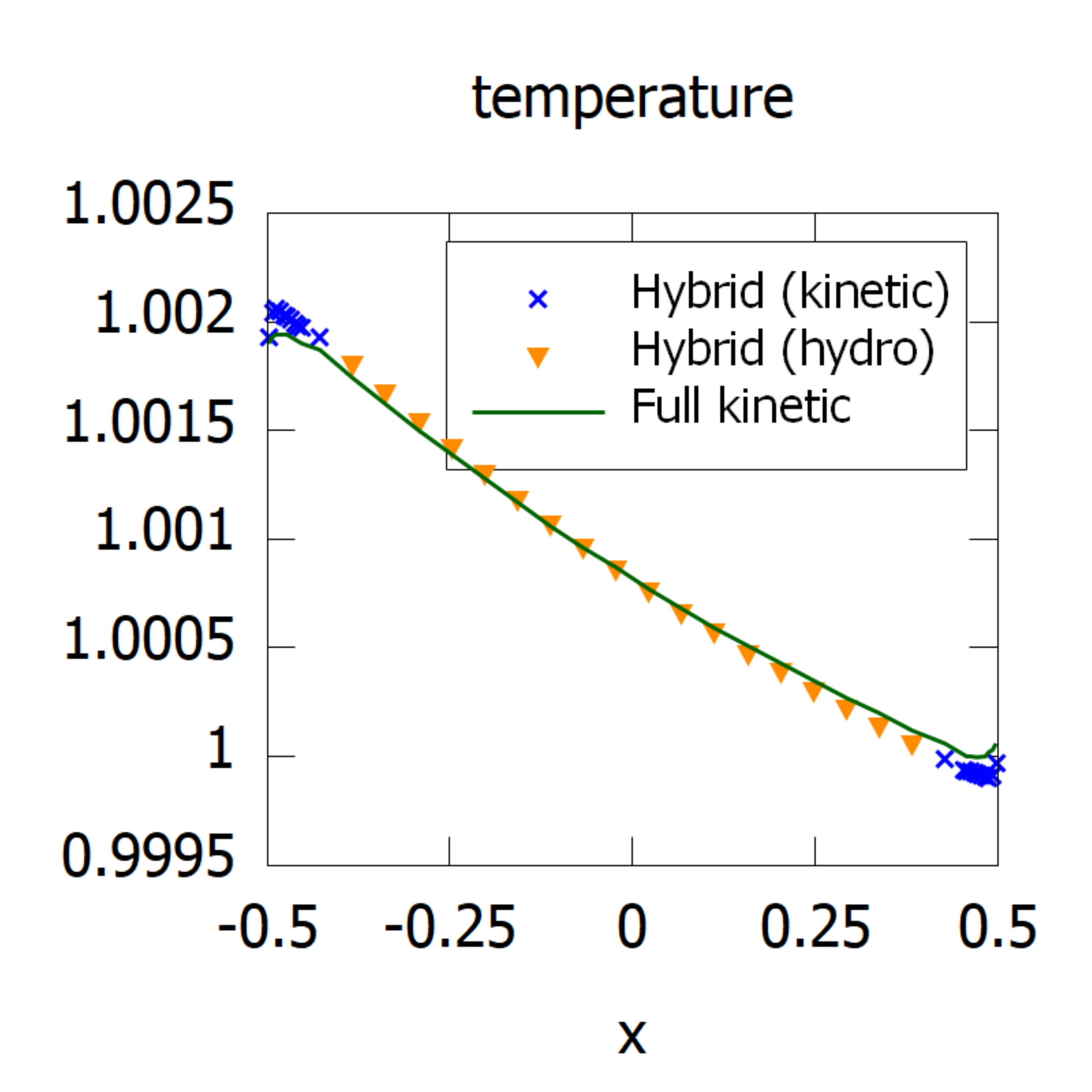}
\includegraphics[width=2in]{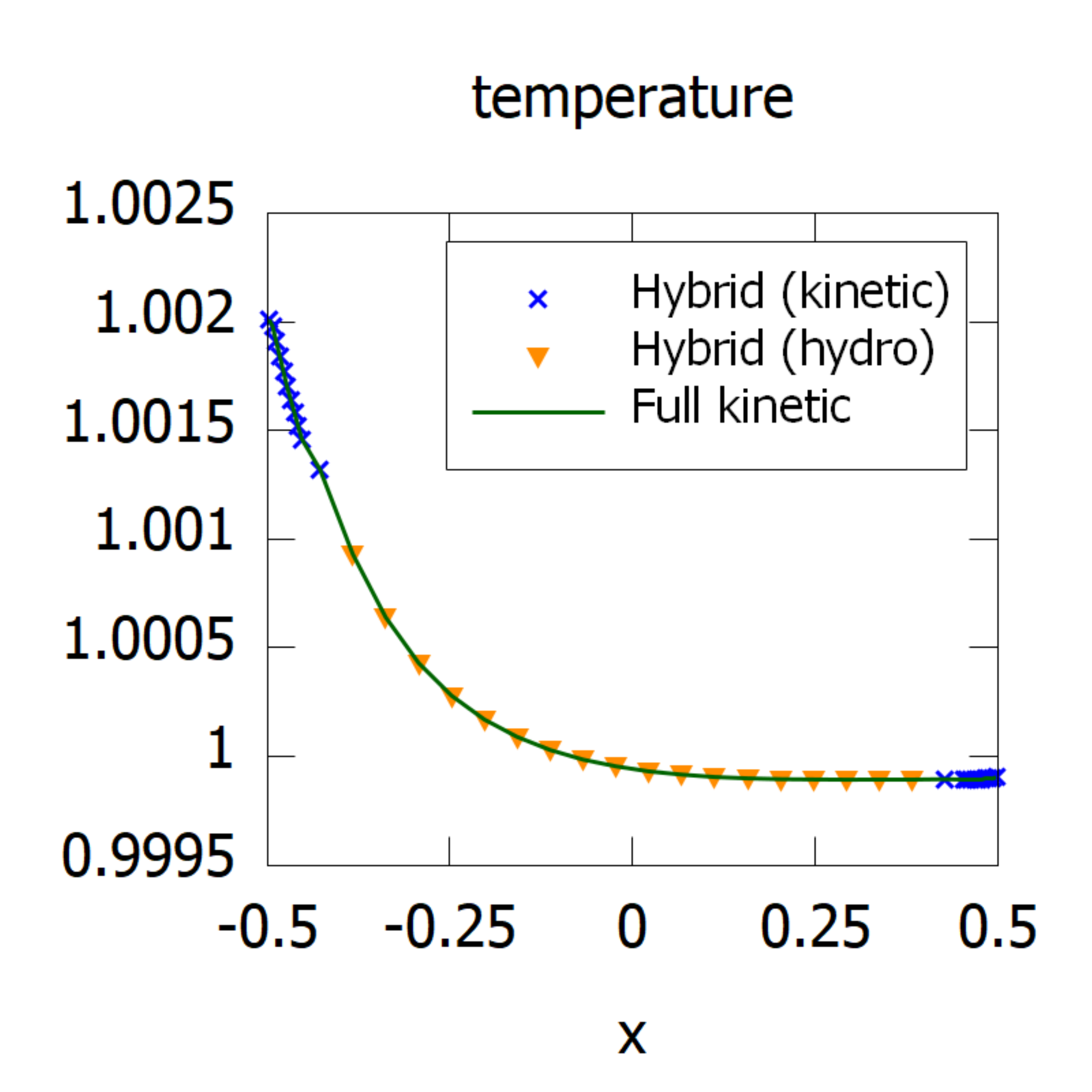}
\caption{{\bf Flow caused by weak evaporation and condensation with
    $p_{wr}/p_{wl}=1.02$, $T_{wr}/T_{wl}=1.002$.} Second order
  discontinuous Galerkin scheme. From left to right:
  $\eps=10^{-1},\,10^{-2},\,10^{-3}$ with a nonuniform mesh with $N_x=40$, $N_x/4$ cells in a width of $0.05$ at the boundary. }
\label{fig21}
\end{figure}

In Figure \ref{fig21}, we show the pressure and temperature with
$N_x=40$ for the weak evaporation and condensation \eqref{eq:evapw},
at $t=100$ for $\eps=10^{-1}, 10^{-2}, 10^{-3}$ respectively. The
results are comparable to those in \cite[Figure 2]{aoki1994gas}
obtained with a much refined mesh $N_x=600$ for $10^{-3} \le \eps
<1$. The mean velocity is nearly a constant up to the spatial
discretization error and is omitted.  For this example, our hybrid
method can well match the full kinetic scheme and identify the flat
part of the steady state solution in the middle to be in the
hydrodynamic region, while the boundary layers are in the kinetic
region. We notice that even for a relatively large $\eps=10^{-1}$, the
hybrid scheme has a result very comparable to the full kinetic
scheme. Furthermore we compare the CPU cost between the hybrid scheme and the full kinetic scheme. We take a uniform mesh with $N_x=40$ and time step $1/5000$. 
We run the code up to $t=40$ for three times. For the full kinetic scheme, the averaged CPU cost
(real time) is $2$ minutes and $45$ seconds, while the hybrid scheme costs about $1$ minute $50$ seconds
for $\eps=10^{-2}$ and $1$ minute $33$ seconds for $\eps=10^{-3}$, so the hybrid scheme saves about $1/3$ and $2/5$ of the CPU cost for $\eps=10^{-2}$ and $\eps=10^{-3}$ respectively. We note that for the nonuniform mesh, since we put half of mesh points around the kinetic boundary, the CPU costs for the hybrid scheme and the full kinetic scheme are almost the same. However, we would emphasize that here we use the reduced BGK system as the kinetic model, the savings of a hybrid scheme would be more significant if we consider the kinetic equation with the full Boltzmann collision operator and full velocity in 3D. Nevertheless the hybrid scheme seems to be only very efficient when the kinetic region (and kinetic cells) is very locally, which may usually happen when $\eps$ is relatively small, e.g., $\eps < 10^{-1}$. 
Besides, for relatively large $\eps$, e.g., $\eps \ge 10^{-1}$, very less hydrodynamic cells can be identified for a hybrid scheme. It will be more efficient to directly use a full kinetic scheme, to avoid the computation of the domain indicator. In the following, we will use $\eps=10^{-1}$ as a reference and only consider our hybrid scheme for $\eps < 10^{-1}$.

\begin{figure}[ht]
\centering		
\includegraphics[width=2in]{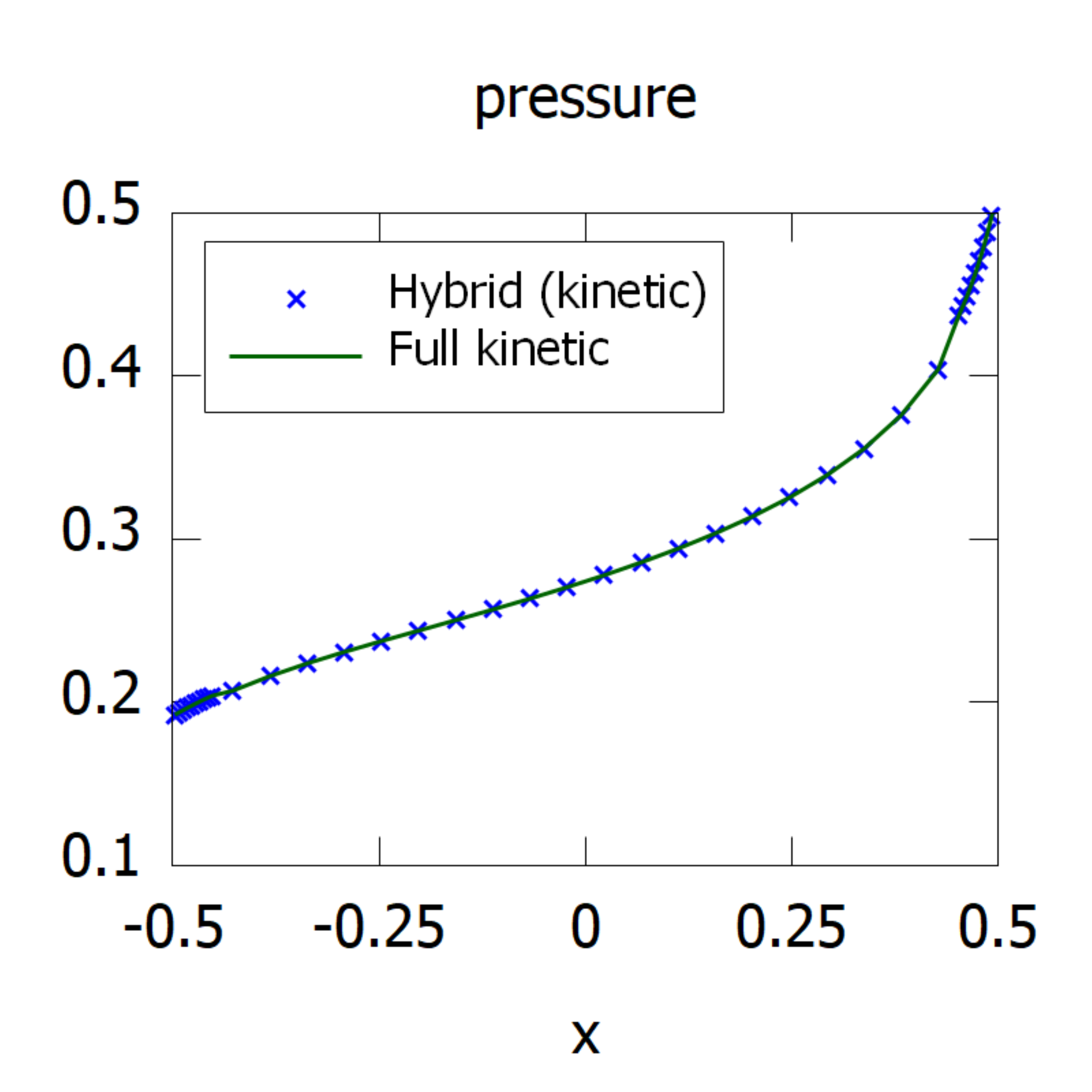}
\includegraphics[width=2in]{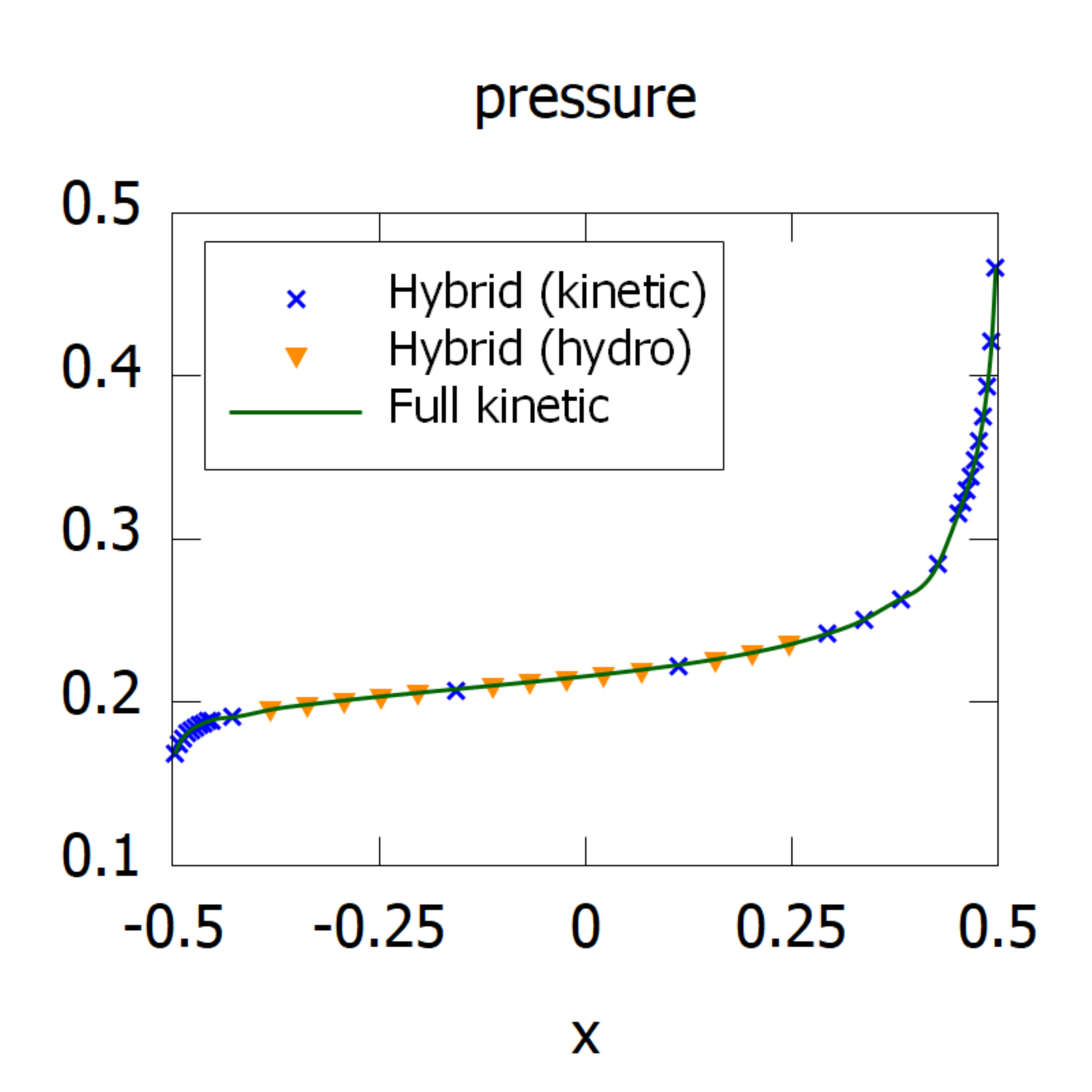}
\includegraphics[width=2in]{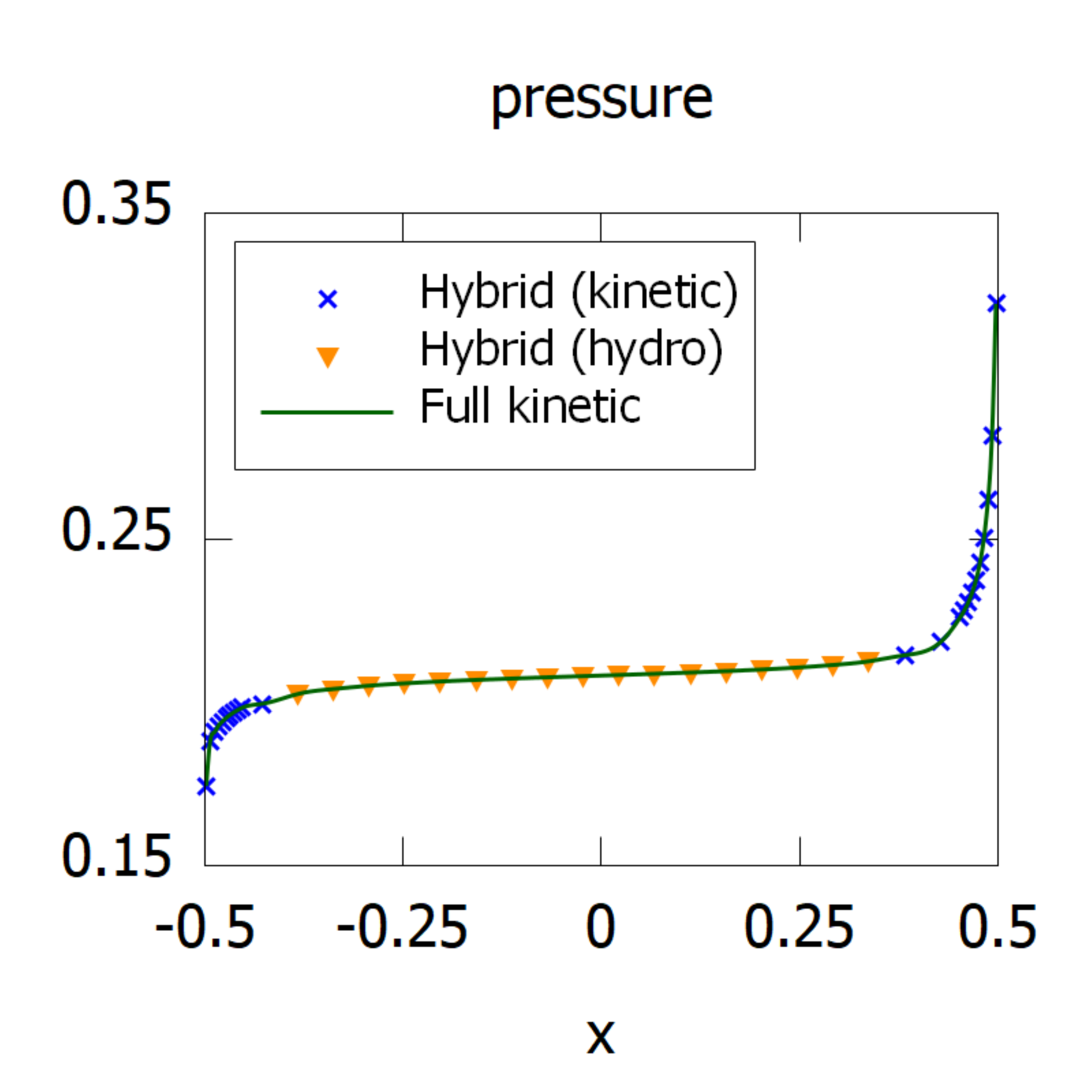}\\
\includegraphics[width=2in]{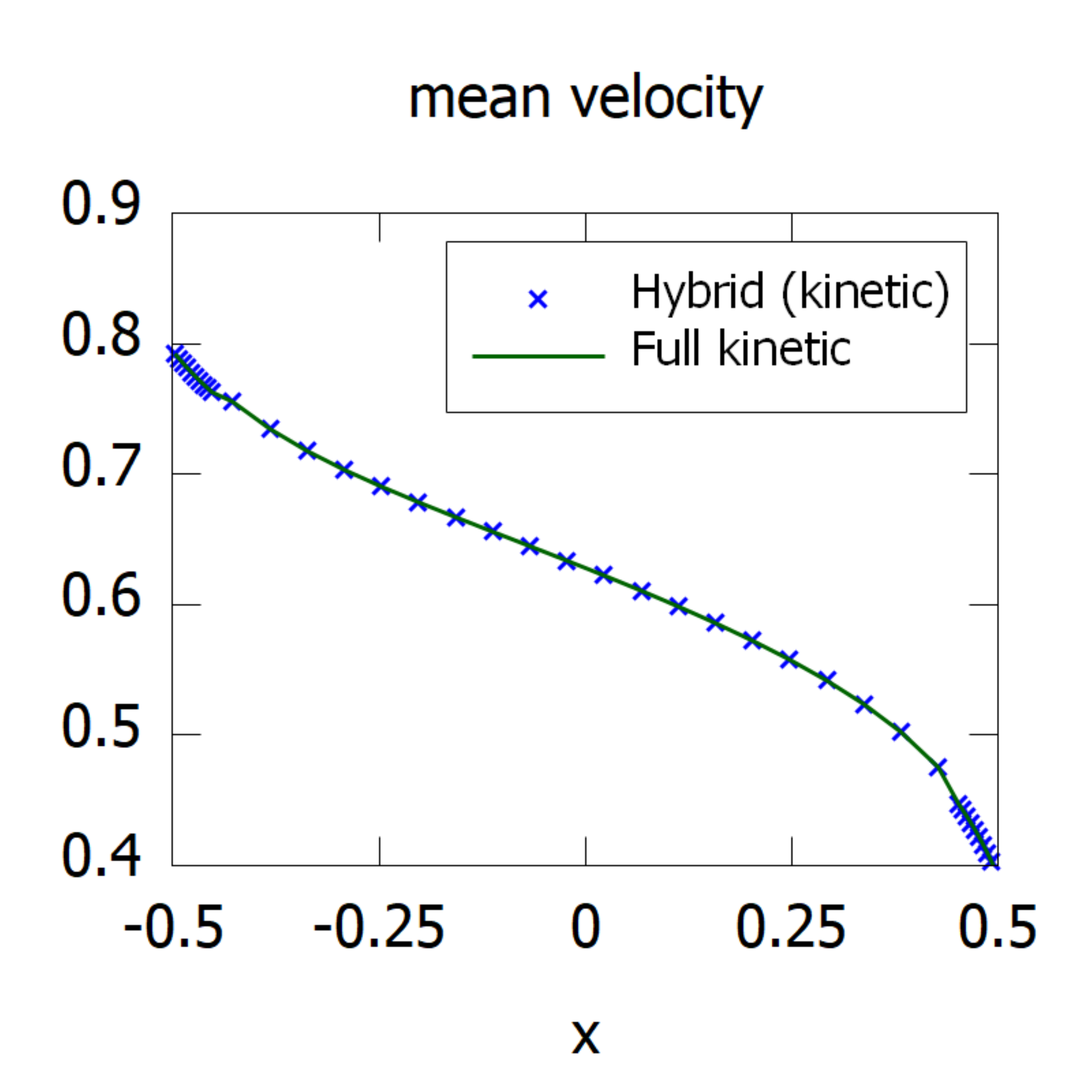}
\includegraphics[width=2in]{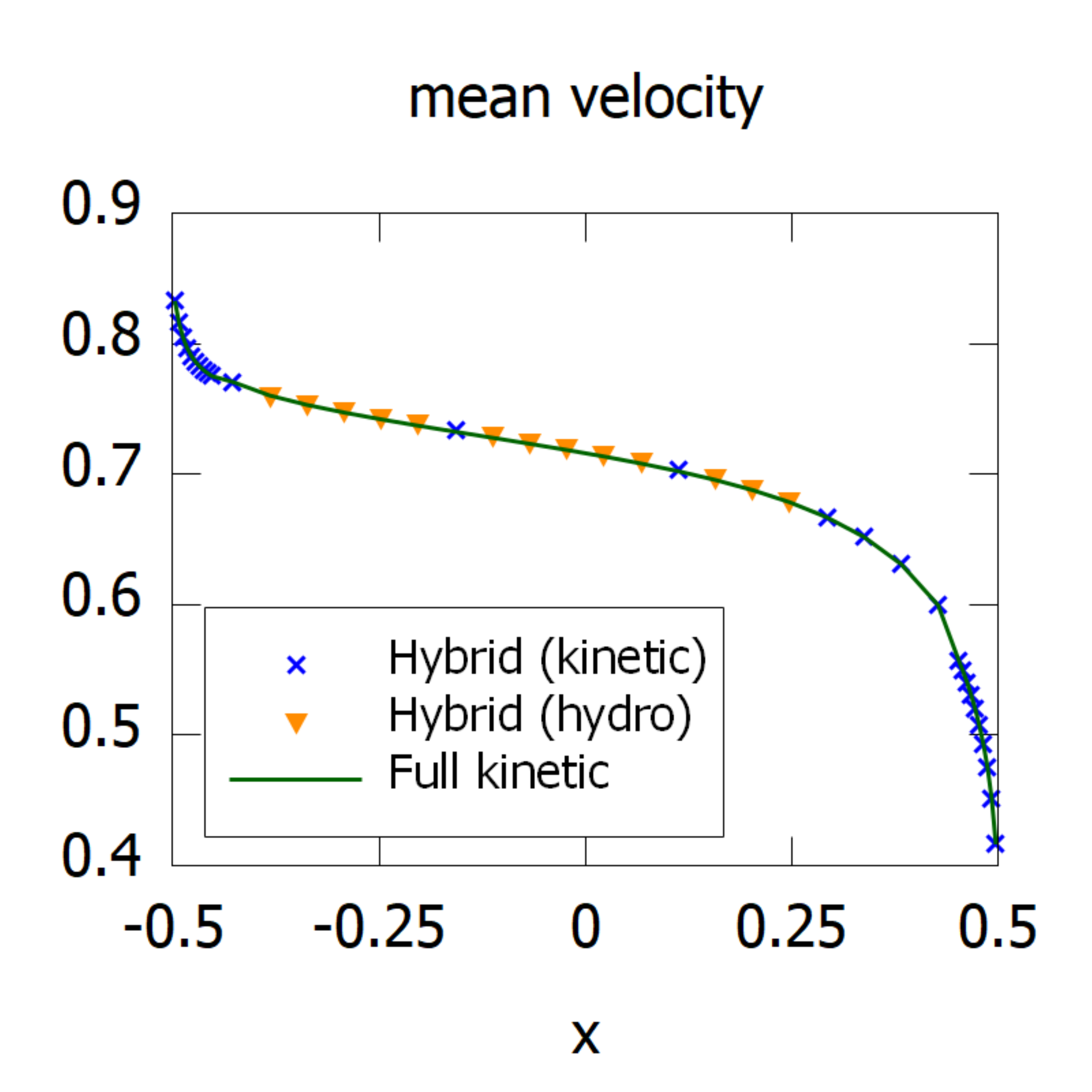}
\includegraphics[width=2in]{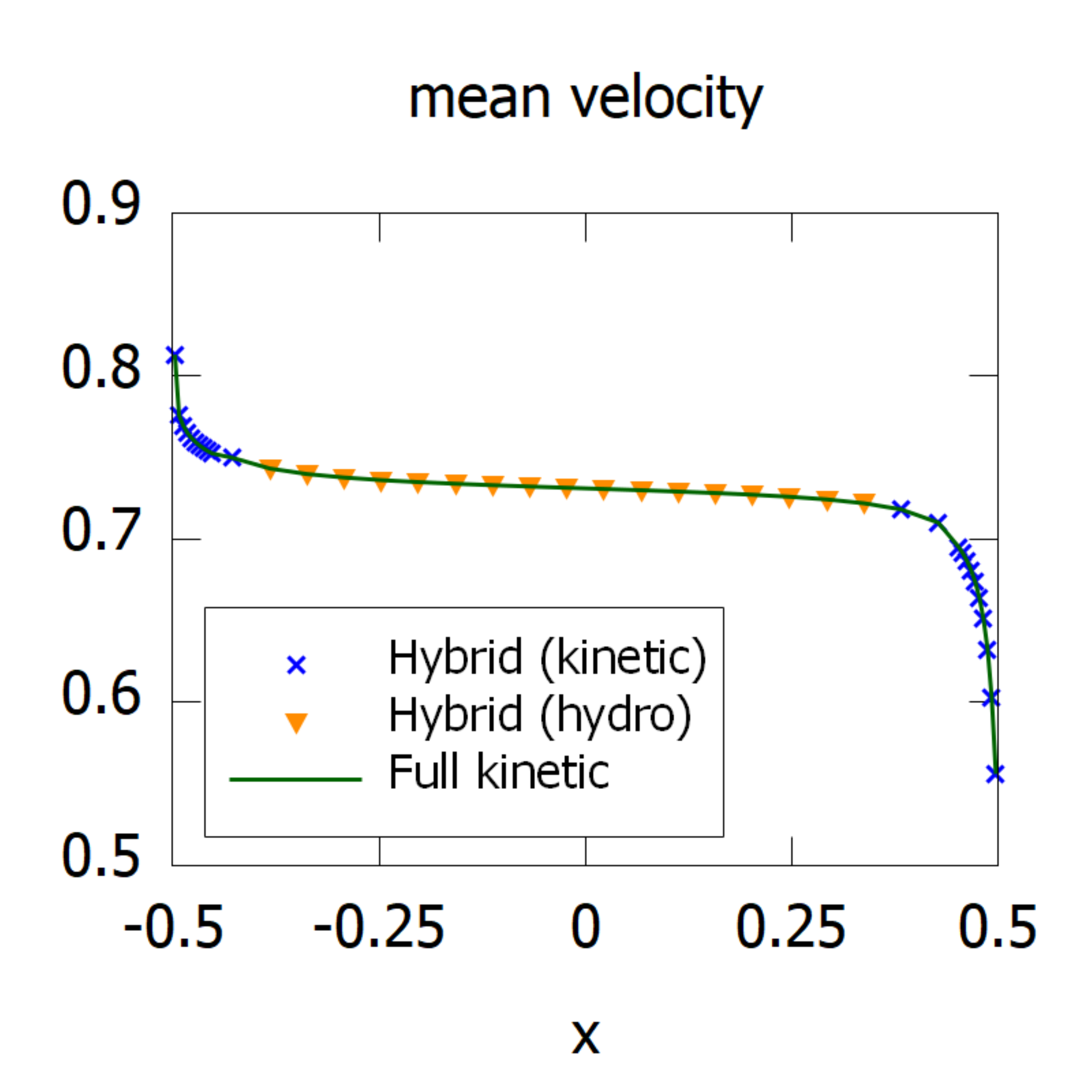}\\
\includegraphics[width=2in]{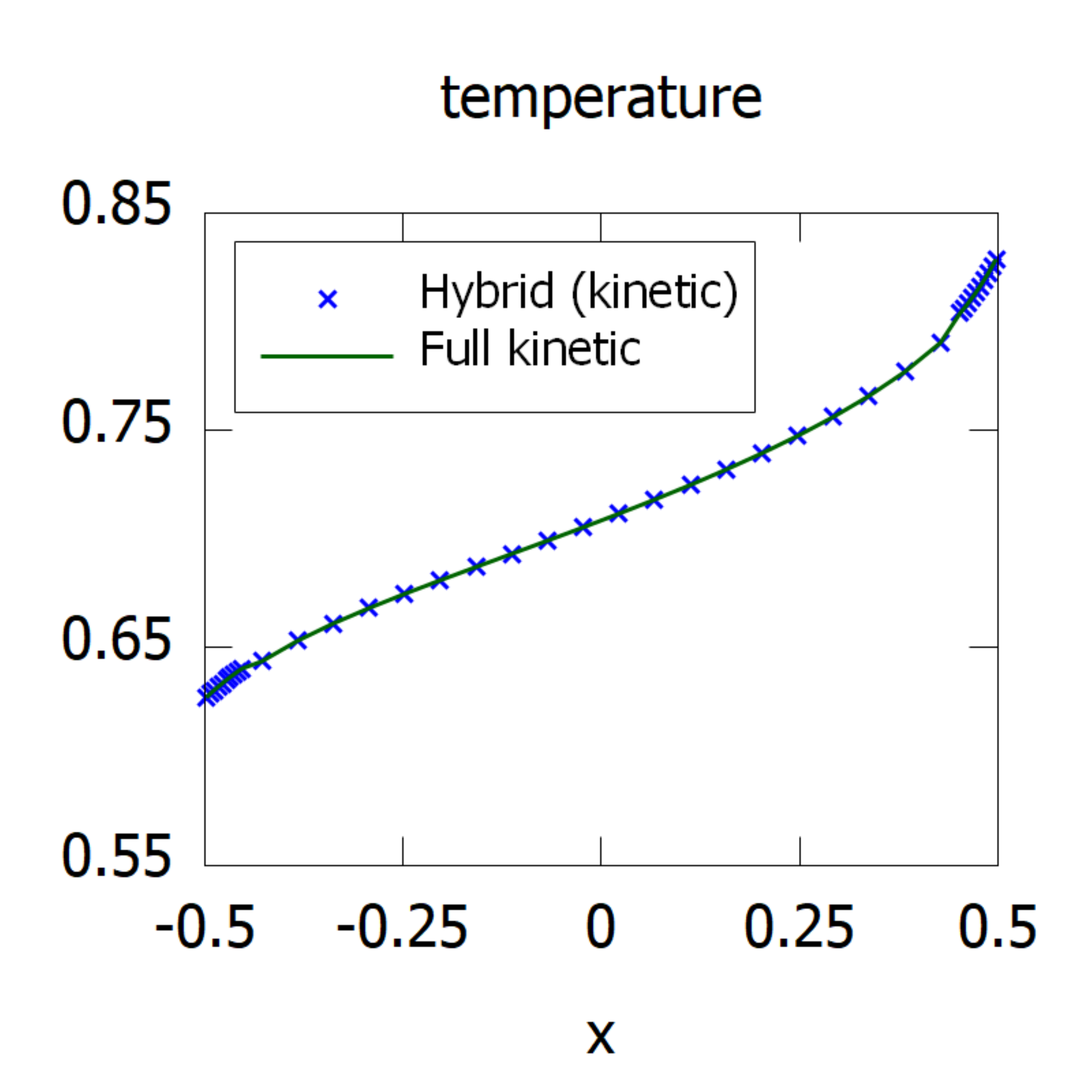}
\includegraphics[width=2in]{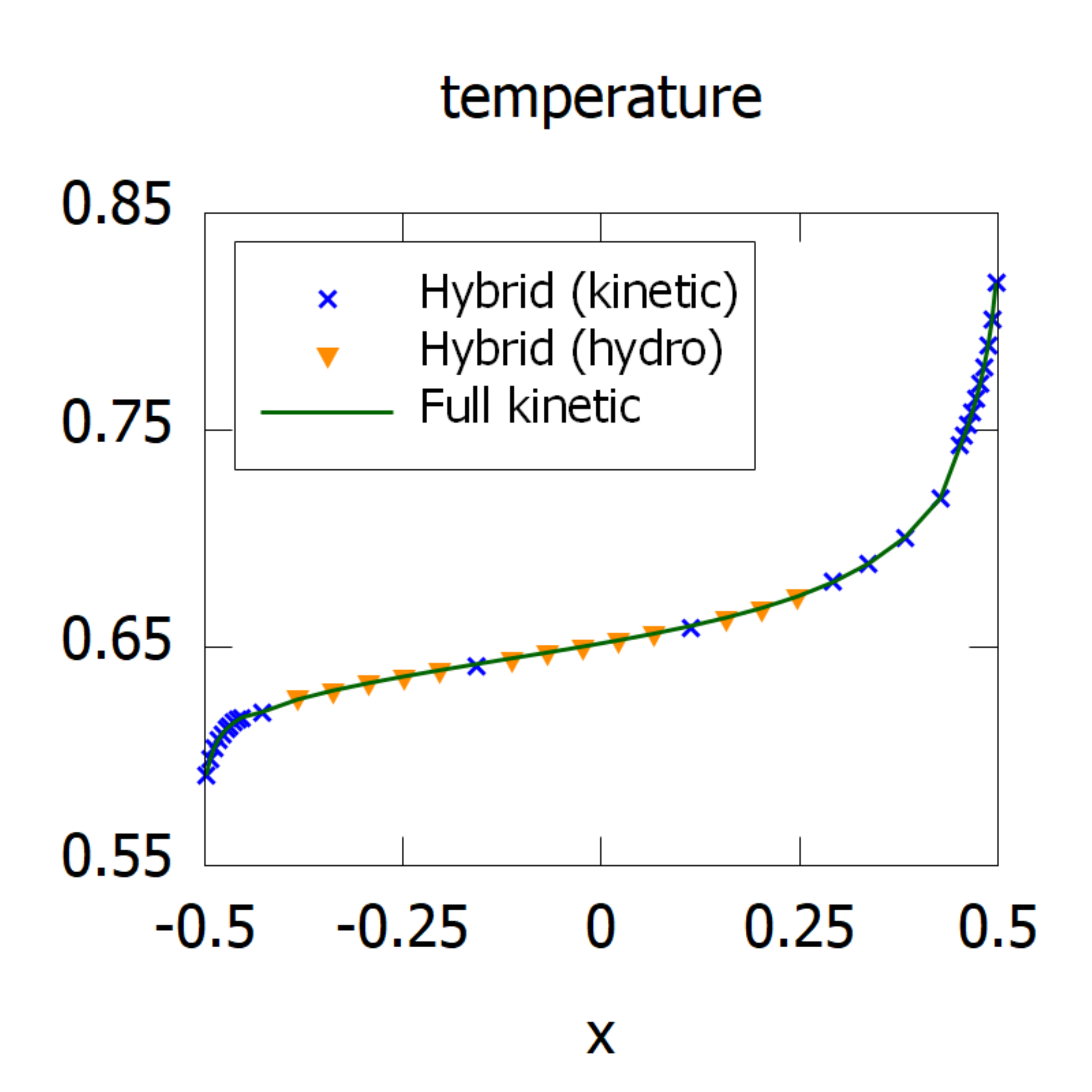}
\includegraphics[width=2in]{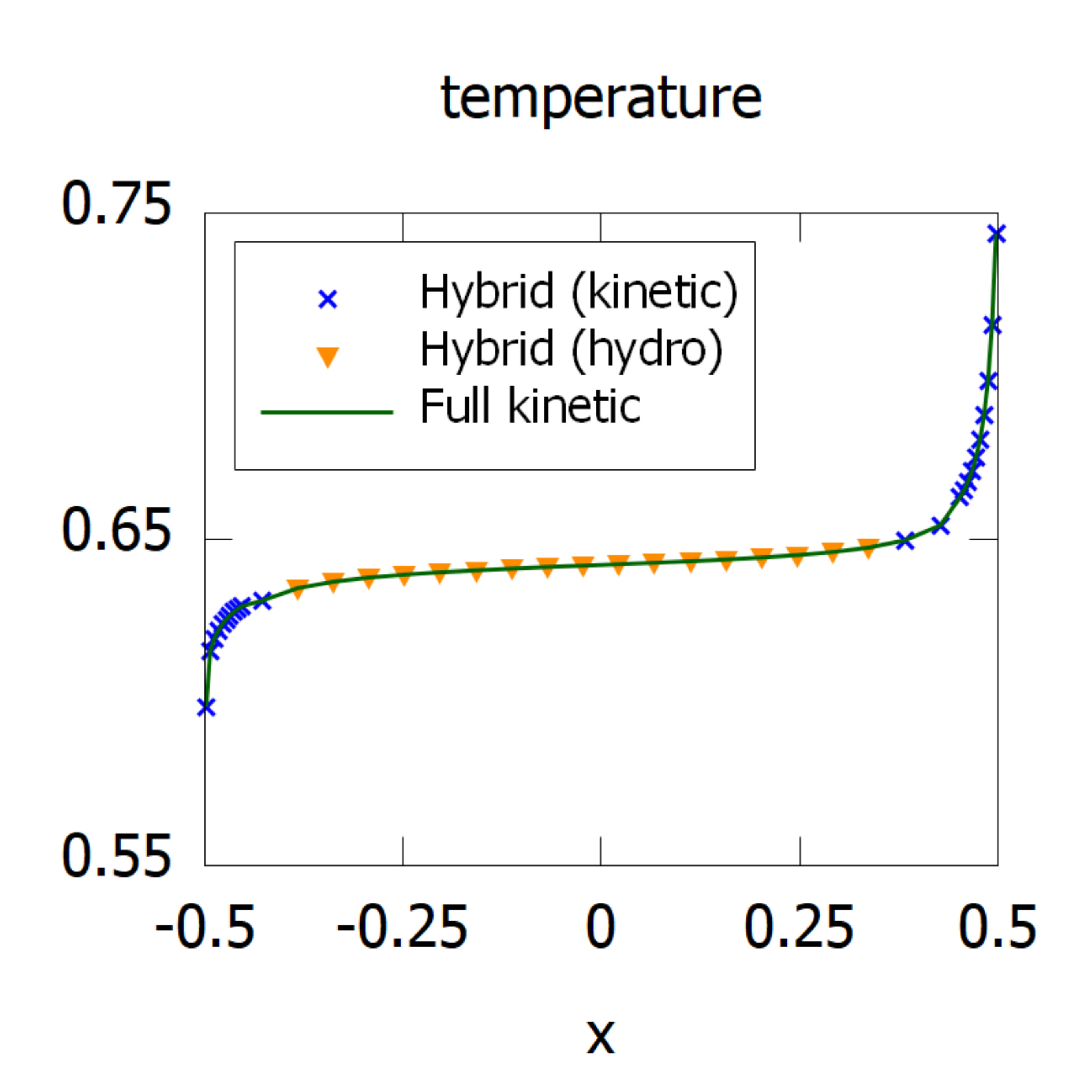}
\caption{{\bf Flow caused by strong evaporation and condensation with
    $p_{wr}/p_{wl}=100$, $T_{wr}/T_{wl}=2$.} Second order
  discontinuous Galerkin scheme.  The mean velocity is scaled by
  ${-1}/{\sqrt{2}}$. From left to right:
  $\eps=10^{-1},\,10^{-2},\,10^{-3}$ with a nonuniform mesh with $N_x=40$, $N_x/4$ cells in a width of $0.05$ at the boundary. }
\label{fig22}
\end{figure}

Next, in Figure \ref{fig22}, we show the pressure, mean velocity and
temperature with $N_x=40$ for the strong evaporation and condensation
\eqref{eq:evaps}, at $t=100$ for $\eps=10^{-1}, 10^{-2}, 10^{-3}$
respectively. The mean velocity is scaled by a factor of
$-{1}/{\sqrt{2}}$. The results are comparable to those in
\cite[Figure 11]{aoki1994gas}  for the results on a much refined
nonuniform mesh. This example is more demanding, as it requires a good
resolution for small values like $\eps=10^{-3}$ in order to well
capture the boundary layer on both sides. We can observe that for
large $\eps=10^{-1}$, the hybrid scheme automatically becomes a full
kinetic scheme, whereas when $\eps$ goes to zero, the numerical scheme
well identifies the hydrodynamic region away from the boundary
layers. In conclusion, the above two cases show the capability of our hybrid scheme, especially the good capturing of
our defined domain decomposition indicator.

\subsection{2D Riemann problem.} 
\label{ex3}
In this example, we now consider a 2D Riemann problem for polytropic gas with initial datum in four quadrants \cite{schulz1993numerical},
 $$
(\rho,p,u,v)(x,y,0) \,=\,
\left\{
\begin{array}{l}
    \ds  (1.5, 1.5, 0, 0), \, \textrm{if \,} x\geq 0 \textrm{ and } y\geq
                         0,   \\ \,\\              
     \ds  (0.6429, 0.3,1.0328, 0), \, \textrm{if \,} x\leq 0 \textrm{ and }
                                 y\geq 0,   \\ \, \\
      \ds (0.1891, 0.0143, 1.0328, 1.0328),  \, \textrm{if \,} x\leq 0 \textrm{ and }
                                 y\leq 0,   \\ \, \\
      \ds (0.6429, 0.3, 0, 1.0328),  \,\textrm{if \,} x\geq 0 \textrm{ and }
                                 y\leq 0, 
\end{array}\right.
$$
here we have $\gamma=5/3$.

To avoid numerical oscillations from high order discontinuous Galerkin
discretization for the shock problem (usually limiters are needed in
order to control numerical oscillations), we take first order discontinuous Galerkin discretizations along both $x$ and $y$ directions. The mesh size is uniform and $N_x=N_y=80$. For the velocity, we take a cut-off domain $\Omega_v = [-8,8]$ with $N_v=64$ along each direction, due to discontinuous macroscopic density and temperature in the Maxwellian distribution function. We run the code with MPI parallelization along $y$ direction with $4$ processors.

\begin{figure}[ht]
	\centering
	\includegraphics[width=2in]{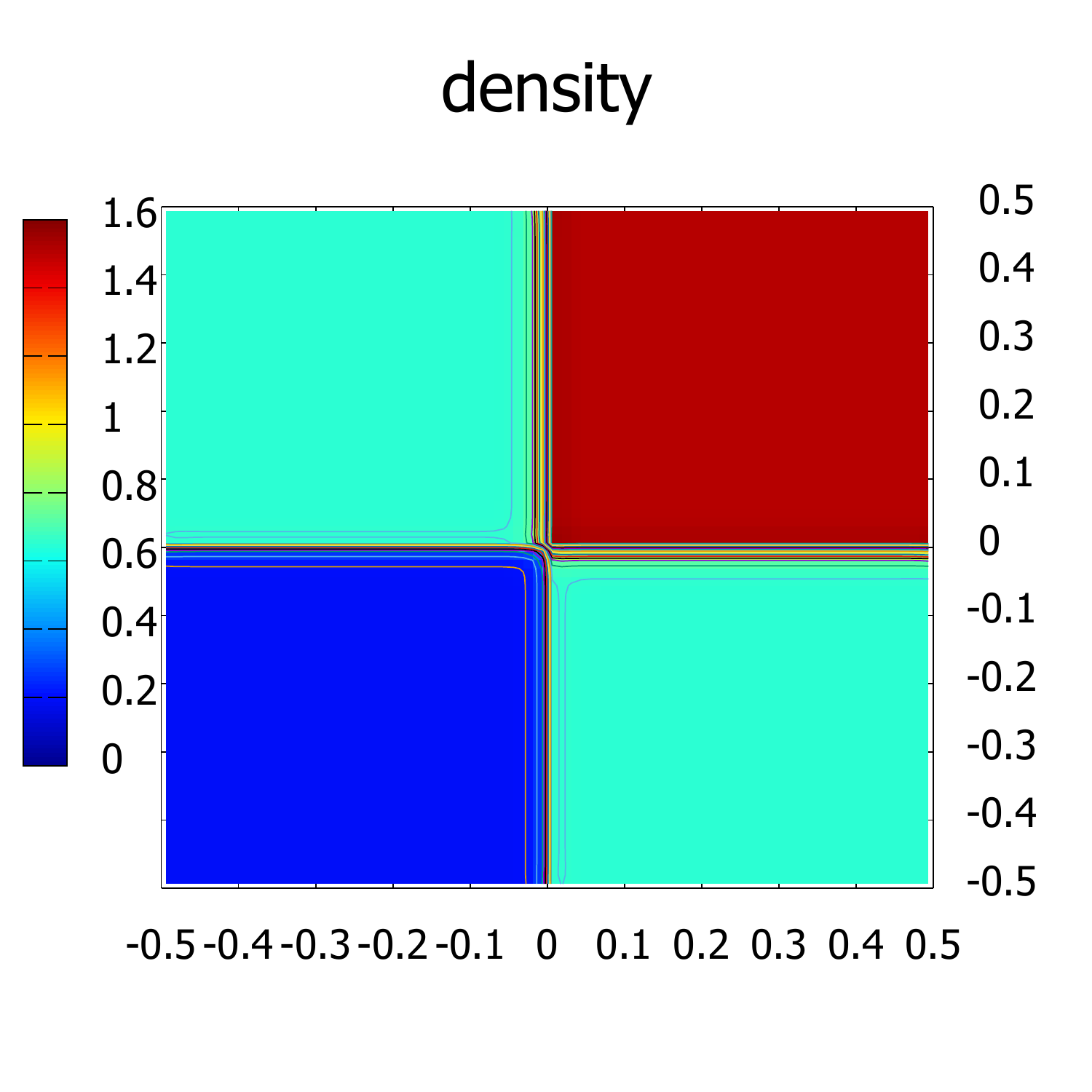}
	\includegraphics[width=2in]{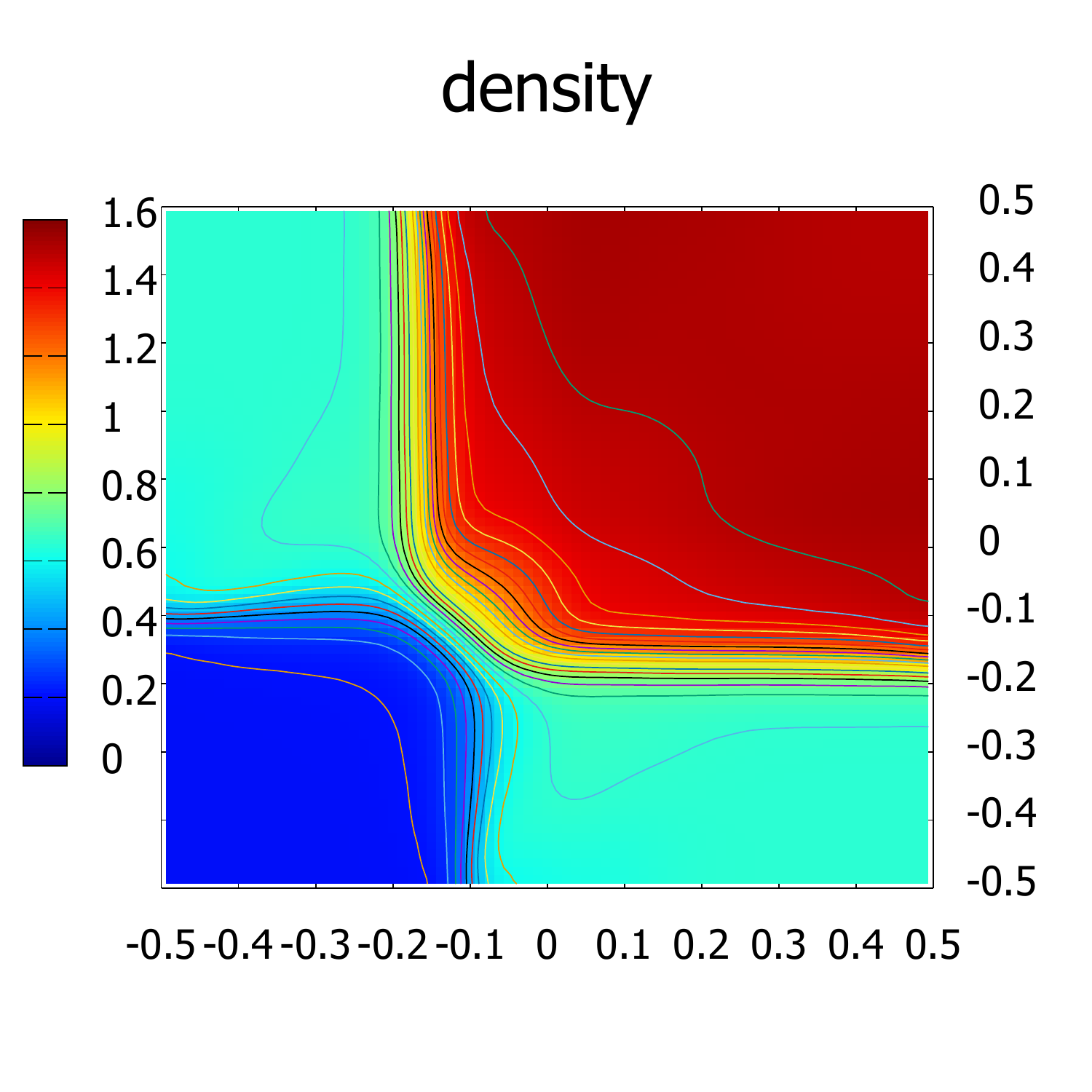}
	\includegraphics[width=2in]{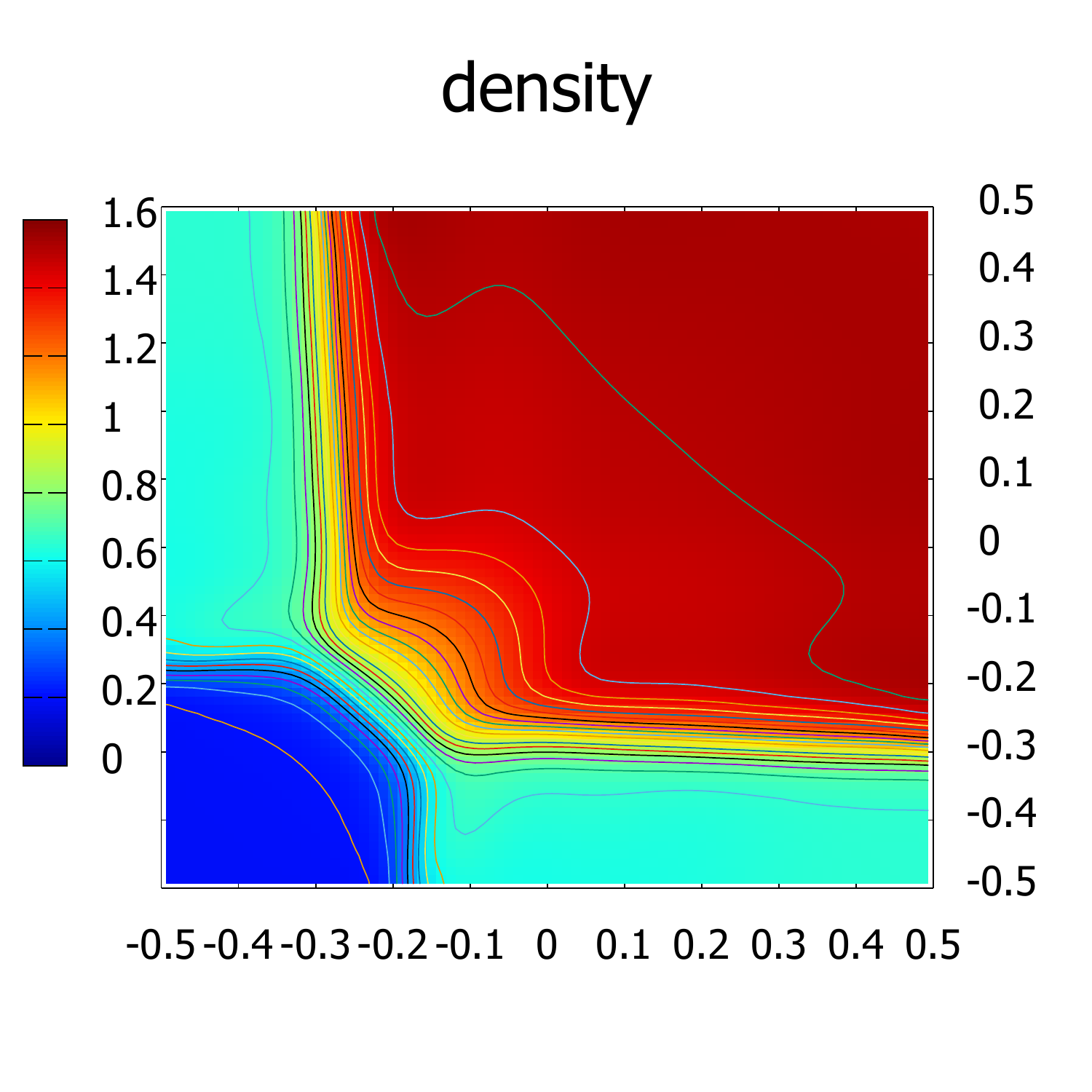}\\
	\includegraphics[width=2in]{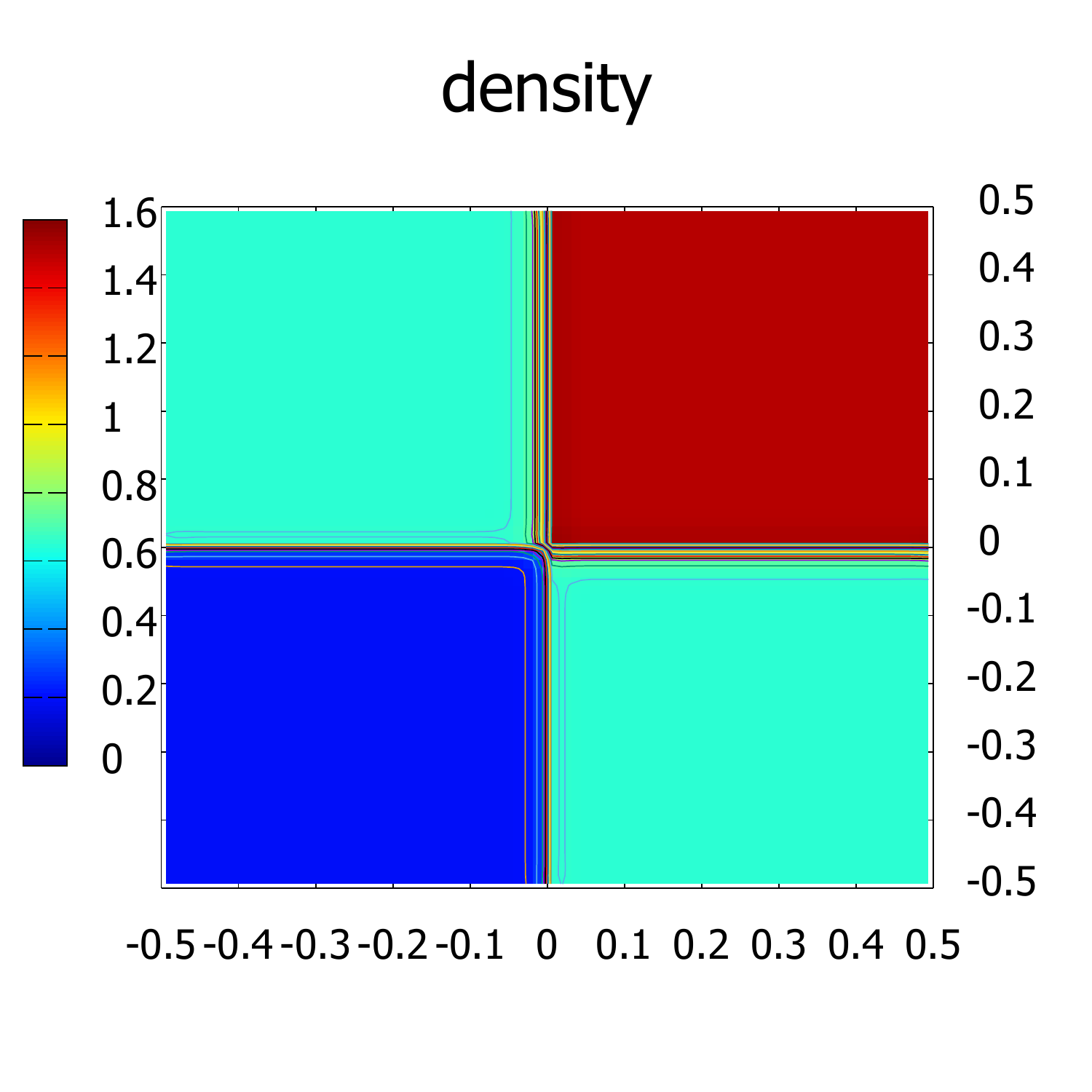}
	\includegraphics[width=2in]{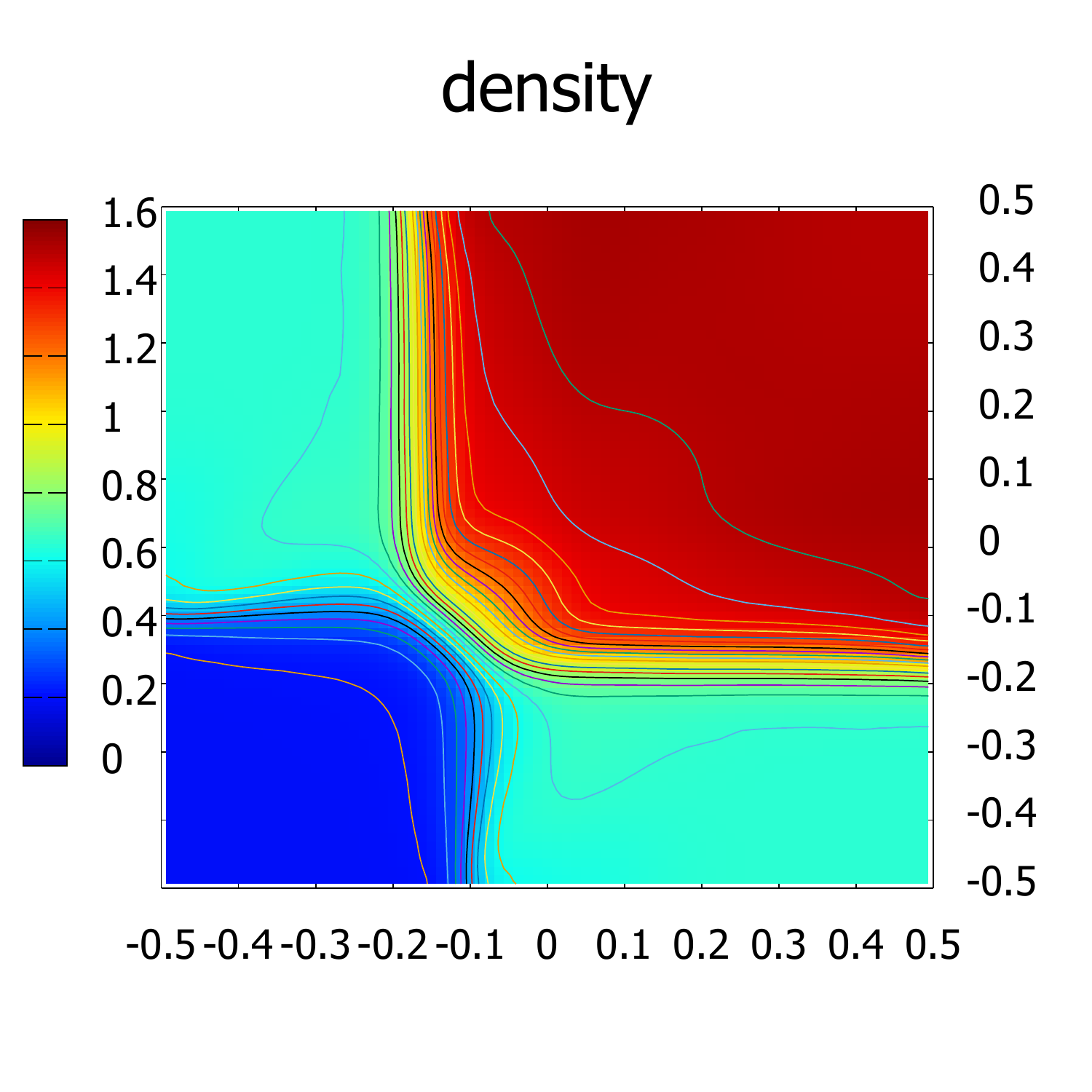}
	\includegraphics[width=2in]{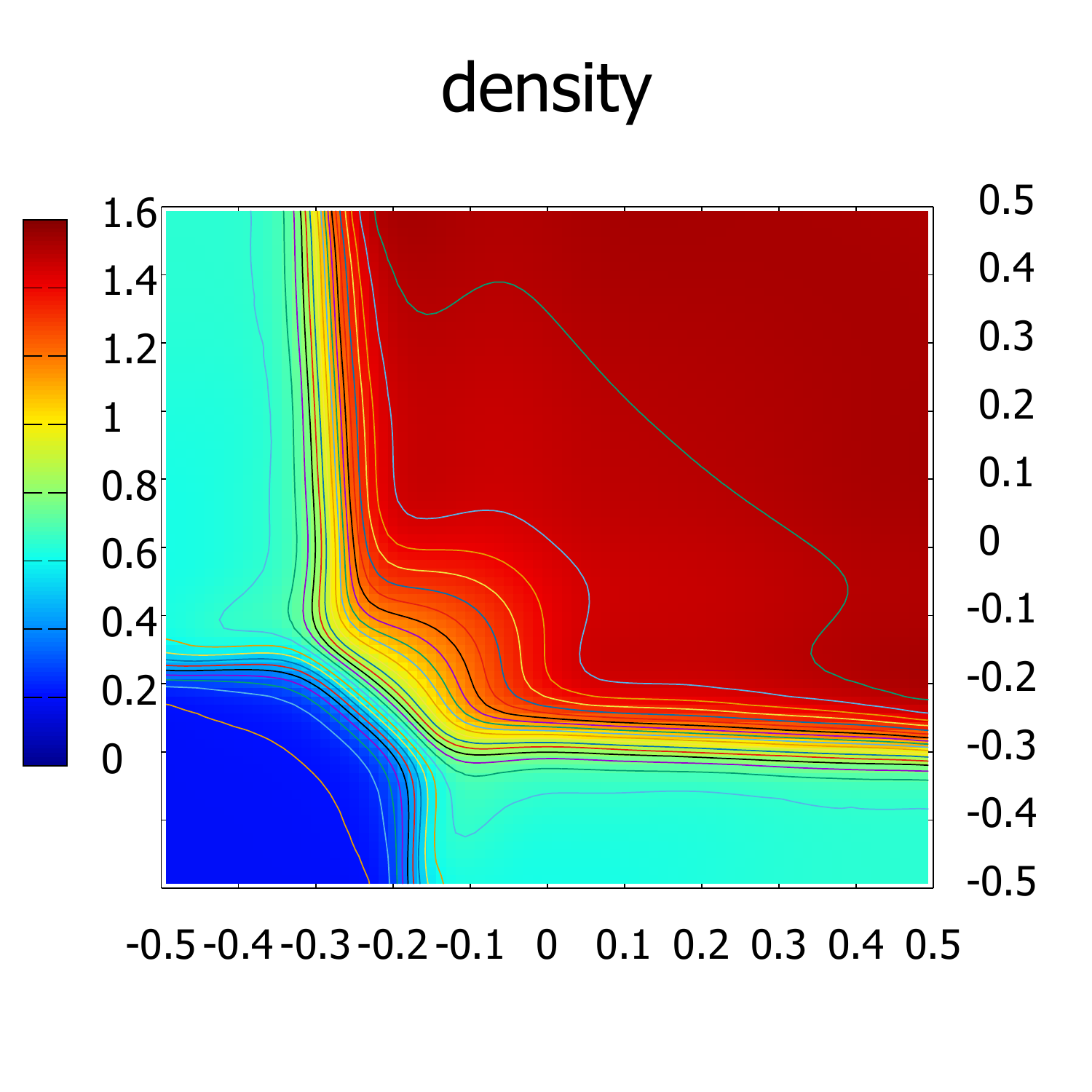}\\
	\includegraphics[width=2in]{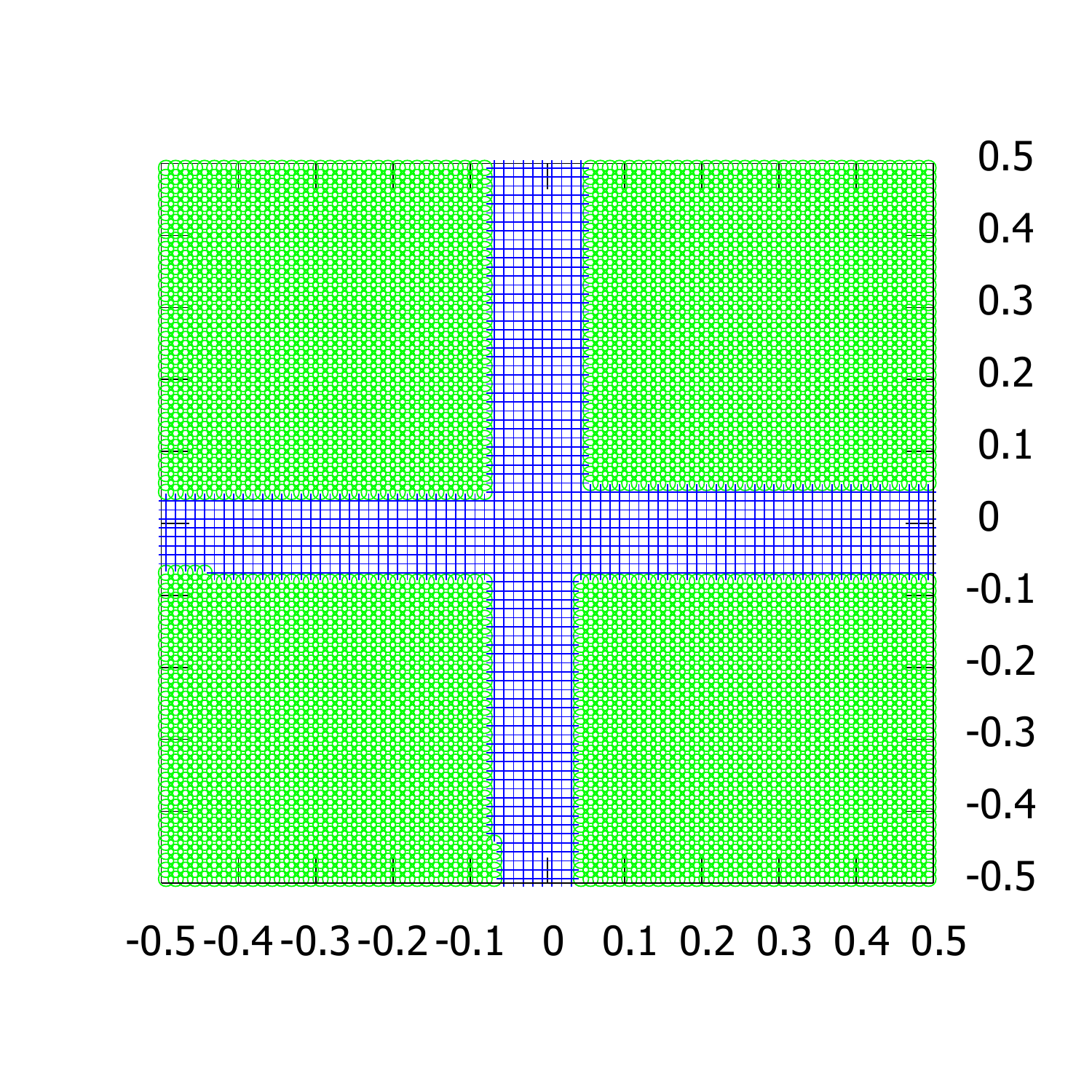}
	\includegraphics[width=2in]{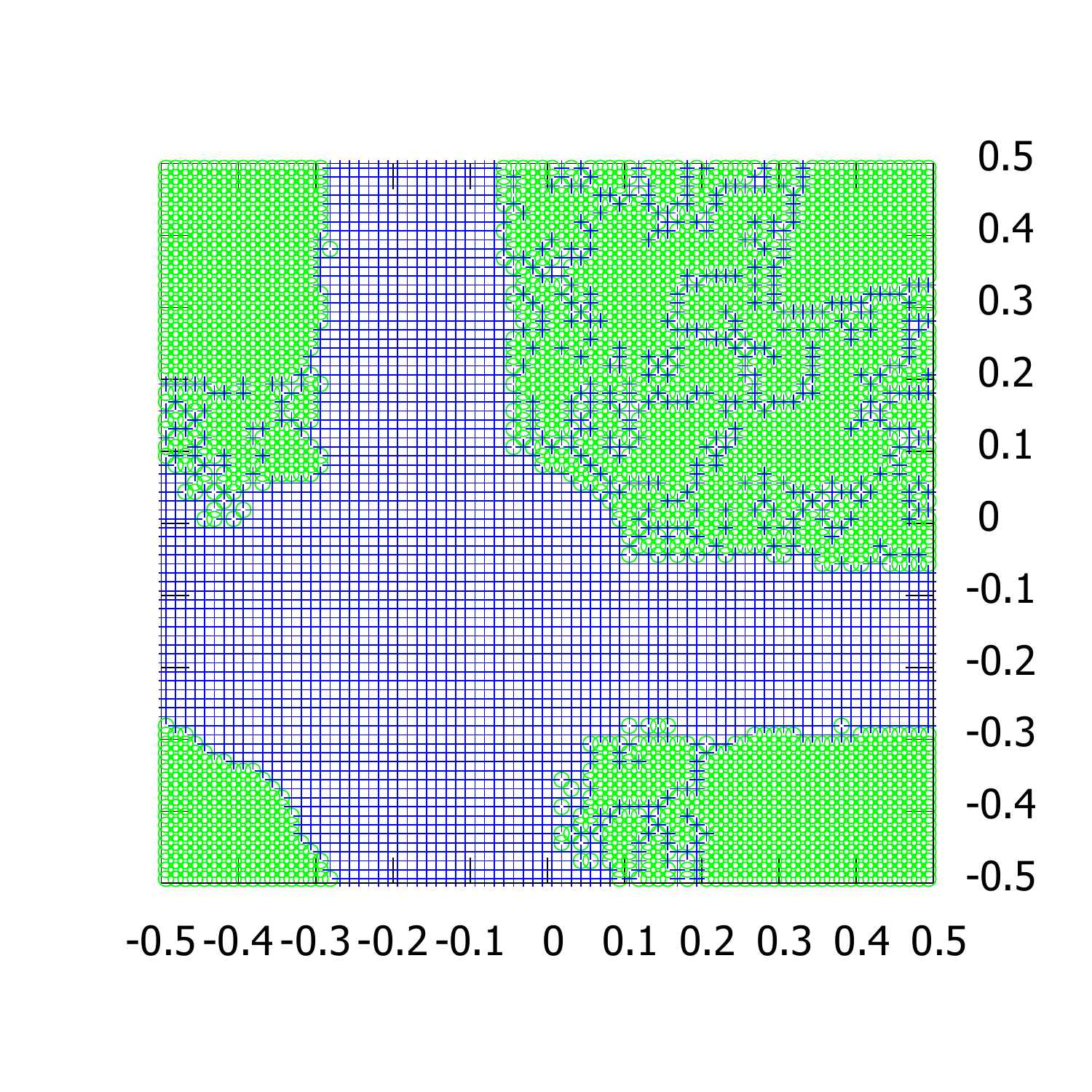}
	\includegraphics[width=2in]{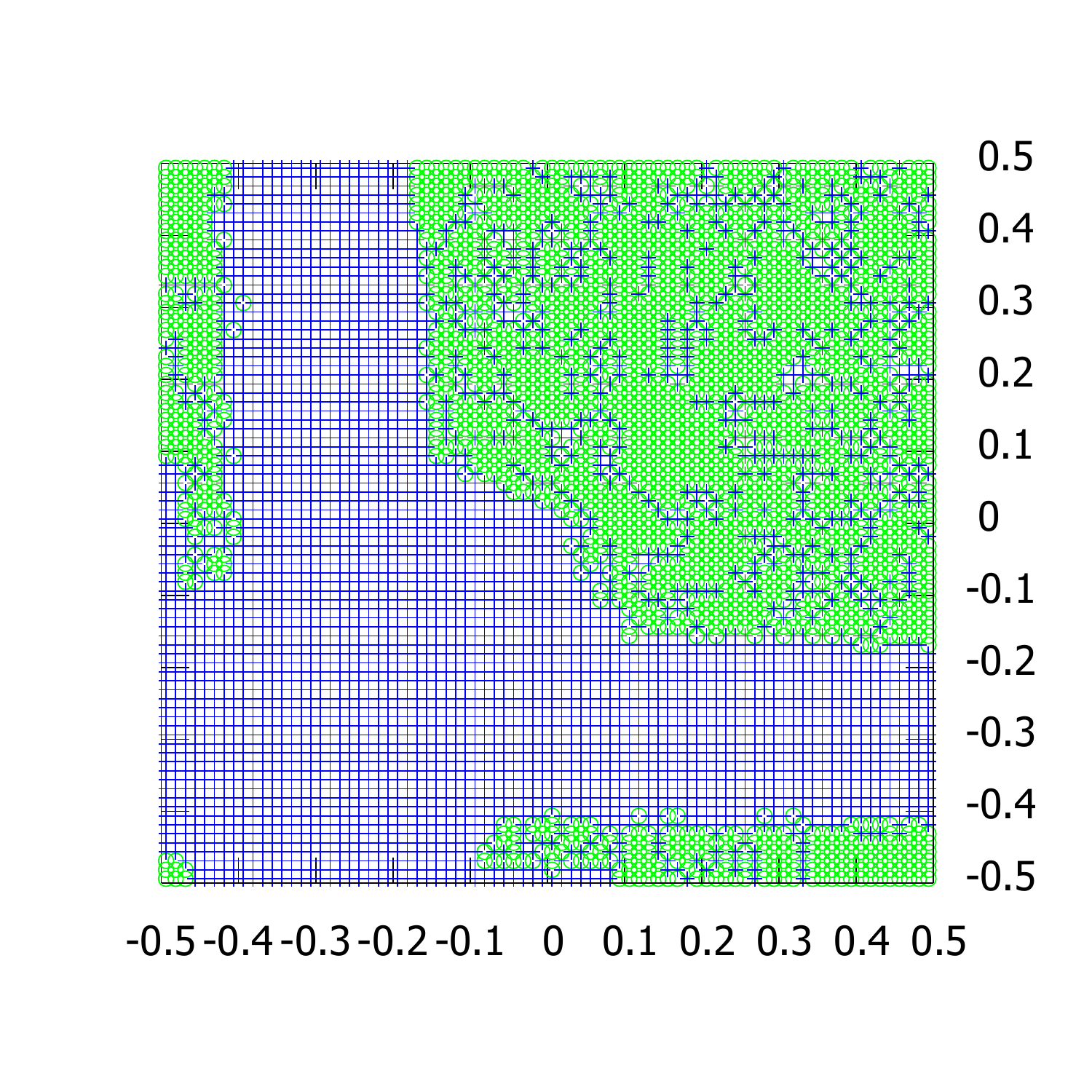}
	\caption{{\bf 2D Riemann problem.} Density profile obtained
          with  a first order discontinuous Galerkin scheme using
          uniform grids with $N_x=N_y=80$ and $\eps=10^{-2}$. From left to right: $t=0.01, 0.2, 0.35$. From top to bottom: the full kinetic scheme, the hybrid scheme, the domain indicator for the hybrid scheme. In the domain indicator, symbol ``+'' denotes kinetic cells, symbol ``o'' denotes hydrodynamic cells. 29 contour lines on the range $[0,1.6]$. }
	\label{fig33}
\end{figure}

\begin{figure}[ht]
	\centering
	\includegraphics[width=2in]{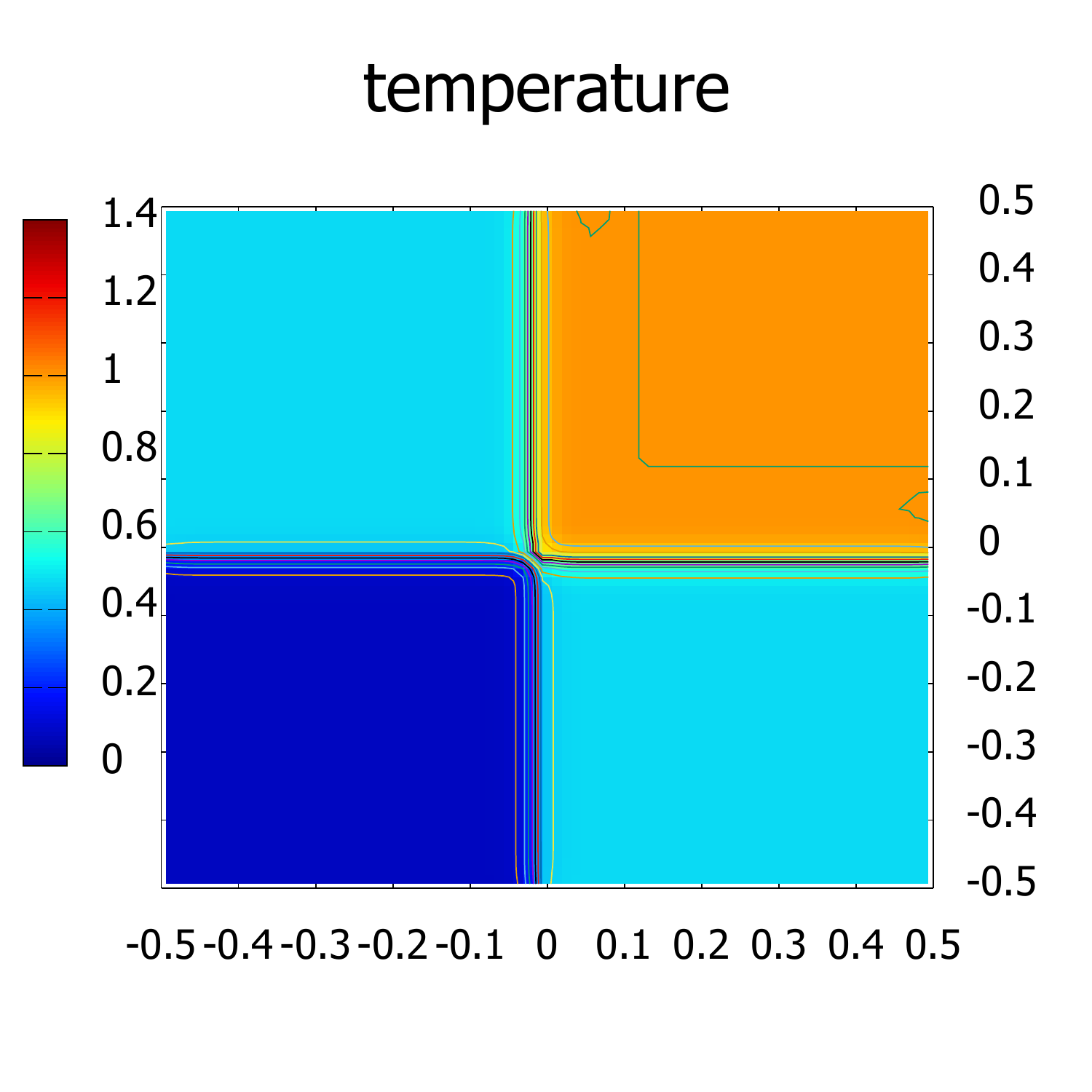}
	\includegraphics[width=2in]{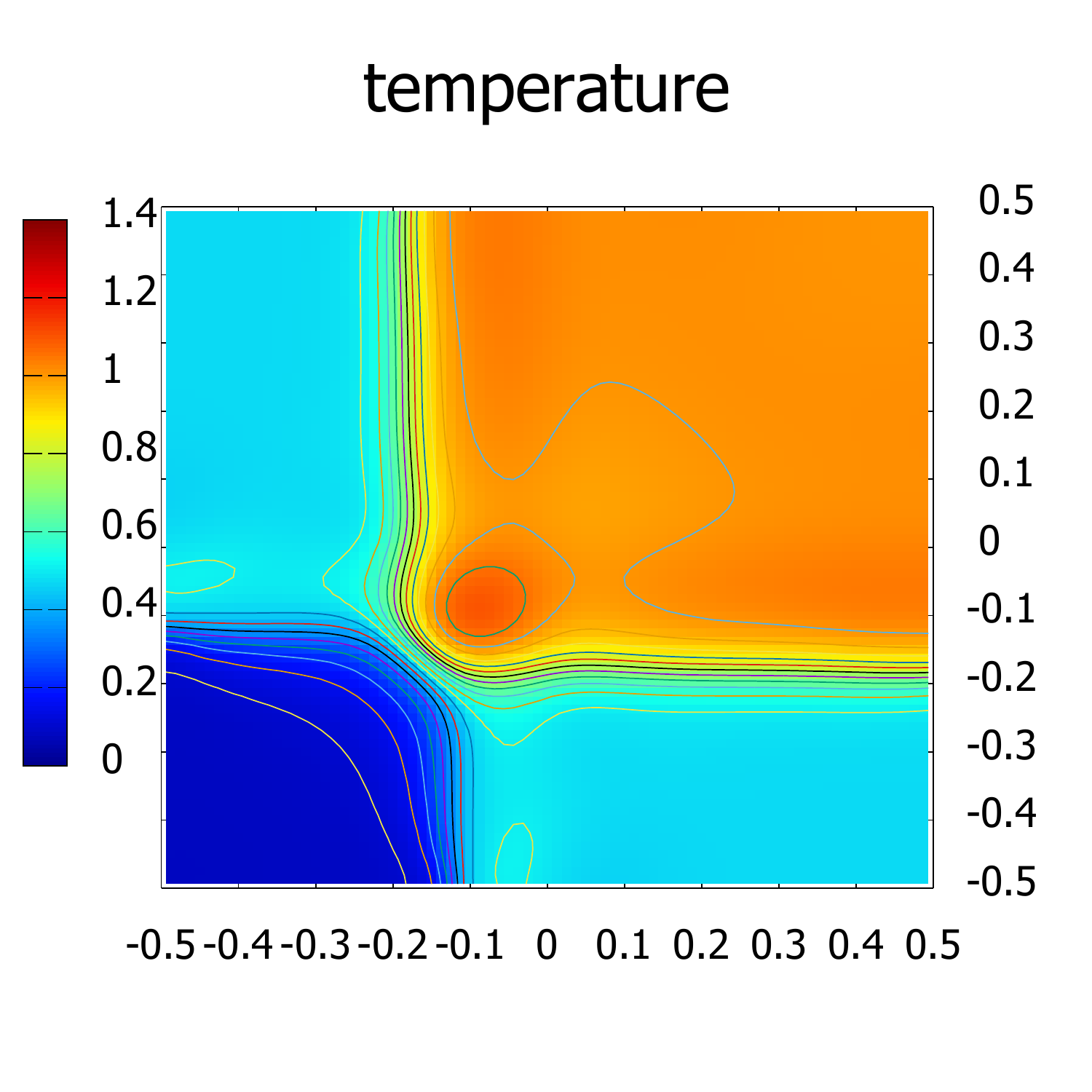}
	\includegraphics[width=2in]{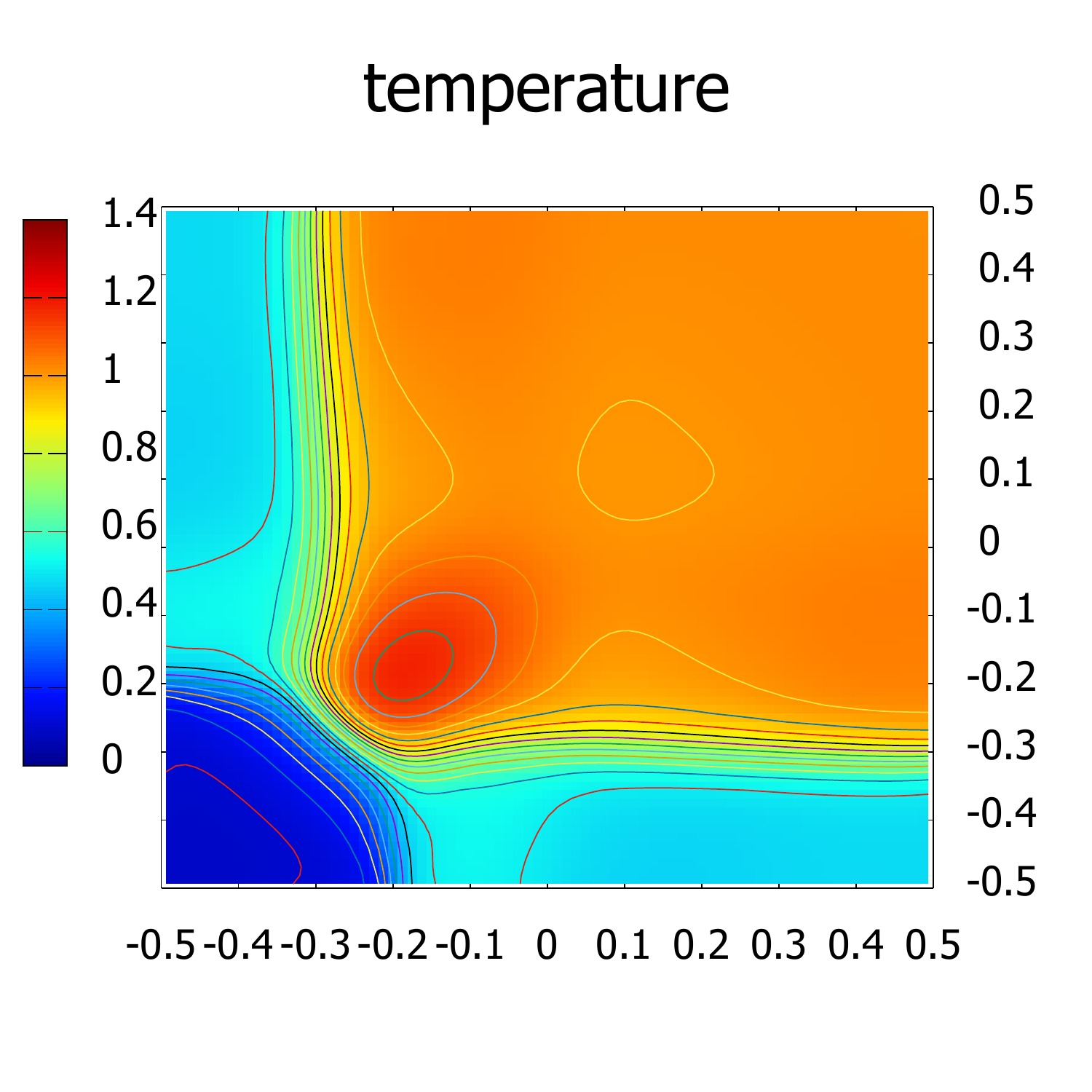}\\
	\includegraphics[width=2in]{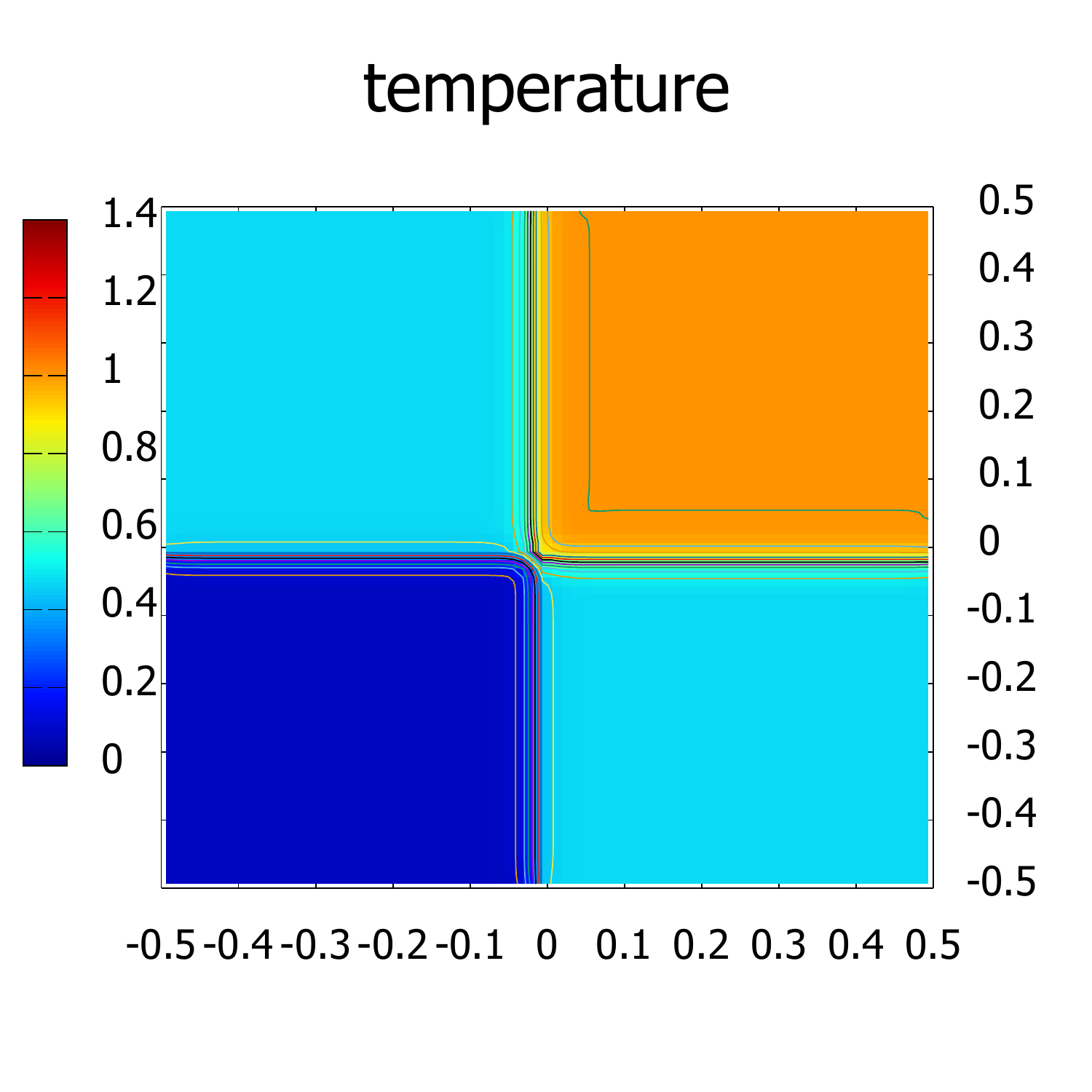}
	\includegraphics[width=2in]{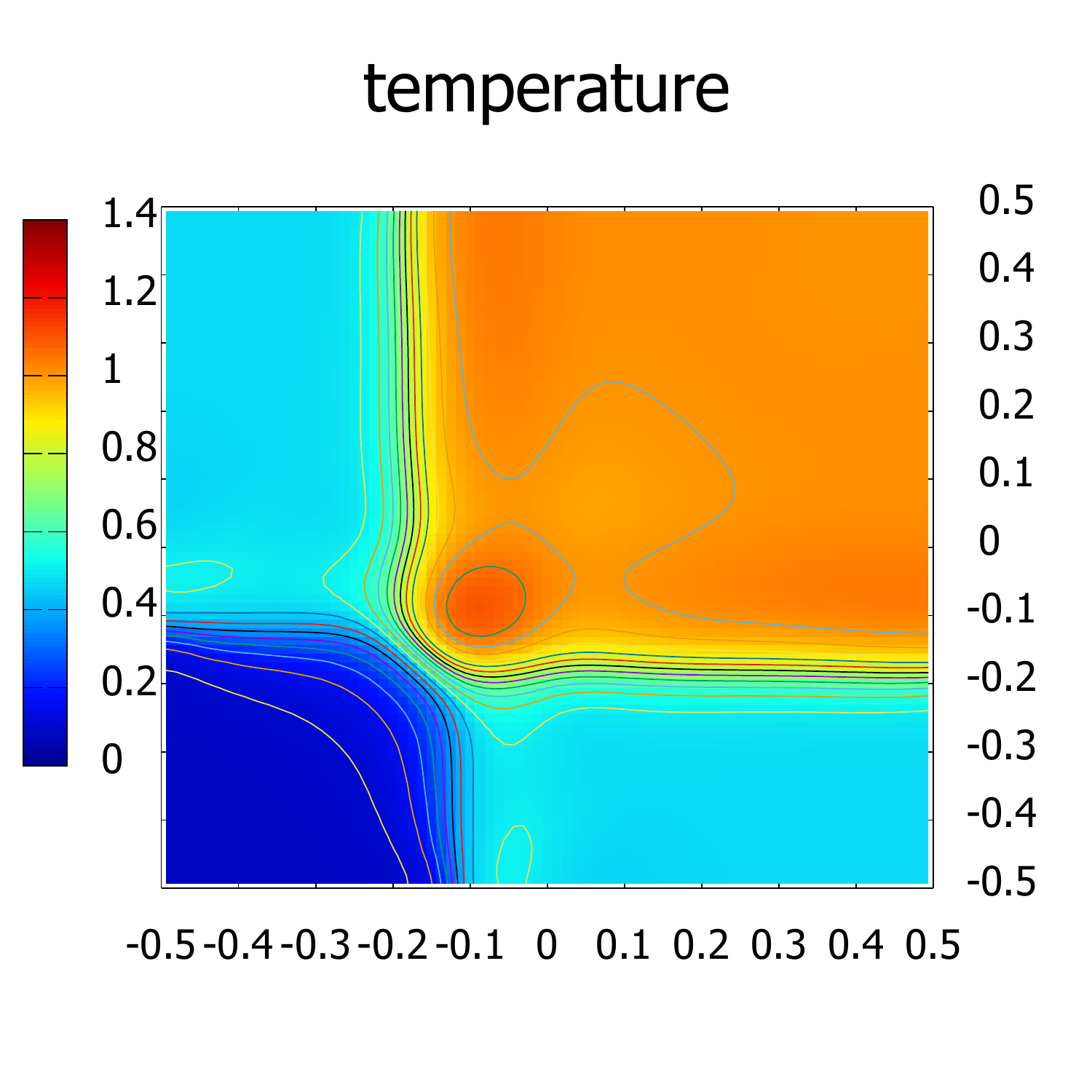}
	\includegraphics[width=2in]{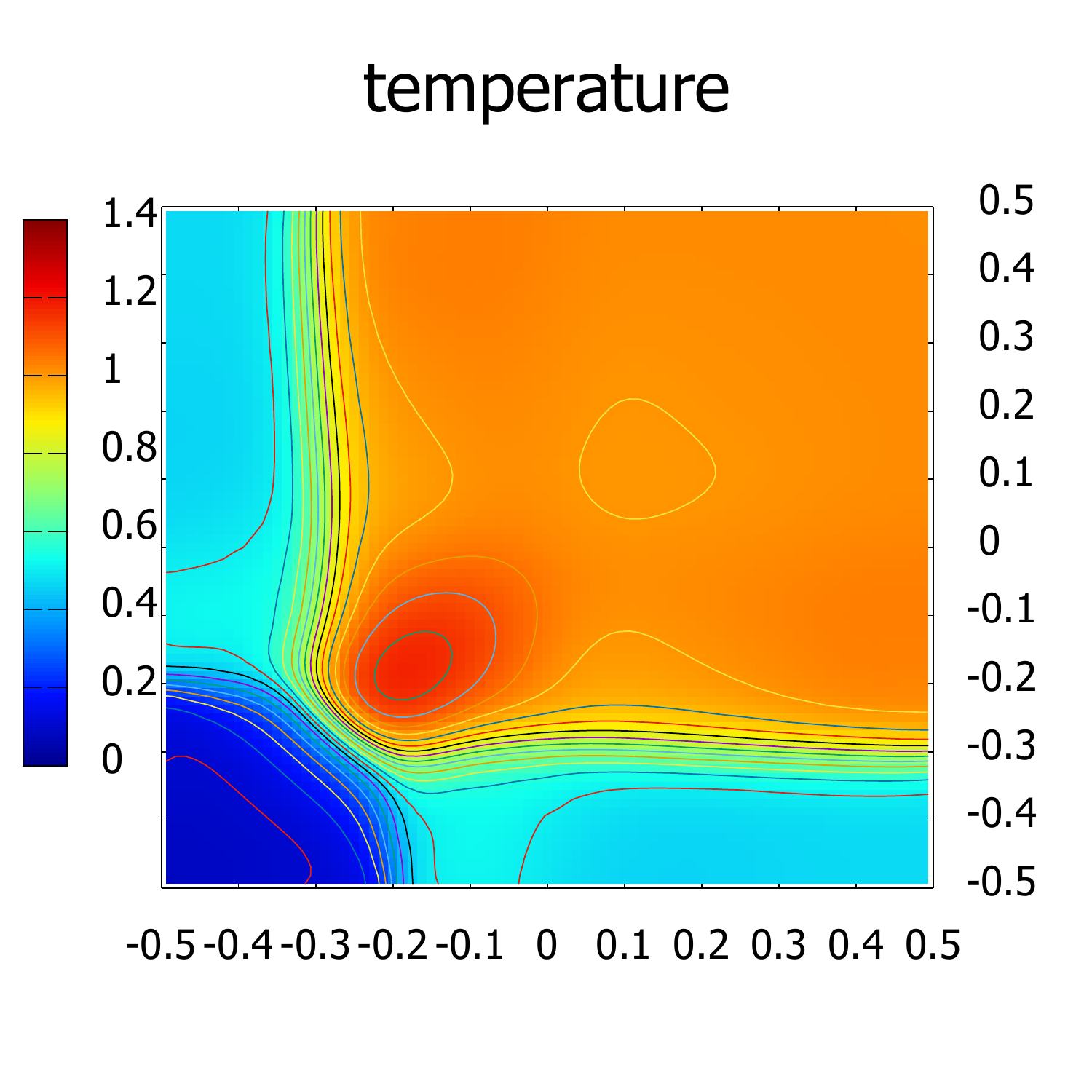}
	\caption{{\bf 2D Riemann problem.} Temperature profile
          obtained with a first order discontinuous Galerkin scheme
          with uniform grids with $N_x=N_y=80$ and $\eps=10^{-2}$. From left to right: $t=0.01, 0.2, 0.35$. Top: the full kinetic scheme; bottom: the hybrid scheme. 29 contour lines on the range $[0,1.4]$. }
	\label{fig34}
\end{figure}

In Figure \ref{fig33},  we first show the density profiles with $\eps=10^{-2}$, at time $t=0.01, 0.2, 0.35$ for the hybrid scheme and the full kinetic scheme, as well as the domain indicators for the hybrid scheme.
By comparing the contour lines, we can see that the hybrid scheme and the full kinetic scheme have almost the same results. However, due to a little large viscosity of $\eps=10^{-2}$, the four shock lines are smeared as time evolves. Our domain decomposition indicator can well capture the kinetic region
around the smeared shock lines. The corresponding temperature profiles are displayed in Figure \ref{fig34}, similarly almost the same results for the hybrid scheme and the full kinetic scheme can be observed.

\begin{figure}[ht]
\centering
\includegraphics[width=2in]{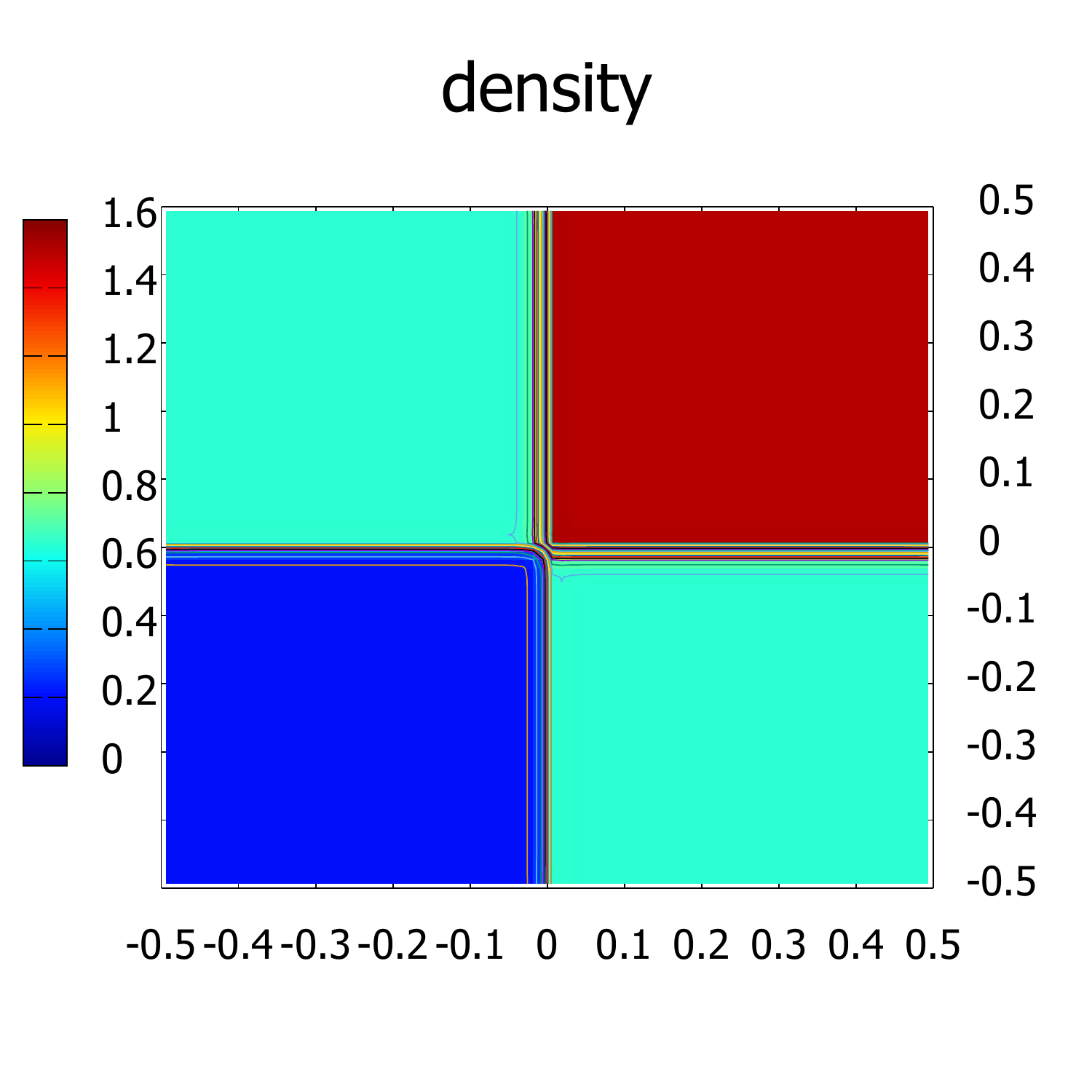}
\includegraphics[width=2in]{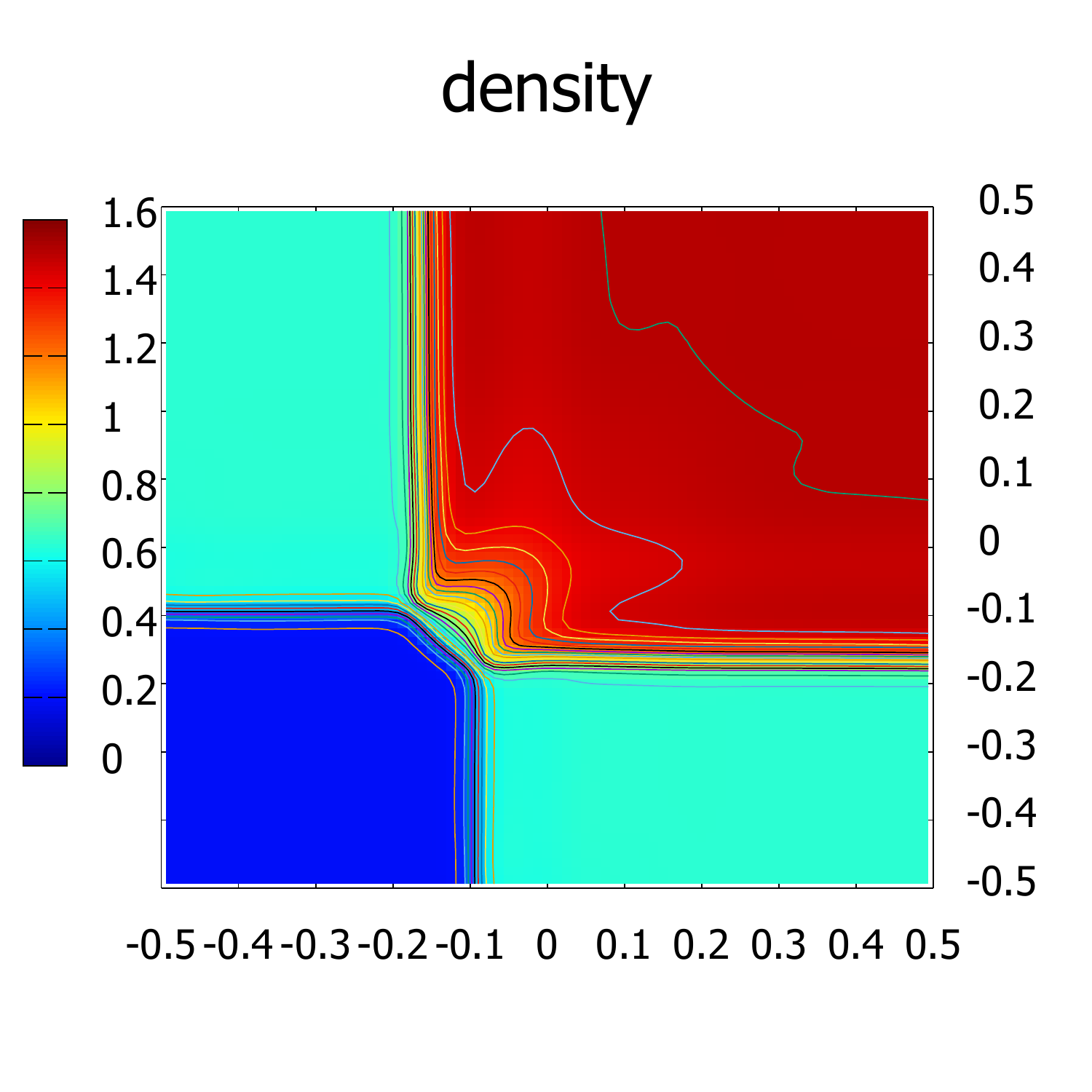}
\includegraphics[width=2in]{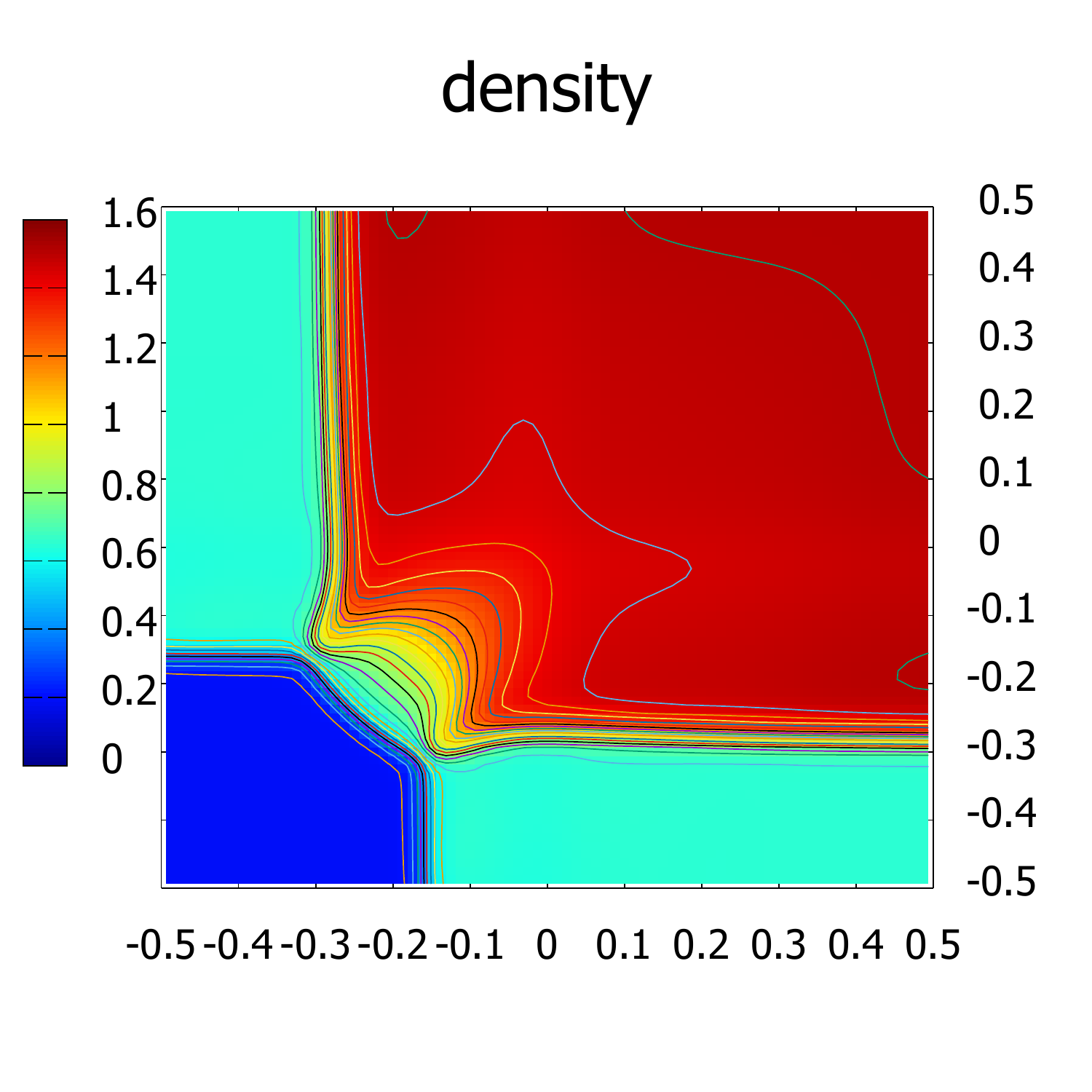}\\
\includegraphics[width=2in]{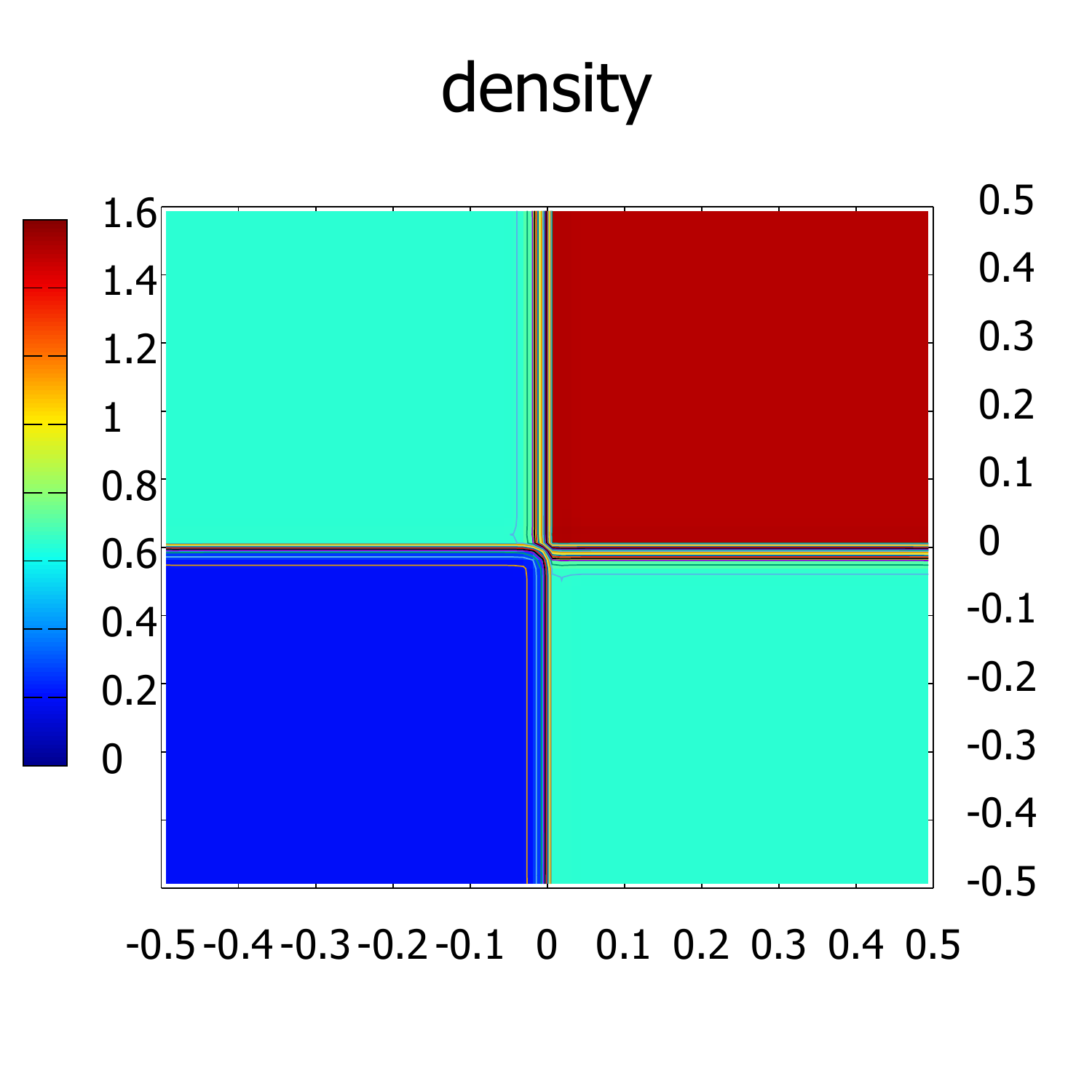}
\includegraphics[width=2in]{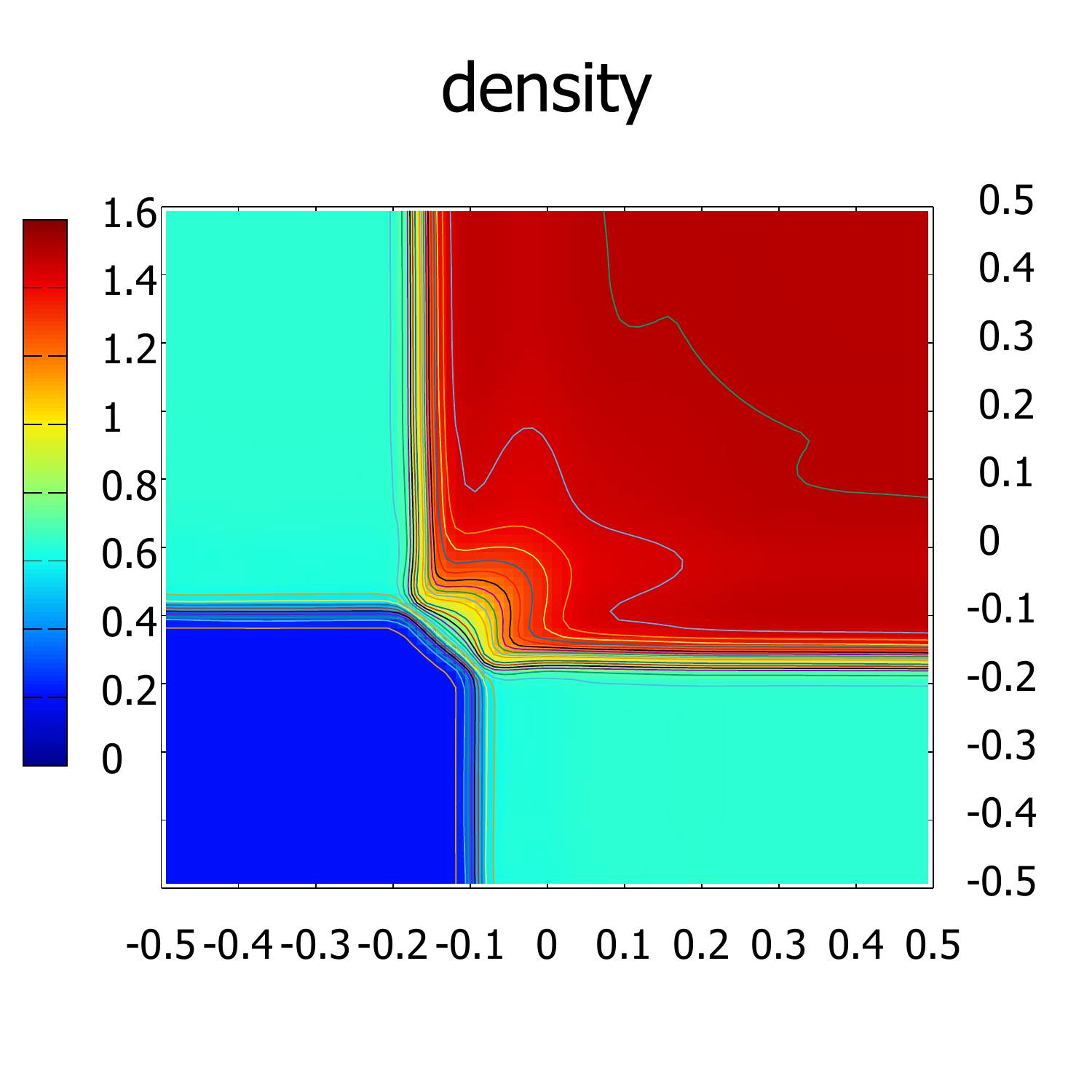}
\includegraphics[width=2in]{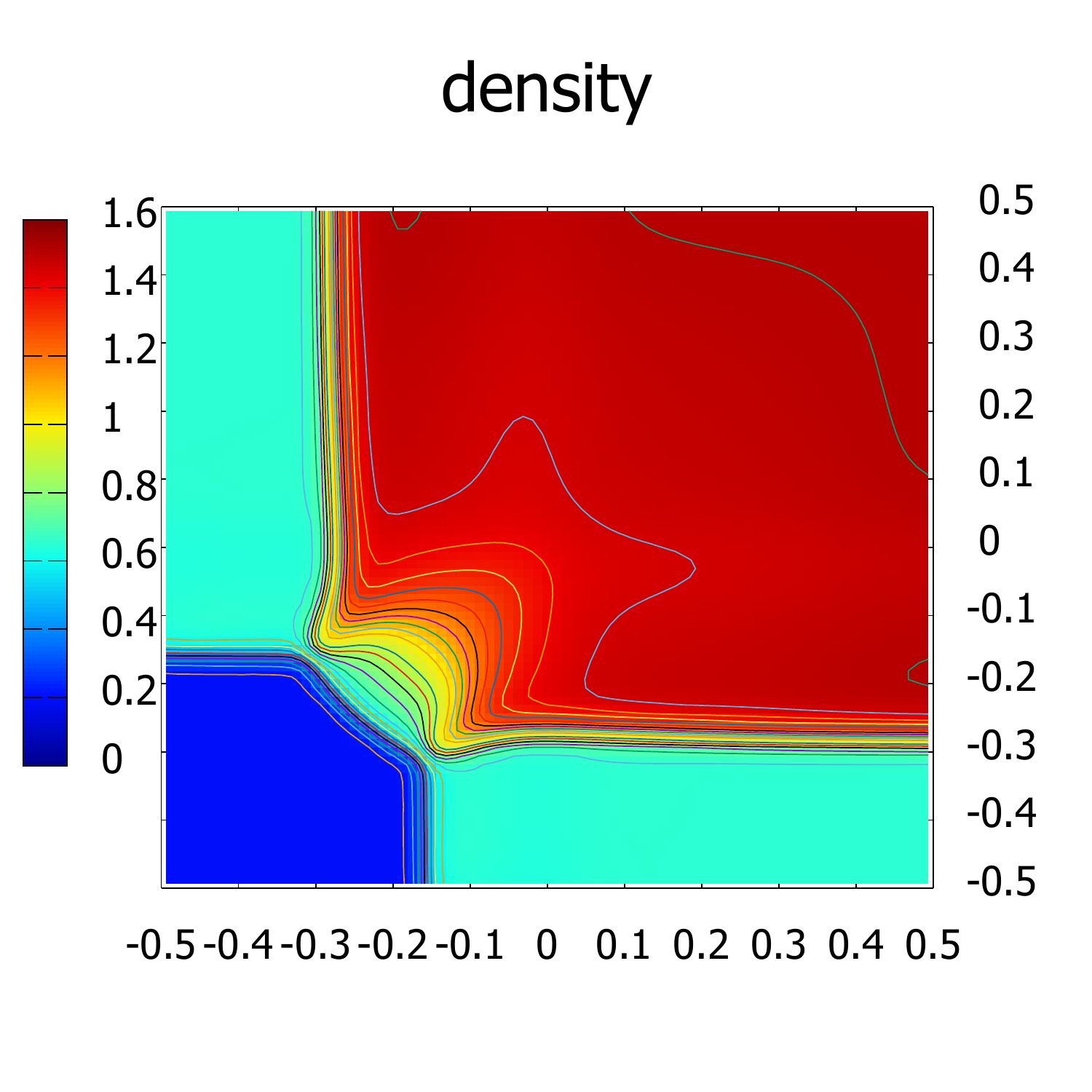}\\
\includegraphics[width=2in]{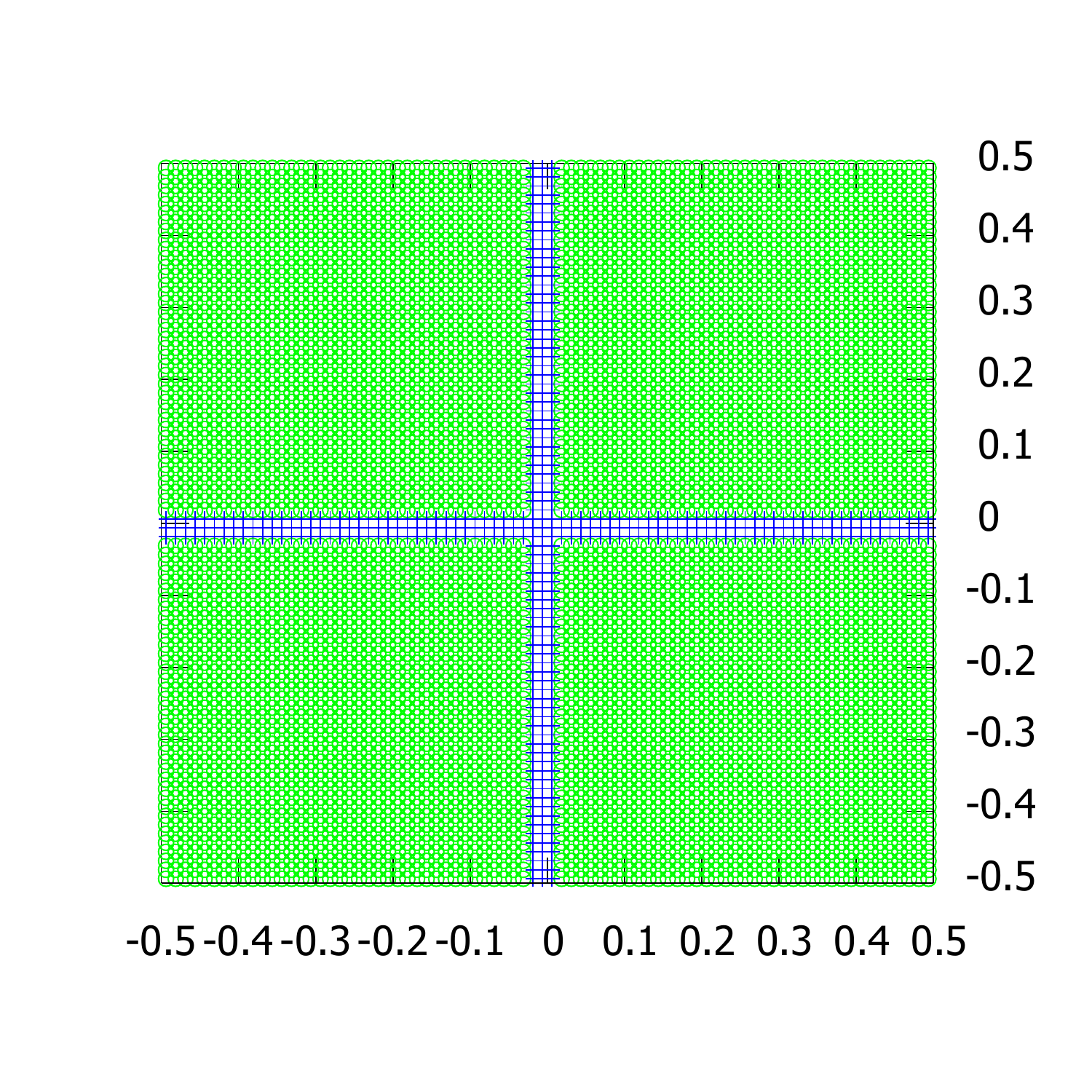}
\includegraphics[width=2in]{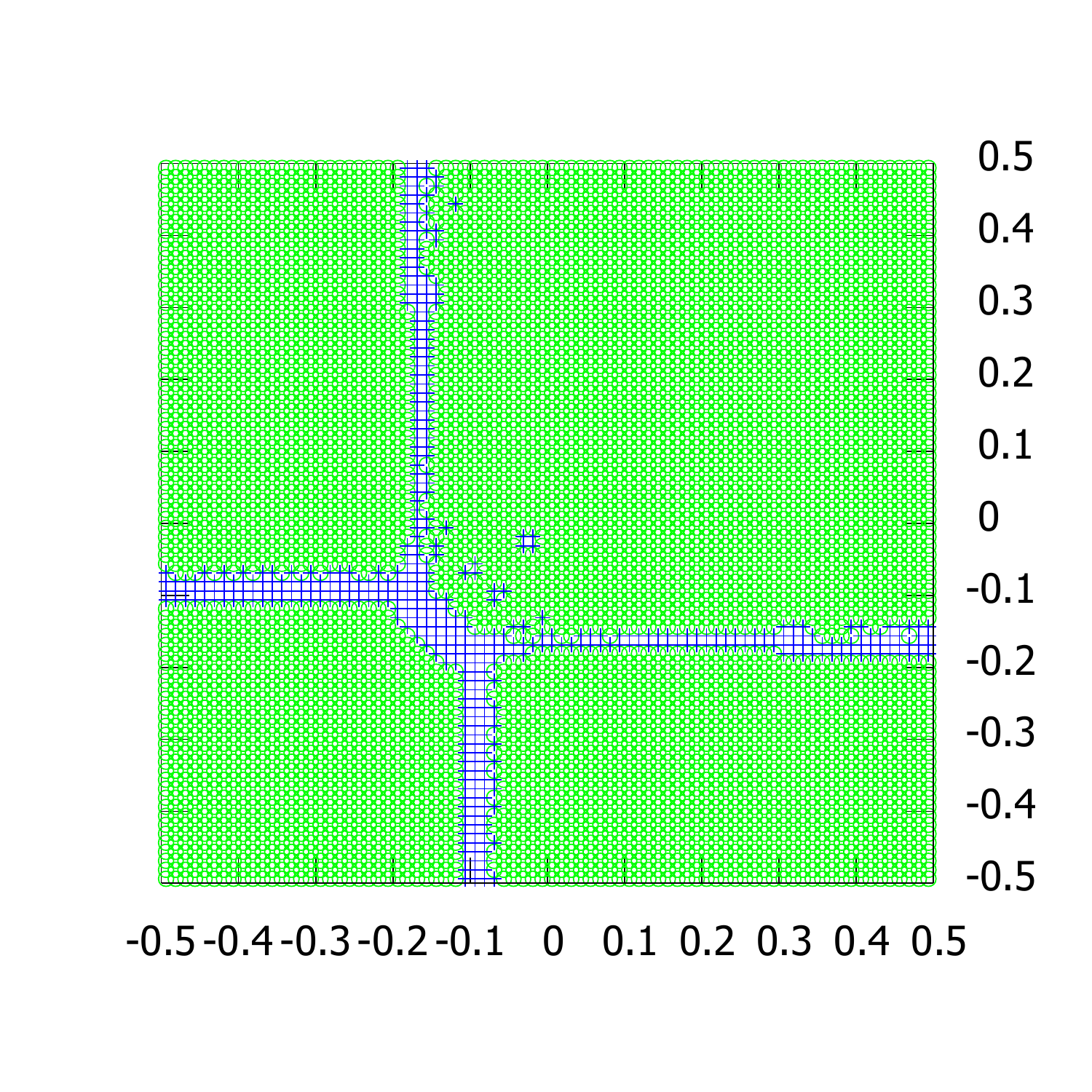}
\includegraphics[width=2in]{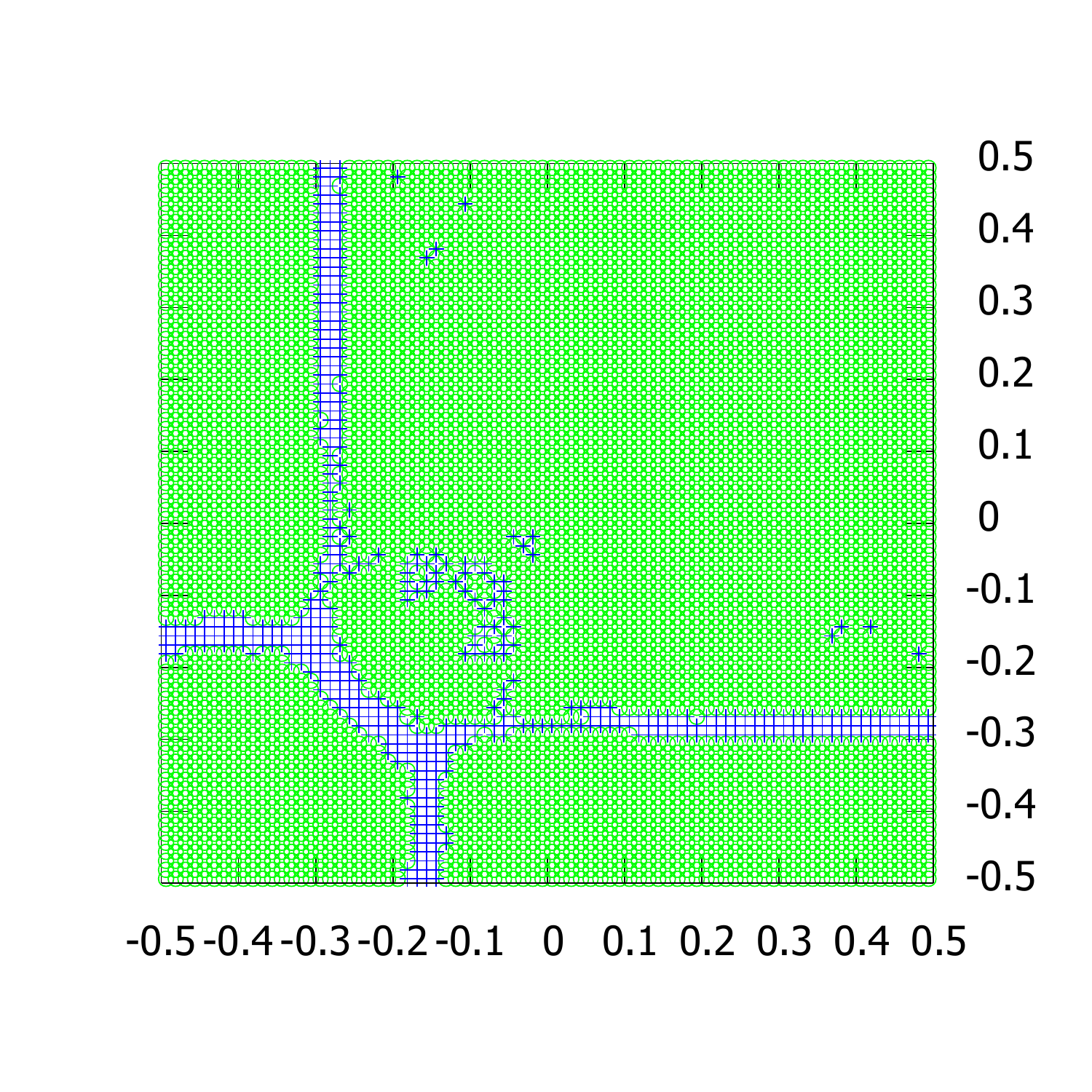}
\caption{{\bf 2D Riemann problem.} Density profile obtained with a
  first order discontinuous Galerkin scheme with uniform grids with
  $N_x=N_y=80$ and $\eps=10^{-3}$. From left to right: $t=0.01, 0.2, 0.35$. From top to bottom: the full kinetic scheme, the hybrid scheme, the domain indicator for the hybrid scheme. In the domain indicator, symbol ``+'' denotes kinetic cells, symbol ``o'' denotes hydrodynamic cells. 29 contour lines on the range $[0,1.6]$. }
\label{fig31}
\end{figure}

\begin{figure}[ht]
\centering
\includegraphics[width=2in]{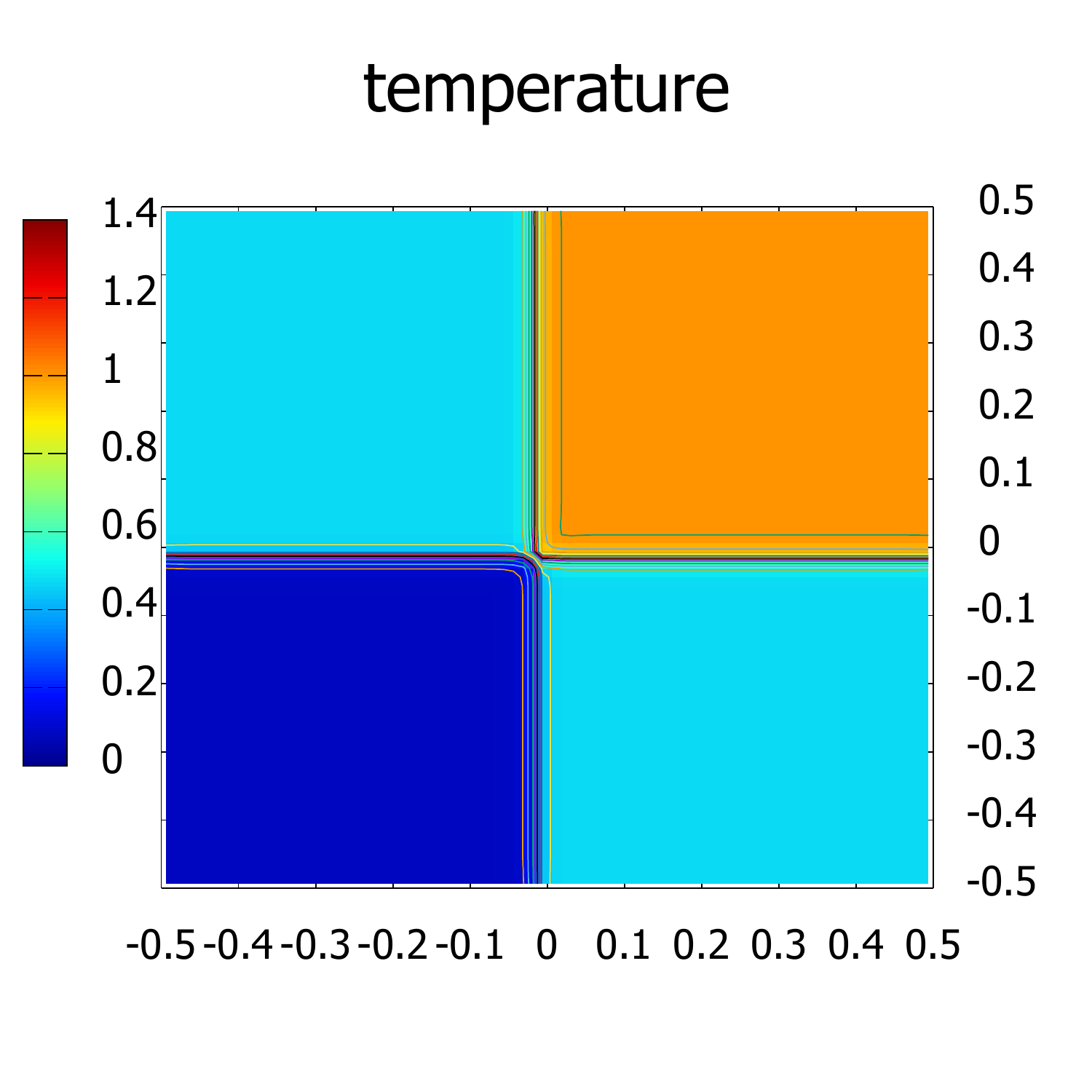}
\includegraphics[width=2in]{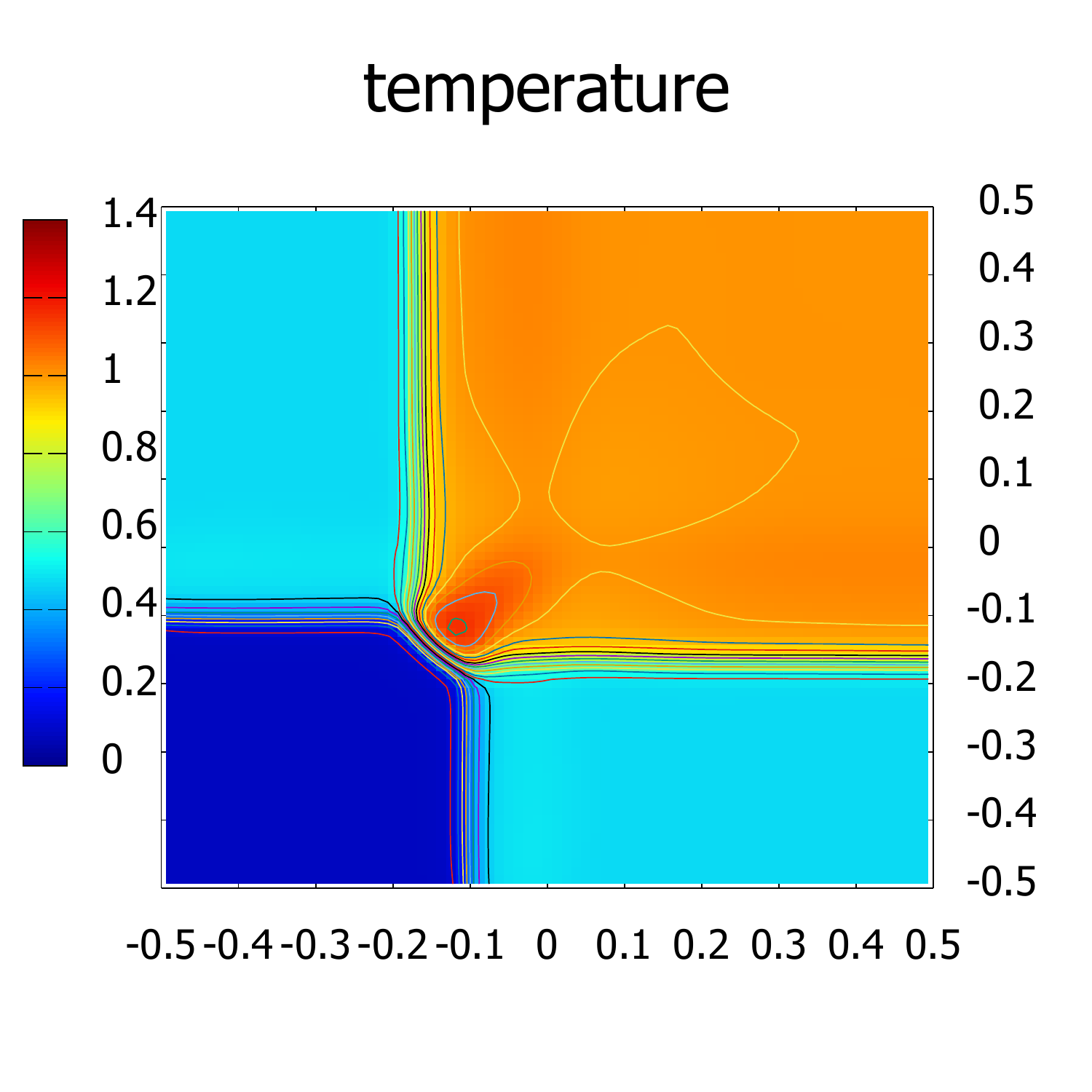}
\includegraphics[width=2in]{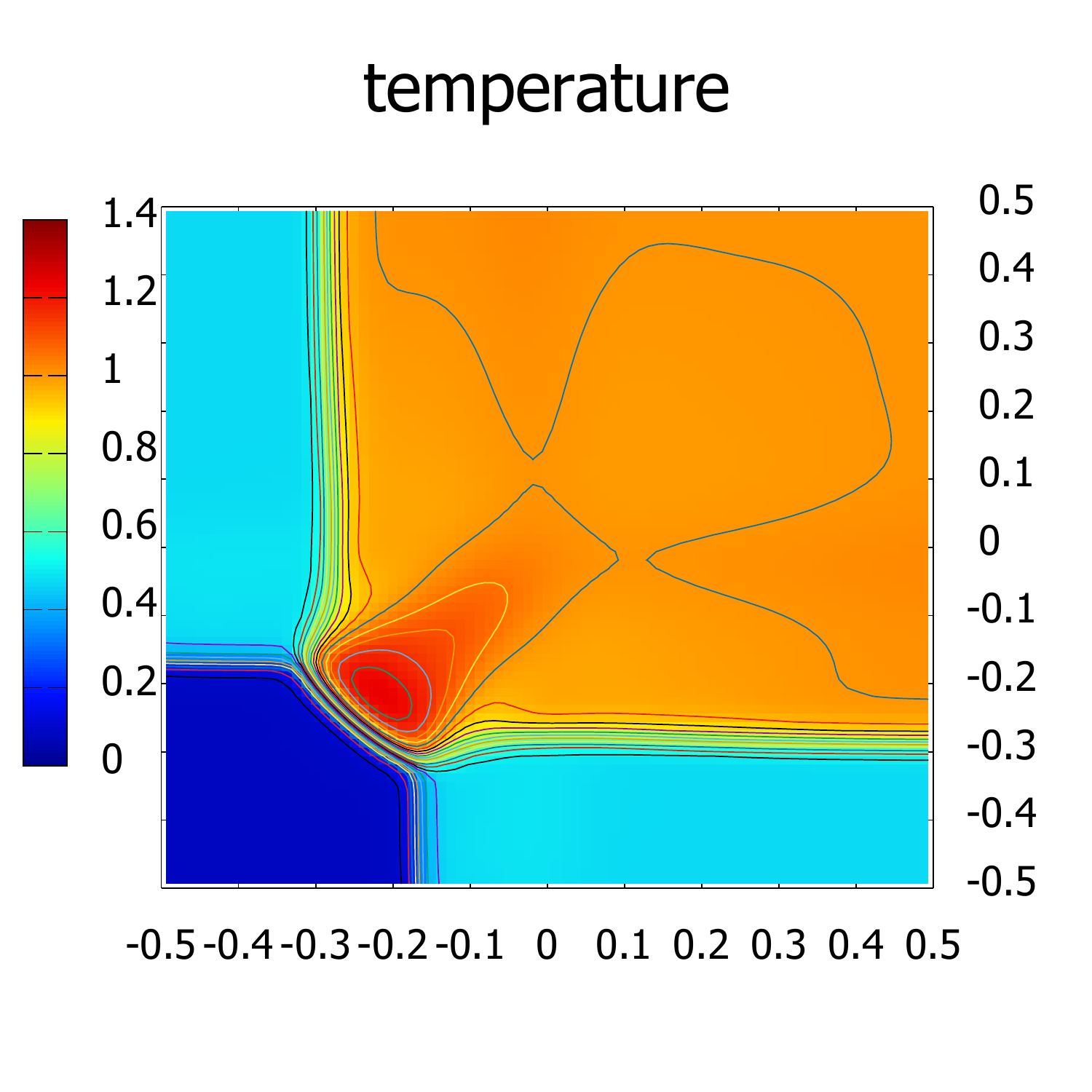}\\
\includegraphics[width=2in]{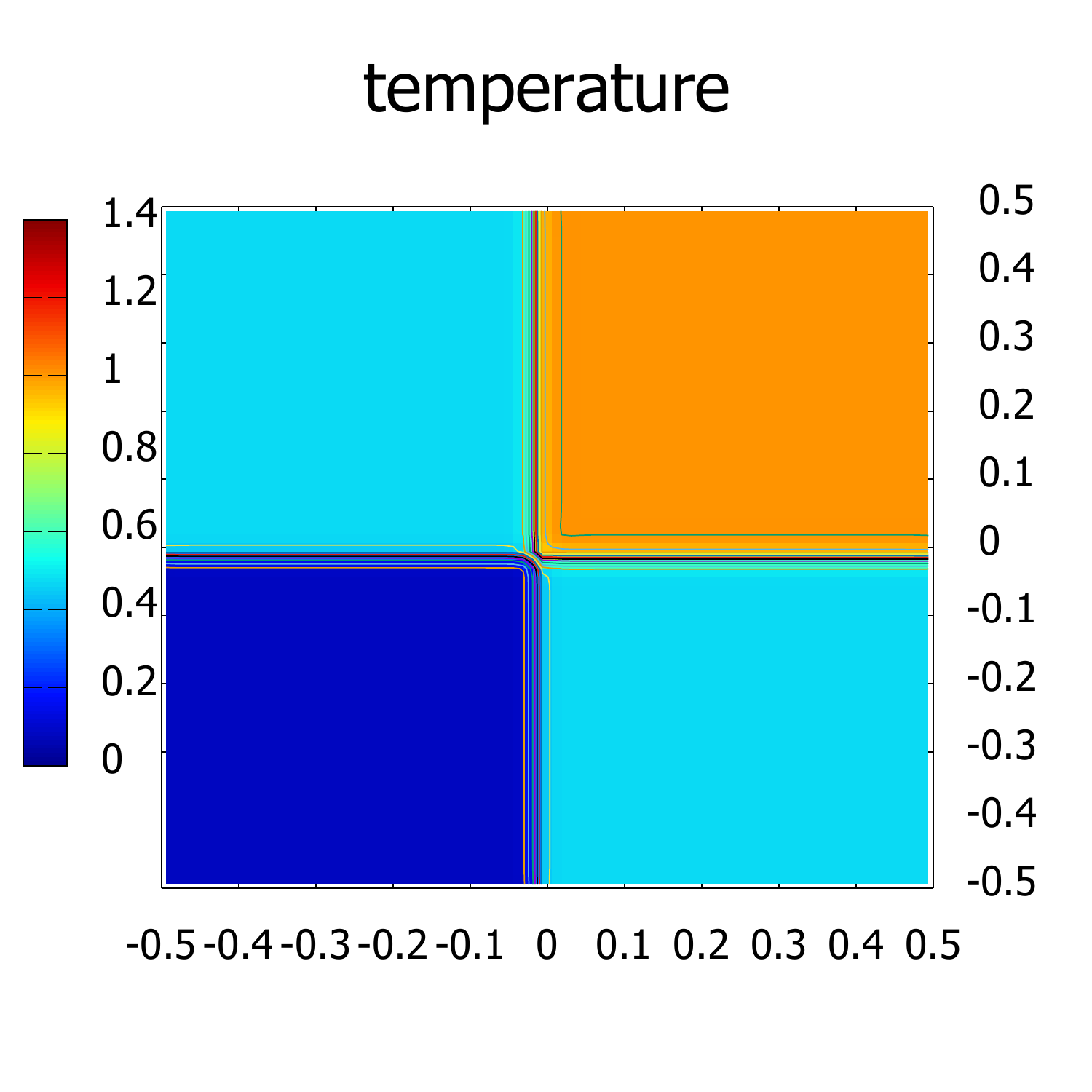}
\includegraphics[width=2in]{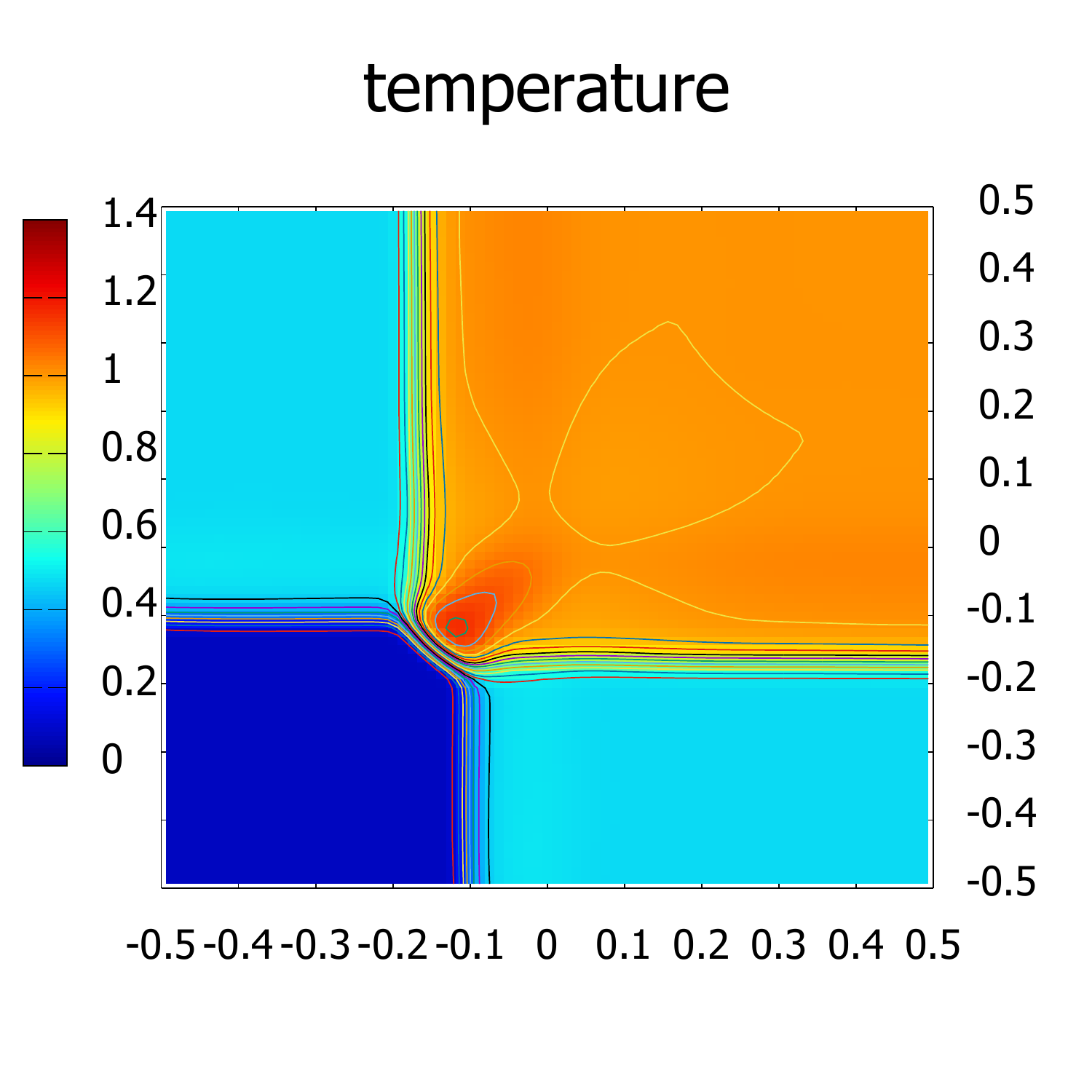}
\includegraphics[width=2in]{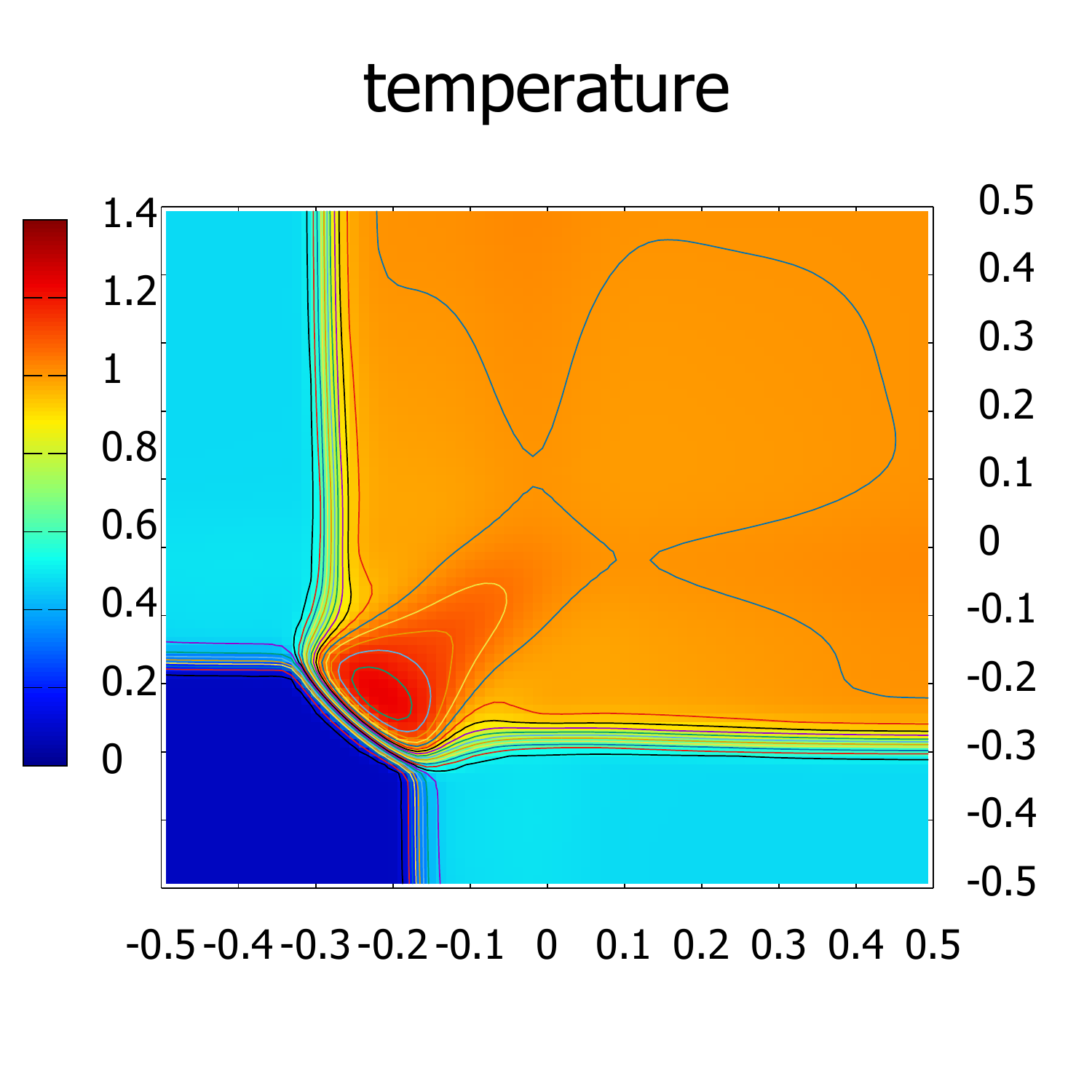}
\caption{{\bf 2D Riemann problem.} Temperature profile obtained with a
  first order discontinuous Galerkin scheme with uniform grids with
  $N_x=N_y=80$ and $\eps=10^{-3}$. From left to right: $t=0.01, 0.2, 0.35$. Top: the full kinetic scheme; bottom: the hybrid scheme. 29 contour lines on the range $[0,1.4]$. }
\label{fig32}
\end{figure}

In Figures \ref{fig31} and \ref{fig32}, we now show the density and temperature profiles with $\eps=10^{-3}$, at time $t=0.01, 0.2, 0.35$ for the hybrid scheme and the full kinetic scheme, as well as the domain indicators for the hybrid scheme. Similarly the hybrid scheme has almost the same results as the full kinetic one. However, due to small viscosity and sharp shock lines, our domain indicator tightly follows
the moving of the shock lines. This problem well demonstrates the good performance of our hybrid scheme and the domain indicator in the 2D case.

Concerning the CPU cost for the 2D problem, we run the code with MPI
parallelization on 4 processors up to $t=0.02$ with a time step
$\Delta t  \,=\, 2.5 \,\times\, 10^{-5}$ and $800$ time steps for $\eps=10^{-2}$, and
$t=0.2$ with a time step $\Delta t = 10^{-4}$ with $2000$ time steps for $\eps=10^{-3}$ respectively. We run the code three times, and the averaged CPU cost (real time) for the hybrid scheme is about $7$ minutes $33$ seconds for $\eps=10^{-2}$ and $15$ minutes and $4$ seconds for $\eps=10^{-3}$, while a full kinetic scheme costs about $14$ minutes and $54$ seconds for $\eps=10^{-2}$ and $37$ minutes for $\eps=10^{-3}$, so that the hybrid scheme saves about $1/2$ for $\eps=10^{-2}$ and $3/5$ for $\eps=10^{-3}$ of the CPU cost.
For the case of $\eps=10^{-2}$, due to large smearing of shock profile, more and more kinetic cells are identified as time goes on, see Figure \ref{fig33}. Our hybrid scheme although still is effective on capturing the kinetic cells, but is getting less efficient. However, we would also mention that there are some communication costs from MPI parallelization, generally the hybrid scheme may save even more as compared to the full kinetic scheme. 

\subsection{2D ghost effect.} 
\label{ex4}
On the basis of kinetic theory, the heat-conduction equation from the
stationary state solution of the classical Navier-Stokes equations is
not suitable for describing the temperature field of a gas in the
continuum limit in an infinite domain without flow at infinity, where the flow vanishes in this limit. By the asymptotic theory, as the Knudsen number of the system approaches zero, the temperature field should be obtained by the kinetic equation, this phenomenon is called the ghost effect \cite{filbet2012deterministic, sone1996inappropriateness, bobylev1995quasistationary}. In this example, we will numerically study this effect based on our hybrid scheme.
	
We consider a rarefied gas between two parallel plane walls at $x=0$ and $x=1$. The walls both have a common periodic temperature distribution $T_w$
\[
	T_w(y) \,=\, 1-0.5\cos(2\pi x),\,\, \forall y\in (0,1),
\]
and move in a common small mean velocity $u_w$ of order $\eps$ 
\[
        u_w(y) \,=\, (\eps,0).
\]
	
\begin{figure}[ht]
\centering
\includegraphics[width=2.8in]{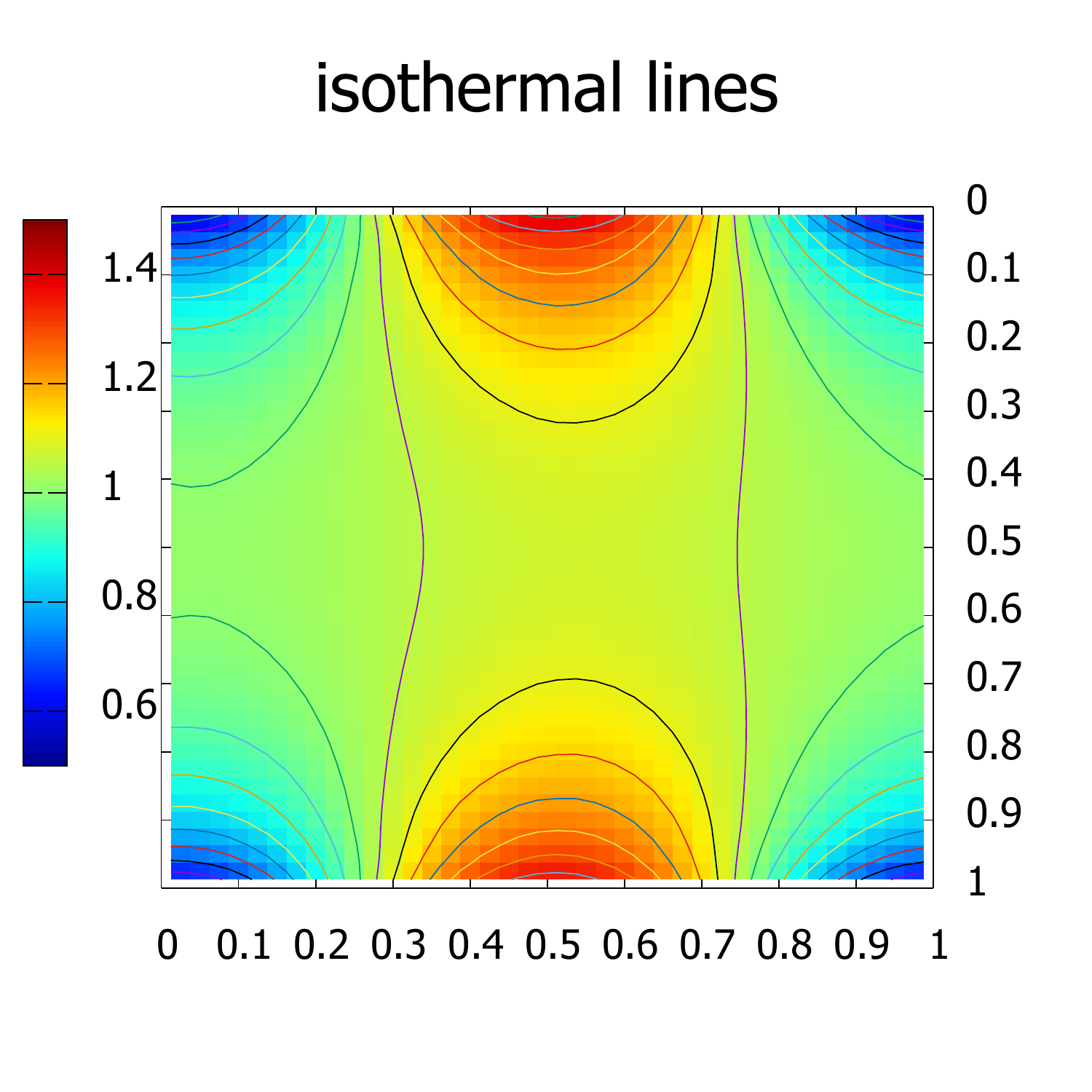}
\includegraphics[width=2.8in]{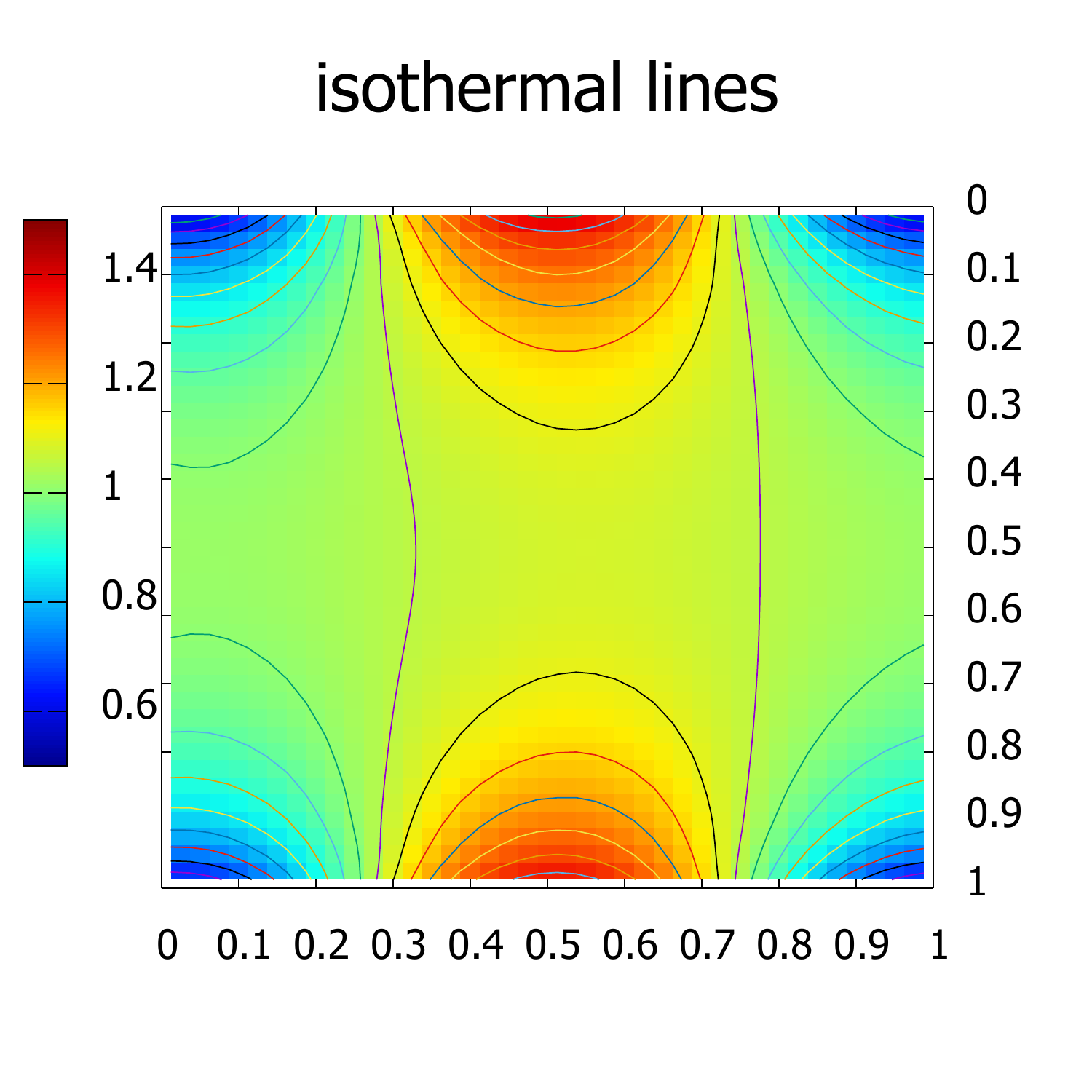}\\
\includegraphics[width=2.8in]{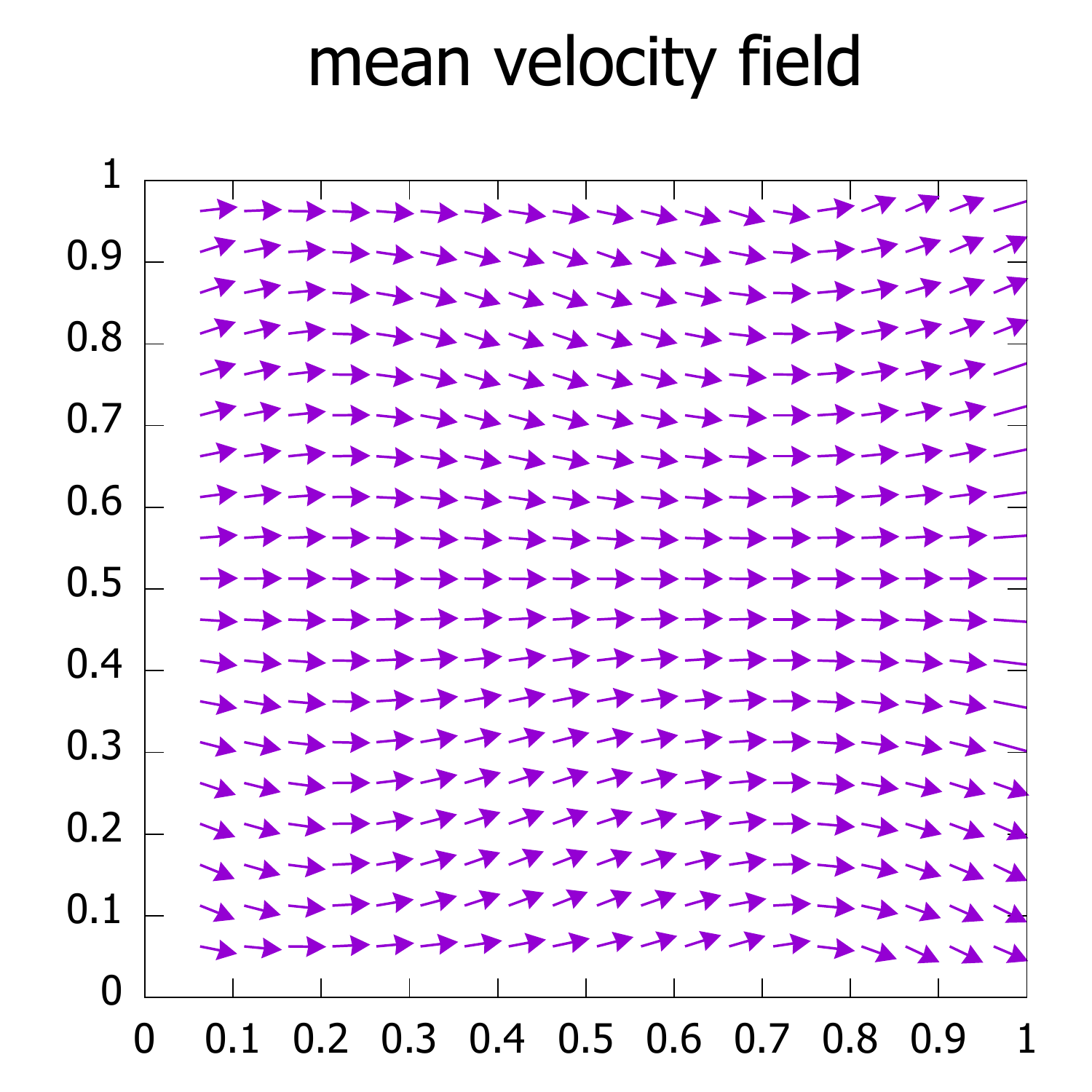}
\includegraphics[width=2.8in]{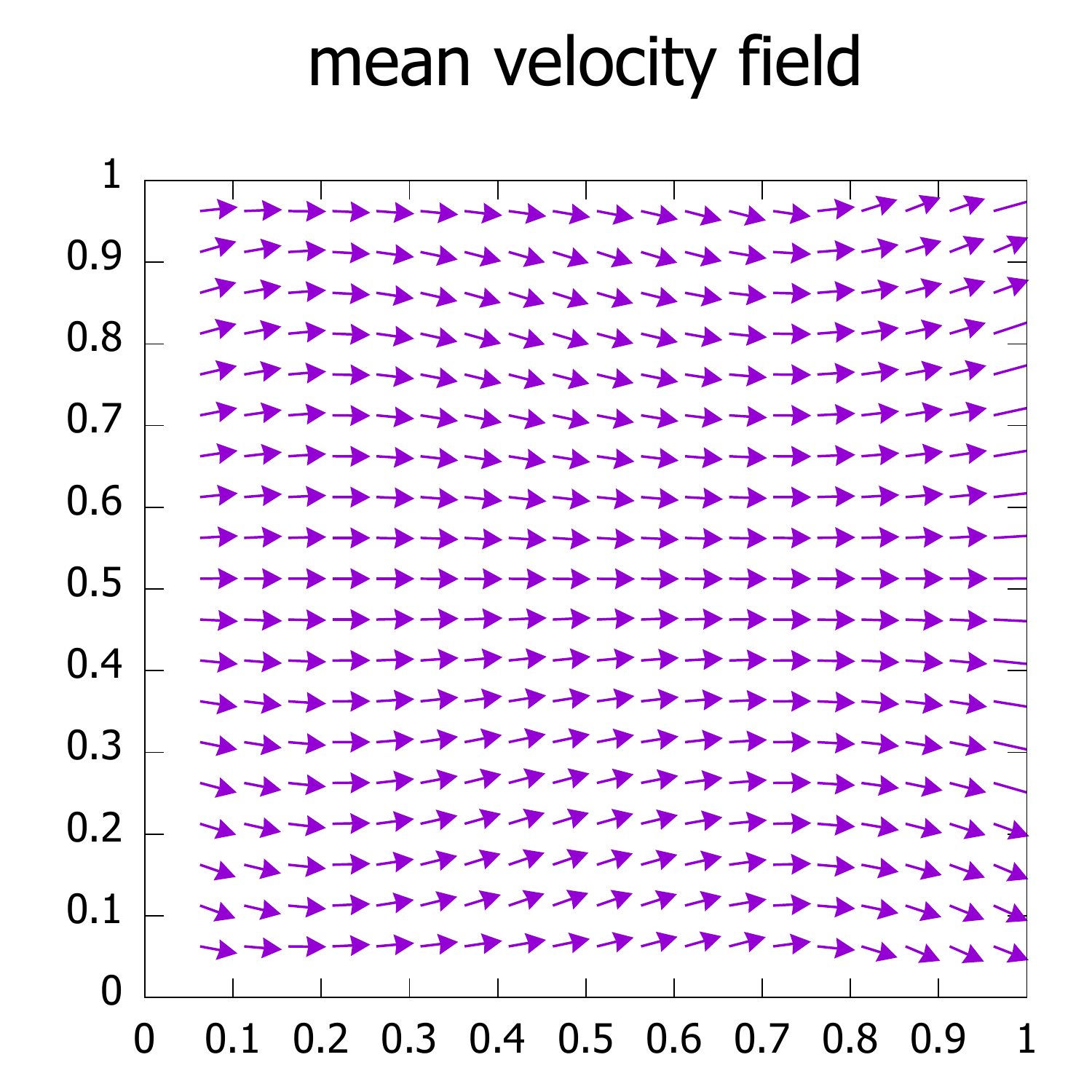}\\
\includegraphics[width=2.8in]{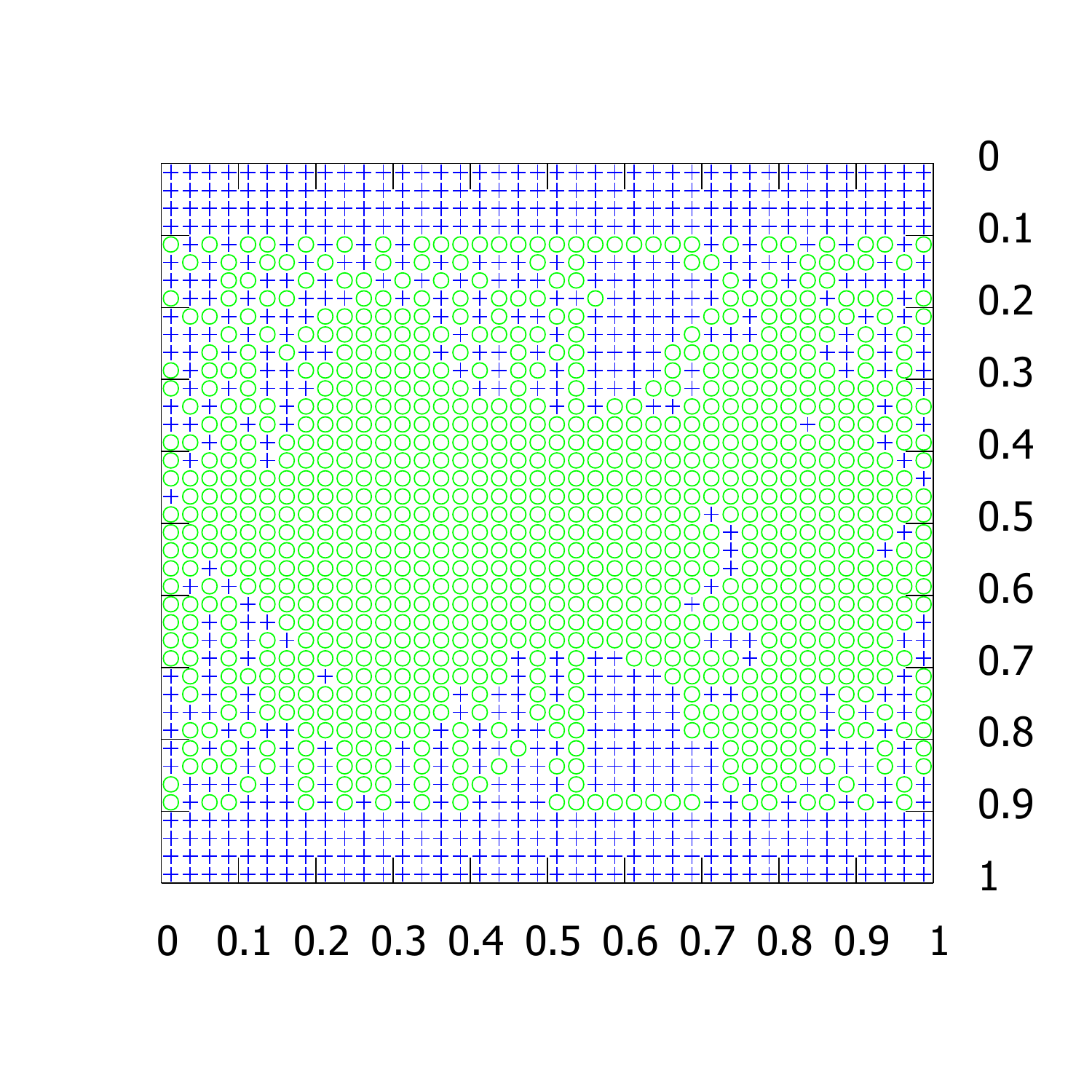}
\includegraphics[width=2.8in]{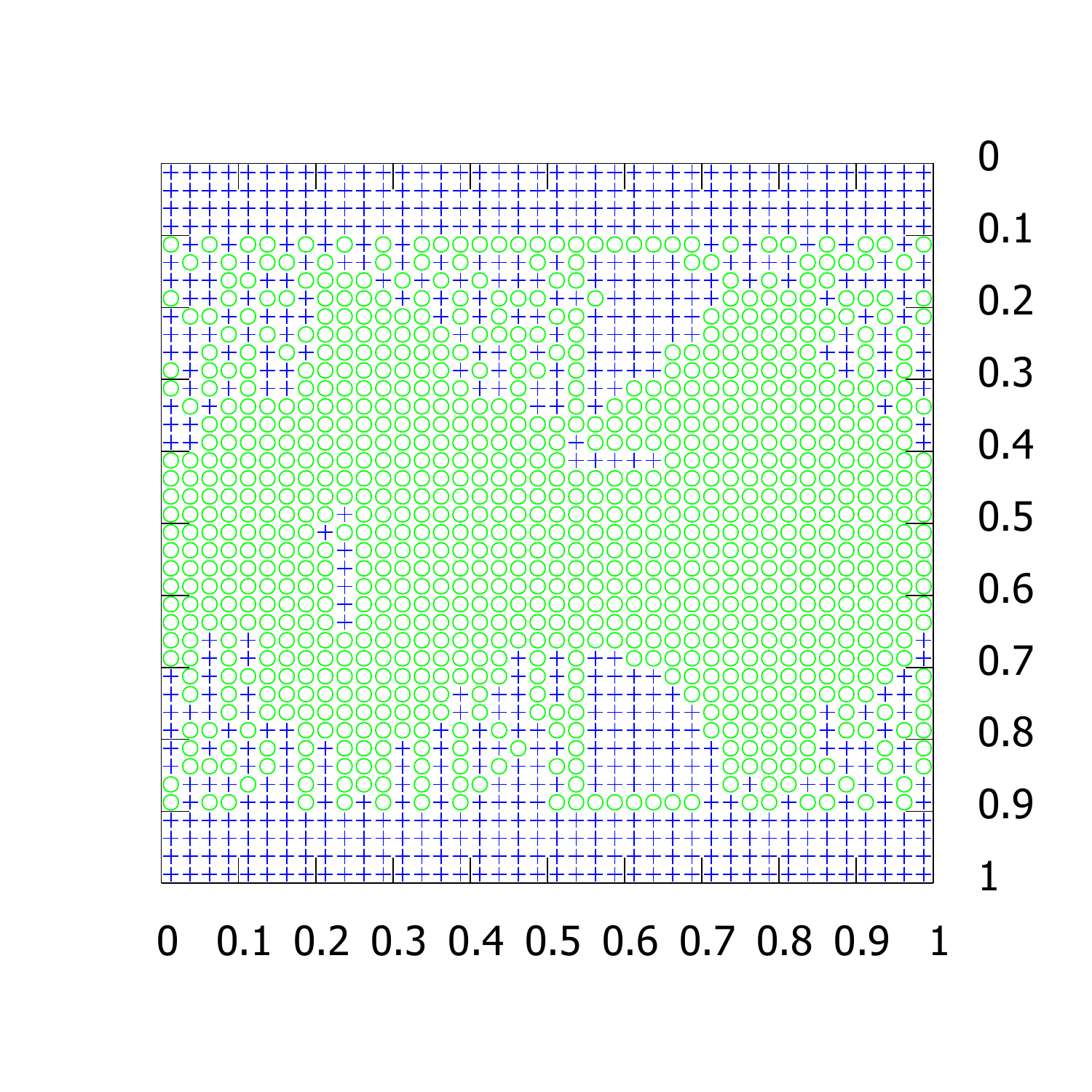}
\caption{{\bf The ghost effect.} Second order discontinuous Galerkin scheme. Uniform grids with $N_x=N_y=40$. $\eps=0.02$. From left to right: $t=10, 80$. From top to bottom: the isothermal lines, the mean velocity field and the domain indicator. Figures are rotated by $90^\circ$. In the domain indicator, symbol ``+'' denotes kinetic cells, symbol ``o'' denotes hydrodynamic cells. The isothermal lines increase by $0.05$ on the range $[0.45,1.55]$.}
\label{fig41}
\end{figure}

We take a uniform mesh with $N_x=N_y=40$ and time step $\Delta t=1/5000$. For velocity, since the solution is smooth and does not vary too much, we take a cut-off domain $\Omega_v=[-8,8]$ and discretize it with $N_v=16$ points along each direction. We apply a second order nodal discontinuous Galerkin scheme and force the grids inside the width of $0.1$ along $x$ direction around the walls always to be kinetic. We show the results (rotated by $90^{\circ}$) for $\eps=0.02$ in Figure \ref{fig41} for the isothermal lines, the mean velocity field as well as the domain indicators for the hybrid scheme, at time $t=10, 80$ respectively.  The solution at $t=80$ is assumed to be close to a steady state. We can observe that the solution at $t=80$ is similar to the results in \cite{sone1996inappropriateness} which is obtained from solving the BGK equation by a finite difference scheme, and also very similar to the results obtained from a $2d_x\times 2d_v$ Boltzmann system \cite{filbet2012deterministic}. 
	
From the domain indicators at different time in Figure \ref{fig41}, we can find that the indicator well defines the kinetic region where the solution has large variants (left and right sides) due to the effect from the moving walls, while other parts are almost in the hydrodynamic region (note that top and bottom are forced to be in the kinetic region). This example with long time simulation well demonstrates the good property of our hybrid scheme, in which we have almost only solved the expensive kinetic equation in the middle region where it is necessary and use the cheap hydrodynamic equations elsewhere. Due to the compactness of our discontinuous Galerkin scheme, we may observe some isolated hydrodynamic or kinetic cells, which show the robustness of the hybrid discontinuous Galerkin method in the multi-dimensional case, while this cannot be easily achieved by a finite volume scheme due to the requirement of wide stencils.

%%%%%%%%%%%%%%%%%%%%%%%%%%%%%%%%%%
%
%%%%%%%%%%%%%%%%%%%%%%%%%%%%%%%%%%
\section{Conclusion}
\label{sec6}
\setcounter{equation}{0}
\setcounter{figure}{0}
\setcounter{table}{0}

In this paper, we developed a hierarchical of hybrid discontinuous Galerkin scheme for some physically relevant problems in both 1D and 2D space dimensions. The compact discontinuous Galerkin scheme has shown its robustness on the domain decomposition and $h$-$p$ adaptivity on capturing the boundary layers. Although only a first order time discretization is used in the paper for some long time simulations, it can be easily generalized to high order in time based on some high order 
implicit-explicit time discretizations. Extensions to unstructured meshes on more complicated geometries and to ES-BGK or Boltzmann collision operators will be investigated in our future work.

%%%%%%%%%%%%%%%%%%%%%%%%%%%%%%%%%%%
%
%%%%%%%%%%%%%%%%%%%%%%%%%%%%%%%%%%%

\section*{Acknowledgement}
%
%Francis Filbet and Eric Sonnendr\"ucker  were supported by the EUROfusion Consortium and has received funding
%from the Euratom research and training programme 2014-2018 under grant
%agreement No 633053. The views and opinions expressed herein do not
%necessarily reflect those of the European Commission.

T. Xiong acknowledges support from the Marie Sk{\l}odowska-Curie Individual Fellowships H2020-MSCA-IF-2014 of the European Commission, under the project HNSKMAP 654175. This work has also been supported by NSFC grant 11601455, NSAF grant U1630247, NSF grant of Fujian Province 2016J05022 and the Fundamental Research Funds for the Central Universities No. 20720160009.

%%%%%%%%%%%%%%%%%%%%%%%%%%%%%%%%%%%%%%%%%%
%
%%%%%%%%%%%%%%%%%%%%%%%%%%%%%%%%%%%%%%%%%%
\appendix
\section{BGK equation and Chu reduction}
\label{sec3}
\renewcommand{\theequation}{A.\arabic{equation}}
\setcounter{equation}{0}
\setcounter{figure}{0}
\setcounter{table}{0}

We consider three dimensional problems in velocity, but to save the
computational cost, we take the BGK equation
and use the technique of Chu reduction \cite{chu1965kinetic}. It has also
recently been used by M. Groppi, G. Russo and G. Stracquadanio in \cite{groppi2015high}. The main idea is that under suitable geometry symmetry assumptions, we can integrate the distribution function $f$ in the velocity phase space and introduce several auxiliary variables, to transform the BGK equation \eqref{eq:bgk} in 3D in velocity into a system of equations in 2D or 1D in velocity, such that the system has the same dimension in space and in velocity. We refer to \cite{chu1965kinetic} for more details.

\subsection*{Chu reduction}
%\label{sec3.1}
We start with the BGK equation in 3D in both space and velocity, which reads
\beq
\pt f \, + \, \bfv \cdot \nablax f \, = \, \frac{\nu}{\eps} \, \left(\mathcal{M}[f] - f\right).
\label{bgk3d}
\eeq
$\bfx = (x,y,z)$, $\bfv = (v_1, v_2, v_3)$ and $\mathcal{M}[f]$ is the Maxwellian defined as 
\beq
\mathcal{M}[f] \, :=\, \mathcal{M}(t,\bfv,\bfU) \, = \, \frac{\rho}{(2\pi T)^{3/2}}\exp\left(-\frac{|\bfv-\bfu|^2}{2T}\right),
\label{maxw3d}
\eeq
where $\bfu = (u_1, u_2, u_3)$  is the macroscopic mean velocity.  As explained in \cite{chu1965kinetic, groppi2015high}, we consider the BGK equation to problems with axial symmetry, where all transverse spatial gradients vanish and the gas is drifting only in the axial directions. In such cases, $f(t,\bfx, \bfv)$ depends on the full velocity, that is, molecular trajectories are 3D, however all transverse components of the macroscopic velocity $\bfu$ vanish. 

We first consider problems are axial symmetry with respect to two axises $(x_1,x_2)$, so that the problem does not depend on $x_3$ and the macroscopic velocity $\bfu$ vanishes at $u_3=0$.  For simplicity, let us now denote $\tilde\bfx = (x_1,x_2)$, $\tilde\bfv = (v_1, v_2)$, $\tilde\bfu = (u_1, u_2)$ and remove the dependence of $f$ on $x_3$. We introduce two new unknowns
\[
g_1(t,\tilde\bfx,\tilde\bfv) \, := \, \int_{\mathbb{R}} f(t,\tilde\bfx, \tilde\bfv, v_3)\, dv_3,  \quad g_2(t,\tilde\bfx,\tilde\bfv) \, := \, \int_{\mathbb{R}} \frac{v_3^2}{2} f(t,\tilde\bfx, \tilde\bfv, v_3)\, dv_3,
\]
which are integrations on $v_3$ and from assumption we notice that 
\[
\int_{\mathbb{R}} v_3 \, f(t, \tilde\bfx, \tilde\bfv, v_3) \,dv_3 \, = \, \rho\,u_3 \,=\,0.
\]
Multiplying \eqref{bgk3d} by $1$ and $v_3^2/2$, integrating on $v_3$, it yields 
a system for the new unknown vector $\bfg=(g_1, g_2)$, coupled with suitable initial conditions
\beq
\label{bgk2d}
\left\{
\begin{array}{rl}
	\pt g_i \, + \, \tilde\bfv \cdot \nabla_{\tilde\bfx} \,g_i &  = \ds\,\,  \frac{\nu}{\eps} \, \left(\mathcal{M}[g]_i - g_i \right), 
	\\
	\,
	\\
	g_i(t=0) &  = \,\,  g_{i,0}, \quad  i \,=\, 1, \, 2.
\end{array}
\right.
\eeq
The new BGK system \eqref{bgk2d} describes a relaxation process towards the vector function $(\mathcal{M}[g]_1, \mathcal{M}[g]_2)$, which is a Chu reduction of the Maxwellian \eqref{maxw3d}
from 3D in velocity to 2D in velocity and has the form
\[
(\mathcal{M}[g]_1, \mathcal{M}[g]_2) \, = \, (\mathcal{M}[g]_1, \frac{T}{2}\mathcal{M}[g]_1),
\]
where
$$
%\beq
\mathcal{M}[g]_1(\tilde\bfv) \, := \, \int_{\mathbb{R}} \mathcal{M}[f](\tilde\bfv,v_3) \,dv_3 = \frac{\rho}{2\pi T}\exp\left(-\frac{|\tilde\bfv-\tilde\bfu|^2}{2T}\right).
%\label{maxw2d}
%\eeq
$$
We only list the argument $\tilde\bfv$ in $\mathcal{M}[g]_1(\tilde\bfv)$ to emphasize it is $\tilde\bfv$-dependent, while $(t, \tilde\bfx)$ are omitted.
The macroscopic moments $\rho$, $\tilde\bfu$ and $T$ of $f$ can be given in terms of $\bfg$ as
$$
%\beq
%\label{U2d}
\rho = \int_{\mathbb{R}^2}g_1(\tilde\bfv) \, d \tilde\bfv, \quad \tilde\bfu = \frac{1}{\rho} \int_{\mathbb{R}^2}\tilde\bfv \, g_1(\tilde\bfv) \, d \tilde\bfv, \quad 3T  \,=\, \frac{1}{\rho}\left[\int_{\mathbb{R}^2}\frac{|\tilde\bfv-\tilde\bfu|^2}{2} \, g_1(\tilde\bfv) \, d \tilde\bfv + \int_{\mathbb{R}^2} g_2(\tilde\bfv) \, d \tilde\bfv\right].
%\eeq
$$
If we further assume the model is axial symmetry only with respect to axis $x_1$, and $f$ only depends on $x_1$, we can transform \eqref{bgk3d} into a system of 1D problem in both $x$ and $v$. The procedure is similar, that is, we introduce
\[
h_1(t,x,v) \, := \, \int_{\mathbb{R}^2} f(t,x, v, v_2, v_3) \,dv_2\, dv_3,  \quad h_2(t,x,v) \, := \, \int_{\mathbb{R}^2} \frac{v_2^2+v_3^2}{2} f(t, x, v, v_2, v_3)\, dv_2 \,dv_3,
\]
multiplying \eqref{bgk3d} by $1$ and $(v_2^2+v_3^2)/2$, integrating with respect to $v_2$ and $v_3$, it yields a system for the unknown vector $\bfh = (h_1, h_2)$, coupled with suitable initial conditions
\beq
\label{bgk1d}
\left\{
\begin{array}{l}
	\ds\pt h_i \, + \, v \frac{\partial h_i}{\partial x} \,= \, \frac{\nu}{\eps} \, \left(\mathcal{M}[h]_i - h_i \right), 
	\\
	\,
	\\
	h_i(t=0)  \,= \,  h_{i,0}, \quad  i = 1, \, 2.
\end{array}
\right.
\eeq
The vector $(\mathcal{M}[h]_1, \mathcal{M}[h]_2)$ has the form
\[
(\mathcal{M}[h]_1, \mathcal{M}[h]_2) \, = \, (\mathcal{M}[h]_1, T\mathcal{M}[h]_1),
\]
where
$$
%\beq
\mathcal{M}[h]_1(v) \, := \, \int_{\mathbb{R}} \mathcal{M}[f](v,v_2, v_3) \,dv_2 dv_3 \,=\, \frac{\rho}{\sqrt{2\pi T}}\exp\left(-\frac{|v-u_1|^2}{2T}\right).
%\label{maxw1d}
%\eeq
$$
The macroscopic moments $\rho$, $u_1$ and $T$ of $f$  can be given in terms of $\bfh$
$$
%\beq
%\label{U1d}
\rho \,=\, \int_{\mathbb{R}}h_1(v) \, d v, \quad u_1 \,=\, \frac{1}{\rho} \int_{\mathbb{R}} v \, h_1(v) \, d v, \quad 3T  \,=\, \frac{1}{\rho}\left[\int_{\mathbb{R}}\frac{|v-u_1|^2}{2} \, h_1(v) \, d v \,+\, \int_{\mathbb{R}} h_2(v) \, dv\right].
%\eeq
$$
\subsection*{Flux relation}
\label{sec3.2}
Now we establish the relation between the hydrodynamic equations \eqref{eq:cns1} and the reduced BGK systems \eqref{bgk1d} in 1D and \eqref{bgk2d} in 2D in space respectively.  For the hybrid discontinuous Galerkin scheme proposed in the following section, due to its compactness, we only need to define consistent numerical fluxes at the interface of two cells in the hydrodynamic region and the kinetic region respectively. Here we would like to borrow the form \eqref{eq:kflux} and provide explicit formulae for the flux functions in \eqref{eq:cns1} to the unknowns of the reduced BGK systems \eqref{bgk1d} and \eqref{bgk2d}.

First for the 1D BGK system, the deformation tensor is simply
\[
\bfD(\bfu) \, = \, {\rm Diag}\left(\frac43 \px u_{1}, -\frac23 \px
  u_{1},-\frac23 \px u_{1}\right).
\]
The first order truncated distribution function $f^\eps(\bfv)$ \eqref{eq:feps} becomes 
\beq
f_T(\bfv)\,\,=\,\,\mathcal{M}(\bfv
)\,\left[1-\frac{\eps}{\nu}\left(\big(V_1^2-\frac13|\bfV|^2\big) \px
    u_{1}+  \frac12\big(|\bfV|^2-5\big)\,V_1 \frac{\px T}{\sqrt{T}}\right)\right],
\label{eq:fcut1d}
\eeq
where $\bfV = (V_1, V_2, V_3):=\frac{\bfv-\bfu}{\sqrt{T}}$. Corresponding to the reduced BGK system \eqref{bgk1d} in 1D, we have
\beq
\label{eq:FT}
\begin{cases}
	&\ds h_{1,T}(v_1) \,=\,\int_{\mathbb{R}^2} f_T(\bfv) \,dv_2
        dv_3\,=\,\mathcal{M}[h]_1\left[1-\frac{\eps}{\nu}\frac23(V_1^2-1)\,
          \px u_{1} \,-\, \frac{\eps}{\nu} \frac12V_1(V_1^2-3)\frac{\px T}{\sqrt{T}}\right], \\ \, \\
	&\ds h_{2,T}(v_1) \,=\,\int_{\mathbb{R}^2} f_T(\bfv)
        \frac{v_2^2+v_3^2}{2} \;dv_2dv_3 \,=\, T \left[\,
          h_{1,T}(v_1)+\frac{\eps}{\nu} \mathcal{M}[h]_1 \left(\frac23
            \px u_{1} - V_1\frac{\px T}{\sqrt{T}}\right) \right].
\end{cases}
\eeq
The subindex ``T" denotes truncations for the corresponding unknown functions. The nonzero first component of the flux function \eqref{eq:kflux} becomes
\beq
\label{eq:FFLUX}
\bfF_1(\bfU,\px \bfU) \,=\, \int_{\mathbb{R}^{3}} v_1 \, m(\bfv)  \,
f_T(\bfv) \, d\bfv \,=\,  \int_{\mathbb{R}} v_1 \, \left[ \,m(v_1)  \, h_{1,T}(v_1) \, +\, {\bf e}_3 \, h_{2,T}(v_1) \right] \, dv_1,
\eeq
which will be used to define the numerical hydrodynamic flux at the interface of two cells between two regions. Here $m(v_1)=(1,v_1,v_1^2/2)^T$ and ${\bf e}_3=(0,0,1)^T$.

For the 2D reduced BGK system \eqref{bgk2d}, the truncated distribution function $f_T(\bfv)$ is
\begin{align}
\label{eq:fcut2d} 
f_T(\bfv)\,\,= \,\, \mathcal{M} \,\Bigg\{ &1-\frac{\eps}{\nu}\left[
                                            \left(V_1^2-\frac13|\bfV|^2\right)
                                            \px
                                            u_{1}+\left(V_2^2-\frac13|\bfV|^2\right)
                                            \py u_{2} \right.
\\ & \left. + \, V_1V_2(\px u_{2} + \py u_{1}) +
     \frac12\big(|\bfV|^2-5\big) \left(V_1 \frac{\px
     T}{\sqrt{T}}+V_2\frac{\py T}{\sqrt{T}} \right) \right] \Bigg\}, \notag
\end{align}
so that
$$
%\beq
%\label{eq:FT2d}
\left\{
\begin{array}{lll}
\ds	g_{1,T}(v_1,v_2) &= &\ds\int_{\mathbb{R}} f_T(\bfv) dv_3   \\ \, \\
\, &= &\ds \mathcal{M}[g]_1 \, \Bigg\{ 1-\frac{\eps}{\nu}\left[
        \left(V_1^2-\frac13(V_1^2+V_2^2+1)\right)\,\px u_{1}+\left(V_2^2-\frac13(V_1^2+V_2^2+1)\right)\,\py u_{2}\right.  \\ \, \\
\,&+&\ds\left. V_1\,V_2\,\left(\px u_{2} + \py u_{1}\right) \,+\,
      \frac12\left(V_1^2+V_2^2-4\right) \left(V_1 \frac{\px
      T}{\sqrt{T}}+V_2\frac{\py T}{\sqrt{T}} \right) \right] \Bigg\},  
\\ \, \\
\ds	g_{2,T}(v_1,v_2) &= & \ds\int_{\mathbb{R}} f_T(\bfv) \frac{v_3^2}{2} dv_3  \\ \, \\
	\, & = & \ds\frac12\,T \,g_{1,T}(v_1,v_2)+\frac{\eps \,
                 T}{\nu} \mathcal{M}[g]_1 \left(\frac13\left(\px u_{1}
                 + \py u_{2}\right) - \frac12\left(V_1\frac{\px
                 T}{\sqrt{T}}+V_2\frac{\py T}{\sqrt{T}}\right)\right).
\end{array}
\right.
%\eeq
$$
The flux function \eqref{eq:kflux} in 2D becomes
\beq
\label{eq:FFLUX2d}	
\bfF(\bfU,\nablax \bfU) \,=\, \int_{\mathbb{R}^{2}} \bfv \left[ \, m(\tilde\bfv) \, g_{1,T}(\tilde\bfv) \, + \,{\bf e}_3 \, g_{2,T}(\tilde\bfv) \,\right] \, d\tilde\bfv,
\eeq
in this equation we denote $\tilde\bfv=(v_1,v_2)$ and $m(\tilde\bfv)=(1,v_1,v_2,(v_1^2+v_2^2)/2)^T$.

\subsection*{Discontinuous Galerkin scheme}
The discontinuous Galerkin scheme
\eqref{eq2d:DGf:a}-\eqref{eq2d:DGf:c} can be easily adapted to the
reduced BGK system \eqref{bgk2d}, which is for any $\zeta \in \mZ_h^\bfK $ 
\begin{eqnarray*}	
\int_{I_{i,j}} \bfR^{n+1}_h(\bfv) \zeta(\bfx) \, d\bfx &= &
                                                            \int_{I_{i,j}}
                                                            \bfg^n_h(\bfv)
                                                            \zeta(\bfx)
                                                            d\bfx  \,+\,
                                                            \Delta t
                                                                                                                       \int_{I_{i,j}}
                                                            \bfv \cdot
                                                            \nablax
                                                            \zeta(\bfx)
                                                            \,
                                                            \bfg^n_h(\bfv)
                                                            \, d\bfx
 \\   
& -& \Delta t\,\int_{I_i} \widetilde{(v_1 \bfg)}(x_{\iR},y)\zeta(x^-_\iR, y) - \widetilde{(v_1 \bfg)}(x_\iL,y) \zeta(x^+_\iL, y) \, dy  
\\
 & -&  \Delta t\,\int_{I_j} \widetilde{(v_2 \bfg)}(x,y_{\jR})\zeta(x, y^-_\jR) - \widetilde{(v_2 \bfg)}(x,y_\jL) \zeta(x, y^+_\jL) \, dx,
\end{eqnarray*}
where $\widetilde{v\, \bfg}$ is an upwind numerical flux. Then,  $\bfU^{n+1}_h$ is given for any $\beta \in \mZ_h^\bfK $ 
\begin{eqnarray*}	
	\int_{I_{i,j}}\bfU^{n+1}_h \beta(\bfx) \, d\bfx &= & 
                                                                \int_{I_{i,j}}
                                                                \int_{\mathbb{R}^2}
                                                                {\bf
                                                                \Phi}(\bfv)
                                                                \,\bfR^{n+1}_h(\bfv)
                                                                d\bfv\,                                              
                                                                \beta(\bfx)\,
                                                                \,
                                                                d\bfx, 
\end{eqnarray*}
 where the matrix ${\bf\Phi}$ is given by   
\[
{\bf \Phi}(\bfv) \, :=  \,
\begin{pmatrix}
1      & 0 \\
\bfv  & 0 \\
{|\bfv|^2}/{2} & 1
\end{pmatrix}.
\]
Finally, for any $\alpha \in \mZ_h^\bfK $  we have
\begin{eqnarray*}	
\int_{I_{i,j}} \bfg^{n+1}_h(\bfv) \alpha(\bfx) \, d\bfx & = & \,
                                                              \int_{I_{i,j}}\frac{\eps}{\eps+\nu^{n+1}\Delta t}\bfR^{n+1}_h(\bfv)\alpha(\bfx)
                                                              d\bfx  +
                                                              \int_{I_{i,j}}
                                                              \frac{\nu^{n+1}\Delta
                                                              t}{\eps+\nu^{n+1}\Delta
                                                              t}
                                                              \bfM[\bfg](\bfv,
                                                              \bfU^{n+1}_h)
                                                              \,
                                                              \alpha(\bfx)
                                                              \,
                                                              d\bfx,  
\end{eqnarray*}
with $\bfg(\bfv) := (g_1, g_2)^\text{T}$ and $\bfM[\bfg] :=
(\mathcal{M}[g]_1, \mathcal{M}[g]_2 )^\text{T}$, whereas $\bfR(\bfv) =  \bfg(\bfv) - \Delta t \, \bfv \cdot \nablax \,
\bfg(\bfv)$.

\bibliographystyle{plain}
\bibliography{refer}
\end{document}